%% file: CPClassificationCubicActivation_arXiv_2026_07_16.tex
\documentclass[12pt]{article}

\input{commands.tex}

\begin{document}

\title{Landscape analysis for shallow neural networks: \\
Complete classification of critical points for \\
cubic activation and affine target functions}

\author{Shokhrukh Ibragimov$^{1}$, Ilkhom Mukhammadiev$^{2}$, and Diyora Salimova$^{3}$
\bigskip\\
\small{$^1$ Applied Mathematics: Institute for Analysis and Numerics,}\vspace{-0.1cm}\\
\small{University of M\"unster, Germany; e-mail: \texttt{sibragim}\textcircled{\texttt{a}}\texttt{uni-muenster.de}}\smallskip\\
\small{$^2$ Department for Applied Mathematics, University of Freiburg,}\vspace{-0.1cm}\\
\small{Germany; e-mail: \texttt{ilkhom.mukhammadiev}\textcircled{\texttt{a}}\texttt{mathematik.uni-freiburg.de}}\smallskip\\
\small{$^3$ Department for Applied Mathematics, University of Freiburg,}\vspace{-0.1cm}\\
\small{Germany; e-mail: \texttt{diyora.salimova}\textcircled{\texttt{a}}\texttt{mathematik.uni-freiburg.de}}
}

\date{\today}

\maketitle

\begin{abstract}
In this paper, we study the optimization landscape induced by the true loss for shallow polynomial neural networks (PNNs) with $\width \in \N$ neurons on the hidden layer, one-dimensional input and output layers, and a monomial activation of degree $d \in \N$, trained against a non-constant affine linear target function. Our first main result provides for arbitrary activation degree $d$ a sharp existence/non-existence criterion for \emph{global minimizers} with necessary structural conditions. We show that the infimum of the loss is always zero and achievable with at least $d$ \Active\ and \visible\ hidden neurons -- that is, hidden neurons with non-zero inner and outer weights -- with pairwise distinct \pivot s. In contrast, if $\width < d$, then the infimum cannot be attained and any minimizing sequence of parameters necessarily diverges to infinity. In the second main result, we provide a complete classification of all critical points of the loss function for the cubic activation. We show that the loss landscape admits no \emph{local maximizers}, critical points cannot have exactly two distinct \pivot s, global minimizers require at least three distinct \pivot s, critical points with no \Active\ hidden neurons correspond to \emph{saddle points} only, and consequently, \emph{non-global local minimizers} and non-trivial saddle points arise only in networks where all \pivot s coincide. Moreover, non-global local minimizers require all hidden neurons to be \Active\ and \visible\ with exactly one hidden neuron having a \slope\ sign matching that of the target function. Our second main result also guarantees that each hidden neuron of a critical point that is not a global minimizer has either input-dependent or zero contribution, but has no nonzero input-independent contribution, to its corresponding realization function.
\end{abstract}

\newpage
\tableofcontents
\section{Introduction}

\ANNs[\emph{Artificial neural networks}]\ have become a central tool in modern data-driven modeling, owing to their remarkable empirical performance in approximation, prediction, and representation learning across a wide range of applications. From a mathematical perspective, a defining feature of these models is the nonlinear activation function applied at each hidden neuron. While this nonlinearity is crucial for learning complex representations, it fundamentally breaks the convexity of the underlying training problem, leading to highly intricate loss landscapes. Consequently, understanding how the choice of activation function interacts with network architecture to shape the geometry of the optimization landscape remains a fundamental problem at the interface of approximation theory and optimization.

Although classical results \cite{MR1015670,HORNIK1989359,HORNIK1991251,LESHNO1993861} dictate that feedforward networks with locally bo\-un\-ded and piecewise continuous activations are universal approximators if and only if the activation is non-polynomial, this limitation in expressivity is offset by the profound theoretical utility of monomial activations in optimization and specific learning tasks (cf., e.g., \cite{BALDI198953,Kawaguchi_no_loc_min,MannelliQuadraticActivation,DuSimonOverparametrization,MR4002887,TragerShallowPNNs,MR4873488} and the references therein). While modern \DL\ predominantly relies on piecewise-polynomial activations such as the \ReLU, smooth activations such as monomials offer an analytically tractable setting where the loss function reduces to a polynomial in the network parameters. This algebraic structure enables a precise, global investigation of critical points and loss landscape degeneracies -- a notoriously difficult task for non-smooth models. \PNNs[\emph{Polynomial neural networks}]\ -- that is, \ANNs\ with monomial activations -- therefore serve as an ideal testbed for developing a rigorous understanding of optimization phenomena arising in the training of \ANNs.

The structure of the loss landscape is of fundamental importance for gradient-based optimization methods, which constitute the standard paradigm for training \ANNs\ (see, e.g., \cite{MR4505886,LiVisualizingLossLandscape} and the references therein). The convergence behavior of first-order methods such as \GD\ and its variants is governed by the geometry of critical points: \textcolor{magenta}{\emph{global minima}} -- the smallest possible loss value -- correspond to exact representability of the target function in case of zero infimum, \textcolor{magenta}{\emph{non-global local minima}} may obstruct convergence to optimal solutions (cf., e.g., \cite{IbragimovNonGlobalLocMin,SafranShamirLocMin,MR4813240,PascanuLocalMinima,RiekertGlobalMinimaNonConvergence,IbragimovExistenceInfLocMin2022}), \textcolor{magenta}{\emph{saddle points}} can substantially slow down optimization due to flat or nearly flat directions with small gradients (see, e.g., \cite{MR4837253,MR3960812,LeeEscapeSaddles,PanageasEscapeSaddles2016,PanageasEscapeSaddles2019,DauphinAttackSaddles} for avoiding convergence of \GD\ type optimization methods to saddle points), and the presence or absence of \textcolor{magenta}{\emph{local maxima}} reflects fundamental curvature properties of the objective; see, e.g., \cite{IbragimovNonGlobalLocMin,ThangDoNonConvergenceGlobalMinima,RiekertGlobalMinimaNonConvergence} for non-convergence and \cite{BachGlobalMinimaConvergence,MR4205063,NIPS2011_40008b9a,MR3839649,Nesterov2014UniversalGM,MR4716376} for convergence of \GD\ type optimization methods to global minima depending, e.g., on the appearance of non-optimal critical points. A detailed structural characterization of critical points therefore provides a rigorous insight into the optimization geometry and clarifies how architectural parameters influence the existence and stability of suboptimal stationary points.

Recently, there have been many results on the investigation of the geometry of loss landscapes arising in empirical risk minimization for \ANNs. For linear networks under squared loss, classical results \cite{BALDI198953,Kawaguchi_no_loc_min} show that every local minimum is global and all other critical points are saddles. Beyond the linear setting, \cite{FUKUMIZU2000317} derived partial structural characterizations by exploiting network hierarchy, without imposing restrictive assumptions on the target function or activation. For nonlinear activations, however, the landscape becomes more complex: \cite{SafranShamirLocMin} established the existence of spurious local minima for \ReLU\ networks using computer-assisted arguments, while \cite{MannelliQuadraticActivation,MR4873488} showed that, for quadratic activation under squared loss,  the geometry differs substantially from the linear case; nevertheless, non-global local minima may disappear under overparameterization (cf., e.g., \cite{DuSimonOverparametrization,MR4002887}). For deep \ANNs, non-global local minima vanish with high probability by increasing the width of the last hidden layer (cf., e.g., \cite{LivniNNs2014,MR3904911,SoudryNoBadLocMin}). Additional theoretical studies of nonlinear landscapes include \cite{ChoromanskaLossSurface,SoudrySubOptimalLocMin}, complemented by numerical investigations in \cite{DauphinAttackSaddles}, although these analyses rely on certain simplifying assumptions (cf.\ \cite{pmlr-v40-Choromanska15}). These results indicate that the landscape properties known for linear networks do not generally extend to nonlinear activations.

While the preceding works focus on empirical risk, the population loss isolates the architecture's intrinsic optimization geometry from sampling effects, thereby providing a natural framework for theoretical analysis. Within this setting, \cite{MR4450126} classified the critical points for shallow \ANNs\ with \ReLU, leaky \ReLU, and quadratic activations under squared $L^2$ population loss and non-constant affine target functions, showing, in particular, that no non-global local minima occur for leaky \ReLU\ and quadratic activations.

In this work, we study the optimization landscape induced by the true squared $L^2$ loss for shallow \PNNs\ with a monomial activation of degree $d \in \N$, trained against a non-constant affine linear target function. We show, in particular, that the loss function admits no local maximizers due to its structural properties and provide a sharp existence/non-existence criterion for global minimizers in terms of the relation between the activation degree and the network width. Moreover, in the case of cubic activation, we establish a full classification of all critical points in structural terms, describing global minimizers, non-global local minimizers, and saddle points via the activity patterns of neurons and their associated realization functions.

In the following, we describe our main findings in detail in \cref{thm:intro:existence_of_global_minima} and \cref{thm:intro:complete_classification_cubic} below by employing in part the notations and terminologies used throughout the paper. The complete list of notations and terminologies employed in the paper is introduced  later in detail  in \cref{subsec:notations}.

\subsection{First main result: Existence of global minima for monomial ac\-ti\-va\-tions}
\label{subsec:FirstMainResult}

Consider shallow \PNNs\ with $\width \in \N = \{1, 2, \dots\}$ neurons on the hidden layer, total number
\begin{equation}
\textstyle \dimension = \width + \width + \width + 1 = 3 \width + 1
\end{equation}
of parameters consisting of $\width$ inner weights, $\width$ inner biases, $\width$ outer weights, and $1$ outer bias (cf.\ \cref{subsec:notations}), monomial ac\-ti\-va\-tion function
\begin{equation}
\textstyle \R \ni x \mapsto \cA(x) = x^d \in \R
\end{equation}
of degree $d \in \N$, non-constant affine linear target function $f \colon \R \to \R$, and squared $L^2$ loss function $\loss \colon \R^{\dimension} \to \R$ (see \cref{eqn:thm:intro:existence_of_global_minima:loss} below).

\cref{thm:intro:existence_of_global_minima} below then shows that the infimum of the loss function $\loss$ (which is a polynomial) is always zero (see \cref{item1:thm:intro:existence_of_global_minima}) and reveals a sharp dichotomy between the \emph{narrow regime}, where the width $\width$ of the hidden layer of shallow \PNN\ under consideration is strictly smaller than the activation degree $d$, and the \emph{wide regime}, where the width $\width$ is at least equal to the activation degree $d$ (see \cref{item2:thm:intro:existence_of_global_minima,item3:thm:intro:existence_of_global_minima,item4:thm:intro:existence_of_global_minima}). In particular, it shows   in the narrow regime that  non-constant affine linear target functions cannot be represented exactly, global minima fail to exist, and any sequence of parameters along which the loss function vanishes necessarily diverges. In contrast, in the wide regime it demonstrates that exact representability becomes possible and the existence of global minimizers with the necessary \Active-\visible\ \pivot\ structure can be characterized in terms of the number and configuration of \Active\ hidden neurons. Furthermore, it proves that the existence of global minima requires at least as many \Active\ and \visible\ hidden neurons -- that is, hidden neurons with non-zero inner and outer weights --  as the activation degree with pairwise distinct \pivot s (see \cref{item3:thm:intro:existence_of_global_minima}). We now state \cref{thm:intro:existence_of_global_minima}.

\cfclear
\begin{savenotes}
\begin{samepage}
\begin{tcolorbox}[colback=white!95!gray,
                  colframe=black,
                  boxrule=0.5pt,
                  sharp corners,
                  enhanced,
                 ]
\begin{theorem}\label{thm:intro:existence_of_global_minima}
Let $d, \dimension, \width \in \N$, $\scra \in \R$, $\scrb \in (\scra, \infty)$ satisfy $\dimension = 3 \width + 1$, let $f \colon \R \to \R$ be affine linear with $f'(0) \neq 0$, let $\norm{\cdot} \colon  \R^{\dimension} \to \R$ be the Euclidean norm on $\R^{\dimension} $, and let $\loss \colon \R^{\dimension} \allowbreak \to \R$ satisfy for all $\theta = (\theta_1, \dots, \theta_{\dimension}) \allowbreak \in \R^{\dimension}$ that
\begin{equation}\label{eqn:thm:intro:existence_of_global_minima:loss}
\textstyle \loss(\theta) = \int_{\scra}^{\scrb} \Abs{\theta_{\dimension} \allowbreak + \sum_{i = 1}^{\width} \theta_{2 \width + i} (\theta_{\width + i} + \theta_i x)^d - f(x)}^2 \, \d x.
\end{equation}
Then
\begin{enumerate}[label=\textnormal{(\roman*)}]
\item
\label{item1:thm:intro:existence_of_global_minima} $\loss$ is a polynomial with $\inf_{\theta \in \R^{\dimension}} \loss(\theta) = 0$,

\item
\label{item2:thm:intro:existence_of_global_minima} it holds that $\{\theta \in \R^{\dimension} \colon \loss(\theta) = 0\} \neq \emptyset$ if and only if $\width \ge d$,

\item
\label{item3:thm:intro:existence_of_global_minima} it holds for all $\theta \in \R^{\dimension}$ with $\loss(\theta) = 0$ that
\begin{equation}
\textstyle \# \bigl(\bigcup_{i \in \{1, 2, \dots, \width\}, \theta_i \theta_{2 \width + i} \neq 0} \{\nicefrac{\theta_{\width + i}}{\theta_i}\}\bigr) - d \ge 0 = (\nabla \loss)(\theta),
\end{equation}
and
\item
\label{item4:thm:intro:existence_of_global_minima} if $\width < d$ it holds for all $(\phi_n)_{n \in \N} \subseteq \R^{\dimension}$ with $\limsup_{n \to \infty} \loss(\phi_n) = 0$ that
\begin{equation}
\textstyle \liminf_{n \to \infty} \norm{\phi_n} = \infty.
\end{equation}
\end{enumerate}
\end{theorem}
\end{tcolorbox}
\end{samepage}
\end{savenotes}

This result highlights how even in a  simplified and analytically tractable setting, the interplay between activation degree and network width gives rise to rich and nontrivial optimization landscapes. \cref{thm:intro:existence_of_global_minima} is a direct consequence of \cref{prop:existence_of_global_minima} in \cref{subsec:existence_of_global_minima} below. We also refer to \cref{lemma:target_fn_non_reachable} and \cref{lemma:target_fn_reachable} in \cref{subsec:polynomial_target_functions_representability} below for non-representability and representability of polynomial target functions of degree possibly greater than one, which are used in our proof of \cref{prop:existence_of_global_minima} (and hence \cref{thm:intro:existence_of_global_minima}).

Similar results have been obtained in the literature in case of commonly used non-polynomial activations. On the one hand,  in case of smooth activations such as the standard logistic (sigmoid), the softsign, the softplus, the inverse tangent (arctan), and the hyperbolic tangent activations it has been shown, e.g., in \cite{GallonBlowUpPhenomena,MR4243432}, that there, in general,  do not exist global minimizers in the optimization landscape even for smooth target functions, i.e., any sequence $(\phi_n)_{n \in \N} \subseteq \R^{\dimension}$ satisfying $\lim_{n \to \infty} \loss(\phi_n) = \inf_{\theta \in \R^{\dimension}} \loss(\theta)$ must diverge to infinity. On the other hand, it has been proved, e.g., in \cite{MR4832358,MR4716376} in the case of  the \ReLU\ activation and Lipschitz continuous target functions as well as in \cite{KAINEN2000695} in the case of the Heaviside activation and indicator target functions, that any shallow \ANN\ of arbitrary hidden layer width does admit global minimizers in the optimization landscape.

For results on the expressive power of \ANNs\ with monomial-type activations  related to \cref{item2:thm:intro:existence_of_global_minima} in \cref{thm:intro:existence_of_global_minima}, we refer, e.g., to \cite{KileelExpressivePowerPNNs,MR4836468,MR4990261,LiPowerNet} and the references therein.

\cref{thm:intro:existence_of_global_minima} above serves as a complementary result to our main result, \cref{thm:intro:complete_classification_cubic} below, which we now describe in detail in the following subsection.

\subsection{Second main result: Classification of critical points for the cubic ac\-ti\-va\-tion}
\label{subsec:SecondMainResult}

Consider shallow \PNNs\ with $\width \in \N = \{1, 2, \dots\}$ neurons on the hidden layer, total number
\begin{equation}
\textstyle \dimension = \width + \width + \width + 1 = 3 \width + 1
\end{equation}
of parameters consisting of $\width$ inner weights, $\width$ inner biases, $\width$ outer weights, and $1$ outer bias (cf.\ \cref{subsec:notations}), cubic ac\-ti\-va\-tion function
\begin{equation}
\textstyle \R \ni x \mapsto \cA(x) = x^3 \in \R,
\end{equation}
non-constant affine linear target function $f \colon \R \to \R$, realization functions $\cN^{\theta} \colon \R \to \R$, $\theta \in \R^{\dimension}$ (see \cref{eqn:thm:intro:complete_classification_cubic:cN} below), squared $L^2$ loss function $\loss \colon \R^{\dimension} \to \R$ over the integration domain $\emptyset \neq [\scra, \allowbreak \scrb] \allowbreak \subseteq \R$, arbitrary parameter vector $\vartheta = (\vartheta_1, \allowbreak \dots, \allowbreak \vartheta_{\dimension}) \allowbreak \in \R^{\dimension}$, and number $\eta \in \N_0$ of hidden neurons of $\vartheta$ either with a zero \slope\ or with the \slope\ sign matching that of the target function $f$, i.e.,
\begin{equation}
\textstyle \eta = \#\{i \in \{1, 2, \dots, \width\} \colon \allowbreak f'(\scra) \vartheta_i \vartheta_{2 \width + i} \ge 0\}.
\end{equation}

\cref{thm:intro:complete_classification_cubic} then proves that the loss function $\loss$ is a polynomial in the parameter vector without any local maximum points (see \cref{item:loc_max:thm:intro:complete_classification_cubic}) and that there are no critical points with exactly two distinct \pivot s of \Active\ hidden neurons (cf.\ \cref{subsec:notations}), i.e.,
\begin{equation}\label{eqn_intro:no_cp_2_pivot}
\textstyle \bigl\{\theta \in (\nabla \loss)^{-1}(\{0\}) \colon \# (\cup_{i \in \{1, 2, \dots, \width\}, \theta_i \neq 0} \{\nicefrac{\theta_{\width + i}}{\theta_i}\}) = 2\bigr\} = \emptyset
\end{equation}
(which follows from \cref{item:global_min:thm:intro:complete_classification_cubic}). In addition, it shows that $\vartheta$ is
\begin{itemize}
\item a global minimum point of $\loss$ if and only if it is a critical point of $\loss$, i.e., $(\nabla \loss)(\vartheta) = 0$, and has at least three \Active\ and \visible\ hidden neurons with pairwise distinct \pivot s,

\item a non-global local minimum point of $\loss$ if and only if all its hidden neurons are \Active\ and \visible\ with exactly one of them having a \slope\ sign matching that of the target function, i.e., $\eta = 1$, and admits the realization function given by
\begin{equation}\label{eqn:intro:locmin_realization}
\textstyle \Forall x \in \R \colon \cN^{\vartheta}(x) = f(\frac{\scrb + \scra}{2}) + \frac{28 f'(\scra)}{5 (\scrb - \scra)^2} (x - \frac{\scrb + \scra}{2})^3,
\end{equation}
and
\item a saddle point of $\loss$ if and only if each hidden neuron that is not \Active\ ($\vartheta_i = 0$) is \inactive\ ($\vartheta_{\width + i} = 0$) or \invisible\ ($\vartheta_{2 \width + i} = 0$) and in addition one of the following holds:

\begin{enumerate}[label=(\alph*)]
\item It has only \inactive\ or \invisible\ \semiactive\ hidden neurons with $\vartheta_{\dimension} = f(\frac{\scrb + \scra}{2})$.

\item It has all finite \pivot s coinciding with $\scrp_0$ in \cref{eqn:thm:intro:complete_classification_cubic:scrp_scrf} and the same realization function as the non-global local minimum points as in \cref{eqn:intro:locmin_realization} but, in contrast, with at least two hidden neurons having zero \slope\ or \slope\ sign matching that of the target function, i.e., $\eta > 1$.

\item Its realization function is one of the two side cubic critical realizations corresponding to $\scrf_1 = 0$ or $\scrf_{- 1} = 0$, and all \Active\ \pivot s coincide with the corresponding $\scrp_1$ and $\scrp_{- 1}$, respectively, in \cref{eqn:thm:intro:complete_classification_cubic:scrp_scrf}.
\end{enumerate}
\end{itemize}

\cref{thm:intro:complete_classification_cubic} \hence shows, in particular, that \PNNs\ with one hidden neuron only, i.e., when $\width = 1$, and with the realization function as in \cref{eqn:intro:locmin_realization}, can correspond among critical points only to non-global local minimum points. We now present the formal statement of \cref{thm:intro:complete_classification_cubic}.

\cfclear
\begin{savenotes}
\begin{samepage}
\begin{tcolorbox}[colback=white!95!gray,
                  colframe=black,
                  boxrule=0.5pt,
                  sharp corners,
                  enhanced,
                 ]
\begin{theorem}\label{thm:intro:complete_classification_cubic}
Let $\dimension, \width \in \N$, $\scra \in \R$, $\scrb \in (\scra, \infty)$, $\Width \subseteq \N$ satisfy $\dimension = 3 \width + 1$ and $\Width \allowbreak = \{1, \allowbreak 2, \allowbreak \dots, \allowbreak \width\}$, let $f \colon \R \to \R$ be affine linear with $f'(0) \neq 0$, for every $\theta = \allowbreak (\theta_1, \allowbreak \dots, \allowbreak \theta_{\dimension}) \allowbreak \in \R^{\dimension}$ let $\cI^{\theta} \subseteq \Width$ satisfy $\cI^{\theta} \allowbreak = \{i \allowbreak \in \Width \colon \allowbreak \theta_i \neq 0\}$ and let $\cN^{\theta} \colon \R \to \R$ satisfy for all $x \in \R$ that
\begin{equation}\label{eqn:thm:intro:complete_classification_cubic:cN}
\textstyle \cN^{\theta}(x) = \theta_{\dimension} + \sum_{i = 1}^{\width} \theta_{2 \width + i} (\theta_{\width + i} + \theta_i x)^3,
\end{equation}
let $\loss \colon \R^{\dimension} \allowbreak \to \R$ satisfy for all $\theta \allowbreak \in \R^{\dimension}$ that $\loss(\theta) = \int_{\scra}^{\scrb} \abs{\cN^{\theta}(x) - f(x)}^2 \, \d x$, let $\vartheta \allowbreak = (\vartheta_1, \allowbreak \dots, \allowbreak \vartheta_{\dimension}) \allowbreak \in \R^{\dimension}$, and for every $j \in \{-1, 0, 1\}$ let $\scrp_j, \allowbreak \scrf_j \allowbreak \in \R$, $\eta \in \N_0$ satisfy
\begin{equation}\label{eqn:thm:intro:complete_classification_cubic:scrp_scrf}
\textstyle \scrp_j = \frac{\scrb + \scra}{2} + \frac{j (\scrb - \scra)}{2 \sqrt{7}}, \quad \scrf_j = \sup_{x \in [\scra, \scrb]} \! \Abs{f(\frac{\scrb + \scra}{2}) + \frac{2 j (\scrb - \scra) f'(\scra)}{5 \sqrt{7}} + \frac{28 f'(\scra) (x - \scrp_j)^3}{5 (1 + \abs{j}) (\scrb - \scra)^2} - \cN^{\vartheta}(x)},
\end{equation}
and $\eta = \#\{i \in \Width \colon \allowbreak f'(\scra) \vartheta_i \vartheta_{2 \width + i} \ge 0\}$. Then the following statements hold:
\begin{enumerate}[label=\textnormal{(\roman*)}]
\item
\label{item:loc_max:thm:intro:complete_classification_cubic} $\loss$ is a polynomial with \textcolor{magenta}{no local maximum points}.

\item
\label{item:global_min:thm:intro:complete_classification_cubic} The following statements are equivalent:
\begin{enumerate}[label=\textnormal{(\alph*)}]
\item $\vartheta$ is a \textcolor{magenta}{global minimum point} of $\loss$ with $\loss(\vartheta) = 0$.

\item It holds that $\# (\cup_{i \in \cI^{\vartheta}} \{\nicefrac{\vartheta_{\width + i}}{\vartheta_i}\}) - 1 > 0 = (\nabla \loss)(\vartheta)$.

\item It holds that $\# (\cup_{i \in \cI^{\vartheta}, \vartheta_{2 \width + i} \neq 0} \{\nicefrac{\vartheta_{\width + i}}{\vartheta_i}\}) - 3 \ge 0 = (\nabla \loss)(\vartheta)$.
\end{enumerate}

\item
\label{item:loc_min:thm:intro:complete_classification_cubic} $\vartheta$ is a \textcolor{magenta}{non-global local minimum point} of $\loss$ if and only if $\eta - 1 = 0 = \scrf_0$.

\item
\label{item:saddle:thm:intro:complete_classification_cubic} $\vartheta$ is a \textcolor{magenta}{saddle point} of $\loss$ if and only if $\sum_{i = 1}^{\width} (\abs{\vartheta_i} + \abs{\vartheta_{\width + i} \vartheta_{2 \width + i}}) = 0 = f(\frac{\scrb + \scra}{2}) - \vartheta_{\dimension}$ or
\begin{equation}
\prod_{j = -1}^1 \textstyle \Bigl[\indicator{\{0\}}(j) \indicator{\{0\}}(\eta - 1) + \scrf_j + \sum\limits_{i \in \cI^{\vartheta}} \Abs{\vartheta_{\width + i} [\vartheta_i]^{- 1} + \scrp_j}\Bigr] = 0 = \!\displaystyle \sum_{i \in \Width \backslash \cI^{\vartheta}} \abs{\vartheta_{2 \width + i} \vartheta_{\width + i}}.
\end{equation}
\end{enumerate}
\end{theorem}
\end{tcolorbox}
\end{samepage}
\end{savenotes}

\cref{thm:intro:complete_classification_cubic} is a direct consequence of \cref{thm:complete_classification_cubic} in \cref{subsec:complete_classification} below. \cref{thm:complete_classification_cubic} is, in turn, established by combining the structural characterization results presented separately in \cref{prop:global_minima_main}, \cref{prop:saddle_points}, and \cref{prop:local_minima} in \cref{subsec:structural_characterization} below with the non-existence of local maxima (cf.\ \cref{lemma:Hessian_of_Risk} in \cref{subsec:derivatives_of_loss} below).

One noteworthy consequence of \cref{thm:intro:complete_classification_cubic} (and \cref{thm:complete_classification_cubic}) is that for an arbitrary parameter $\theta \in \R^{\dimension}$, in order to be a non-global local minimum point of $\loss$, it is sufficient to require having only \Active\ and \visible\ hidden neurons with exactly one of them admitting a \slope\ sign matching that of the target function, i.e., $\eta = 1$, and a realization function as in \cref{eqn:intro:locmin_realization} above without additionally requiring to be a critical point of $\loss$ with all \pivot s coinciding with $\scrp_0$ in \cref{eqn:thm:intro:complete_classification_cubic:scrp_scrf} which will, in fact, be ensured by the former assumptions. This fact can be seen in \cref{lemma:loc_min} in \cref{subsec:local_minima} below which is used to derive the proof of \cref{thm:intro:complete_classification_cubic} (and \cref{thm:complete_classification_cubic}). In addition, \cref{thm:intro:complete_classification_cubic} guarantees that if a critical point $\theta \in \R^{\dimension}$ is not a global minimum point, then each of its hidden neurons contributes either input-dependent data ($\theta_i \theta_{2 \width + i} \neq 0$) or zero data ($\theta_{2 \width + i} [\abs{\theta_i} + \abs{\theta_{\width + i}}] = 0$), but not input-independent nonzero data ($\theta_i = 0 \neq \theta_{\width + i} \theta_{2 \width + i}$) to its realization function $\cN^{\theta}$.

\begin{figure}[!htb]
\centering
\subfloat{\includegraphics[width=10.0cm]{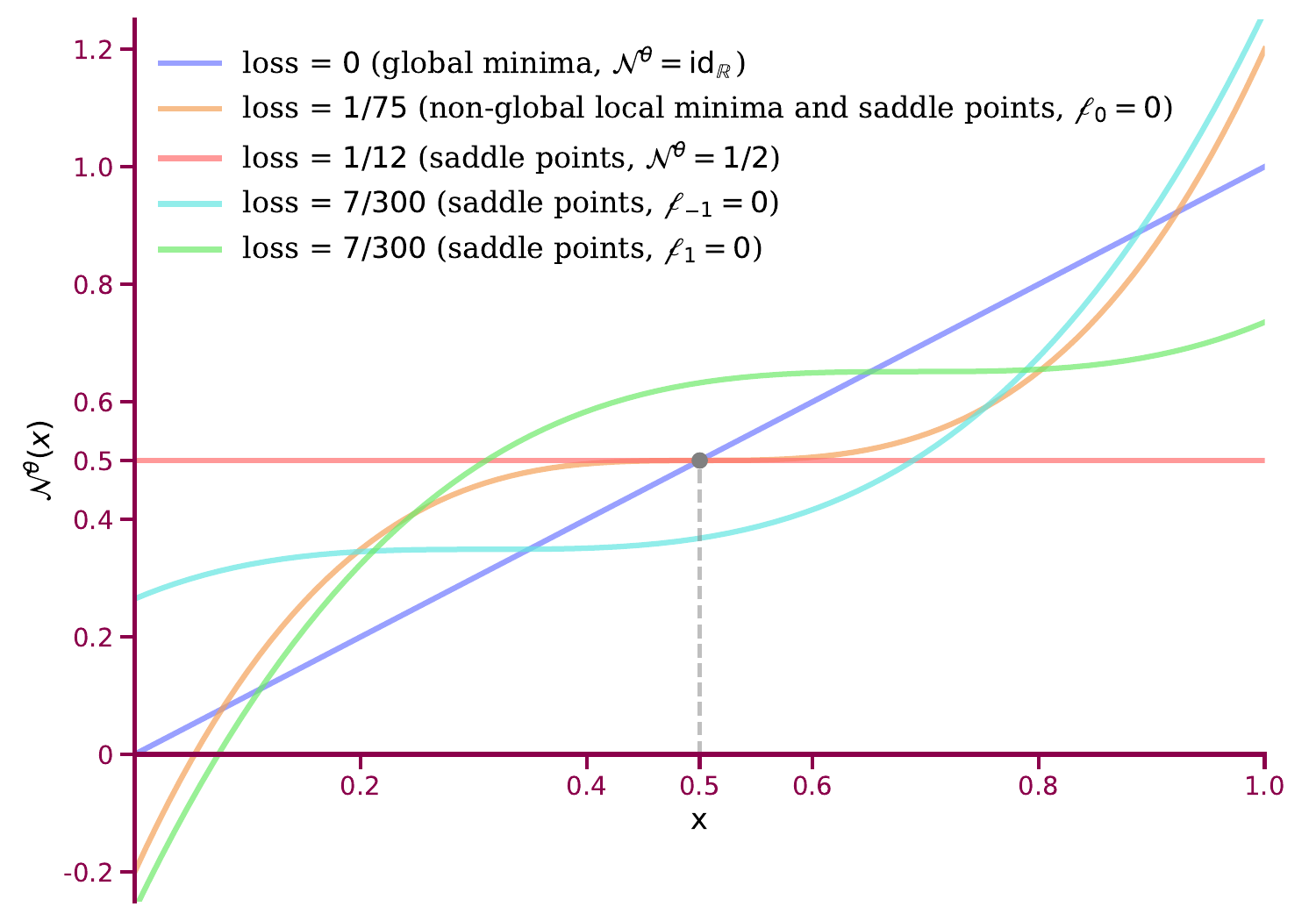}}
\caption{Numerical illustrations associated with \cref{thm:intro:complete_classification_cubic} in the canonical setting, where we consider a unit integration domain for the loss function $\loss$ and an identity target function: $\scra = 0$, $\scrb = 1$, and $f = \id_{\R}$. We plot the realization functions of critical points of the loss function with their corresponding loss values: The loss value $\nicefrac{1}{75}$ of non-global local minimizers corresponds to $\scrf_0 = 0$ and $\nicefrac{7}{300}$ of saddle points corresponds to $\scrf_1 = 0$ and $\scrf_{-1} = 0$ in \cref{eqn:thm:intro:complete_classification_cubic:scrp_scrf} while the loss values $0$ and $\nicefrac{1}{12}$ correspond to global minimizers, whose realization function is $\id_{\R}$, and saddle points with a constant realization function, $\nicefrac{1}{2}$, respectively. The loss value $\nicefrac{1}{75}$ corresponds to saddle points alongside non-global local minimizers, but only when $\width > 1$.}
\label{fig:critical_points}
\end{figure}

We \nobs that while our complete classification of critical points is currently established for the cubic case, we view this result as a fundamental stepping stone rather than an isolated phenomenon. Importantly, the techniques developed for this classification lay the essential mathematical groundwork for the broader setting. These techniques seem to provide a route toward a generalized characterization for \PNNs\ with arbitrary odd and even degree monomial activations.

\subsection{Structure of the article}

The remainder of this article is organized as follows. In \cref{sec:loss_derivatives}, for an arbitrary activation degree $d \in \N$ and an arbitrary continuous target function $f \in C(\R, \R)$, we exploit that the loss function $\loss \colon \R^{\dimension} \to \R$ is a polynomial in the \PNN\ parameter vector and derive in \cref{lemma:derivatives_of_Risk} and \cref{lemma:Hessian_of_Risk} explicit formulas for its gradient and the entries of its Hessian matrix. These results, in particular, show that the loss function is strictly convex with respect to the outer bias of an arbitrary input \PNN\ parameter vector, implying the absence of local maximizers.

In \cref{sec:global_minima}, for an arbitrary activation degree $d \in \N$ and a non-constant affine linear target function $f \in \cP(1)$, we establish in \cref{prop:existence_of_global_minima} that the zero infimum of the loss function $\loss$ is attained only in the wide regime -- that is, when the network width $\width \in \N$ satisfies $\width \ge d$ -- whereas it cannot be attained in the narrow regime $\width < d$, in which case every minimizing sequence necessarily diverges in the parameter space. We note again that \cref{thm:intro:existence_of_global_minima} above is an immediate consequence of \cref{prop:existence_of_global_minima}.

In \cref{sec:determination_cubic_activation}, we restrict to the cubic activation and a canonical setting with the identity target function $f = \id_{\R}$ and the unit integration domain $[0, 1]$ for the loss function $\loss$. Under this setting we determine in \cref{lemma:realizations} all critical points of the loss function -- zeros of the gradient $\nabla \loss$ -- in terms of explicit expressions for the realization functions and the configuration of \pivot s of \Active\ hidden neurons.

Finally, in \cref{sec:classification_cubic_activation}, we combine the results from the previous sections to establish in \cref{thm:complete_classification_cubic} a complete classification of all critical points of the loss function for the cubic activation. We first characterize critical points in the canonical setting with the identity target function $f = \id_{\R}$ and the unit integration domain $[0, 1]$ for the loss function $\loss$. Then we employ the invariance property of the loss landscape under suitable affine transformations, established in \cref{lemma:canonical_extension}, to extend these characterizations, in \cref{prop:global_minima_main}, \cref{prop:saddle_points}, and \cref{prop:local_minima}, to the general setting with an arbitrary domain of integration $\emptyset \neq [\scra, \scrb] \subseteq \R$ and an arbitrary non-constant affine linear target function $f \in \cP(1)$. \cref{thm:complete_classification_cubic} is simply a collection of the results in \cref{lemma:Hessian_of_Risk}, \cref{prop:global_minima_main}, \cref{prop:saddle_points}, and \cref{prop:local_minima}. We note again that \cref{thm:intro:complete_classification_cubic} above is an immediate consequence of \cref{thm:complete_classification_cubic}.

We remark that some of the symbolic computations in \cref{sec:determination_cubic_activation,sec:classification_cubic_activation} are verified using {\sc Wolfram Mathematica} \cite{Mathematica}.

\subsection{Notations employed throughout the article}
\label{subsec:notations}

Throughout the paper we make use of the following notations:
\begin{itemize}[itemsep=0pt]
\item We denote by $\norm{\cdot} \colon ( \cup_{n \in \N} \R^n) \to \R$ the Euclidean norm on $(\cup_{n \in \N} \R^n)$ defined for every $n \in \N$, $x = (x_1, \dots, x_n) \in \R^n$ by $\norm{x} = (\sum_{i = 1}^n \abs{x_i}^2)^{\nicefrac{1}{2}}$.

\item Let $\dimension \in \N$, let $U \subseteq \R^{\dimension}$ be open, let $f \colon U \to \R$ be a function, and let $x \in U$. Then we say that $x$ is a \emph{local minimum point} (or \emph{local minimizer}) of $f$ if and only if there exists $\eps \in (0, \infty)$ such that $f(x) = \inf_{y \in \{v \in U \colon \norm{v - x} \le \eps\}} f(y)$, we say that $x$ is a \emph{local maximum point} (or \emph{local maximizer}) of $f$ if and only if $x$ is a local minimum point of $(- f)$, and if $f$ is differentiable at $x$ then we say that $x$ is a \emph{saddle point} of $f$ if and only if $(\nabla f)(x) = 0$ and $x$ is neither a local minimum nor a local maximum point of $f$.

\item For every set $A \subseteq \R$ we denote by $\#A$ the number of elements of $A$.

\item Throughout the paper we use the convention $\max(\emptyset) = - \infty = - \min(\emptyset)$.

\item We \nobs that in the binomial expansions (see, e.g., \cref{eqn:lemma:non_representability_poly:binomial} in \cref{subsec:polynomials_representability_via_monomials} below) we adopt the convention $0^0 = 1$.

\item We call every constant polynomial a polynomial of degree zero and for every $\scrd \in \N$ we call $\R \ni x \mapsto p(x) \in \R$ a polynomial of degree $\scrd$ if and only if it satisfies for all $x \in \R$ that $p^{(\scrd)}(x) = p^{(\scrd)}(0) \neq 0$. For every $\scrd \in \N_0$ we denote by $\cP(\scrd)$ the set of polynomials of degree $\scrd$ from $\R$ to $\R$.
\end{itemize}
We \nobs that for all $m, n \in \N_0$ it holds that $\cP(m) \cap \cP(n) \neq \emptyset$ if and only if $m = n$.

\subsubsection*{Terminology for polynomial neural networks (PNNs).}

In the following, we introduce the structural terminology for the parameters and hidden neurons of shallow \PNNs\ used throughout the paper, together with intuitive explanations; see \cref{fig:activity_types_of_neurons,figure_shallow_PNNs_illustration} below for activity-type diagram and graphical illustration of shallow \PNNs, respectively. The corresponding symbolic framework is given later in \cref{setting:SNNs} in \cref{subsec:PNNs_setting} and \cref{setting1} in \cref{subsec:PNNs_setting_special} below.

For positive integers $d, \dimension, \width \in \N$ with $\dimension = 3 \width + 1$ consider a shallow \PNN\ with $\width$ neurons on the hidden layer, parameter vector $\theta \allowbreak = (\theta_1, \allowbreak \dots, \allowbreak \theta_{\dimension}) \allowbreak \in \R^{\dimension}$, monomial activation function
\begin{equation}\label{eqn:intro_cA}
\textstyle \R \ni x \mapsto \cA(x) =  x^d \in \R
\end{equation}
of degree $d$, and realization function
\begin{equation}\label{eqn:intro:realization}
\textstyle \R \ni x \mapsto \cN^{\theta}(x) = \theta_{\dimension} + \sum_{\ell = 1}^{\width} \theta_{2 \width + \ell} \cA(\theta_{\width + \ell} + \theta_{\ell} x) \in \R
\end{equation}
associated with $\theta$. Then for every $i \in \Width = \{1, 2, \dots, \width\}$ we call
\begin{itemize}[itemsep=0.1em]
\item $\theta_i$ the $i^{\text{th}}$ \emph{inner weight} of $\theta$ and denote it by $\ww_i^{\theta}$,

\item $\theta_{\width + i}$ the $i^{\text{th}}$ \emph{inner bias} of $\theta$ and denote it by $\bb_i^{\theta}$,

\item $\theta_{2 \width + i}$ the $i^{\text{th}}$ \emph{outer weight} of $\theta$ and denote it by $\vv_i^{\theta}$,

\item $\theta_{\dimension}$ the \emph{outer bias} of $\theta$ and denote it by $\cc^{\theta}$,
\end{itemize}
we call $\vv_i^{\theta} [\ww_i^{\theta}]^d$ the \emph{\slope} of the $i^{\text{th}}$ hidden neuron of $\theta$, we call the sum of the slopes $\sum_{\fii = 1}^{\width} \vv_{\fii}^{\theta} [\ww_{\fii}^{\theta}]^d$ the \emph{\tslope} of $\theta$ and denote it by $\Slope{\theta}$, if $\ww_i^{\theta} \neq 0$ then we call $\nicefrac{- \bb_i^{\theta}}{\ww_i^{\theta}}$ the \emph{\pivot} (the \emph{finite pivot}) of the $i^{\text{th}}$ hidden neuron of $\theta$ and denote it by $\pvt{i}{\theta}$, we call the $i^{\text{th}}$ hidden neuron of $\theta$ \newline
\noindent
\begin{minipage}{0.49\linewidth}
\vspace{0.5cm}
\begin{itemize}[itemsep=0.1em]
\item \emph{\Active} if and only if $\ww_i^{\theta} \neq 0$,
\item \emph{\semiactive} if and only if $\ww_i^{\theta} = 0 \neq \bb_i^{\theta}$,
\item \emph{\inactive} if and only if $\ww_i^{\theta} = 0 = \bb_i^{\theta}$,
\item \emph{\visible} if and only if $\vv_i^{\theta} \neq 0$, and
\item \emph{\invisible} if and only if $\vv_i^{\theta} = 0$,
\end{itemize}
we call the inner weight $\ww_i^{\theta}$ and bias $\bb_i^{\theta}$ the \emph{activity} parameters and the outer weight $\vv_i^{\theta}$ the \emph{visibility} parameter of the $i^{\text{th}}$ hidden neuron of $\theta$, and we call the index set $\{\fii \in \Width \colon \ww_{\fii}^{\theta} \neq 0\}$ the \emph{active index set} of $\theta$ and denote it by $\cI^{\theta}$.
\end{minipage}
\hfill
\begin{minipage}{0.48\linewidth}
	\vspace{0.15cm}
	\centering
	\includegraphics[width=0.75\linewidth]{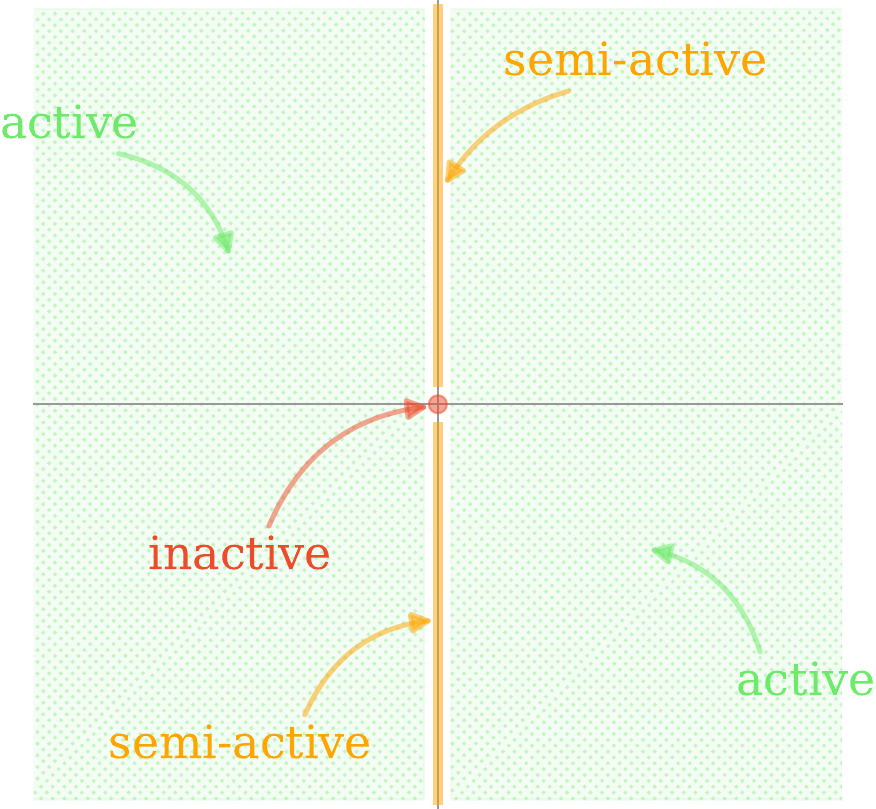}
	\vspace{-0.2cm}
	\captionof{figure}{Regions for activity types of hidden neurons as seen in the $(\ww_i^{\theta}, \bb_i^{\theta})$-plane.}
	\label{fig:activity_types_of_neurons}
\end{minipage}

\medskip
We \nobs that under these notations the realization function in \cref{eqn:intro:realization} can also be expressed as
\begin{equation}\label{eqn:intro:realization_new}
\textstyle \R \ni x \mapsto \cN^{\theta}(x) =\cc^{\theta}+ \sum_{\ell = 1}^{\width} \vv_{\ell}^{\theta} \cA(\bb_{\ell}^{\theta}+ \ww_{\ell}^{\theta} x) \in \R.
\end{equation}

The classification of hidden neurons of a shallow \PNN\ parameter vector $\theta \in \R^{\dimension}$ into ac\-ti\-vi\-ty and visibility types provides an intuitive framework for understanding their functional role.

The activity type of the hidden neuron -- \Active, \semiactive, or \inactive\ -- reflects the nature of the data it generates. When applying a monomial activation function in \cref{eqn:intro_cA} the \Active\ hidden neuron of $\theta$ produces input-dependent data via the surjection
\begin{equation}
\textstyle \R \ni x \mapsto \cA(\ww_i^{\theta} x \allowbreak + \bb_i^{\theta}) \allowbreak \in \cA(\R),
\end{equation}
the \semiactive\ hidden neuron of $\theta$ produces input-independent non-zero data via
\begin{equation}
\textstyle \R \ni x \mapsto \cA(\ww_i^{\theta} x \allowbreak + \bb_i^{\theta}) = \cA(\bb_i^{\theta}) \allowbreak \in \R \backslash \{0\},
\end{equation}
and the \inactive\ hidden neuron of $\theta$ produces input-independent zero data via
\begin{equation}
\textstyle \R \ni x \mapsto \cA(\ww_i^{\theta} x \allowbreak + \bb_i^{\theta}) = \cA(0) = 0.
\end{equation}

The visibility type of the hidden neuron -- \visible\ or \invisible\ -- determines whether the data $\cA(\ww_i^{\theta} x \allowbreak + \bb_i^{\theta}) \in \R$, generated by the corresponding hidden neuron of $\theta$ employing the activation function in \cref{eqn:intro_cA}, is transmitted to the realization function in \cref{eqn:intro:realization_new} associated with $\theta$.

\def\layersep{4cm}
\begin{figure}[H]
\centering
\begin{adjustbox}{width=\textwidth}
\begin{tikzpicture}[shorten >=1pt,->,draw=black!50, node distance=\layersep]
\tikzstyle{every pin edge}=[<-,shorten <=1pt]
\tikzstyle{input neuron}=[very thick, circle,draw=red, fill=red!20, line width=0.8pt, minimum size=15pt, inner sep=0pt]
\tikzstyle{output neuron}=[very thick, circle, draw=green, fill=green!20, line width=0.8pt, minimum size=30pt, inner sep=0pt]
\tikzstyle{hidden neuron}=[very thick, circle,draw=teal, fill=teal!20, line width=0.8pt, minimum size=20pt, inner sep=0pt]
\tikzstyle{annot} = [text width=12em, text centered]
\tikzstyle{annot2} = [text width=4em, text centered]

\node[hidden neuron](H0-1) at (\layersep, -1 cm) {$1^{\text{st}}$};
\node[hidden neuron](H0-2) at (\layersep, -2 cm) {$2^{\text{nd}}$};

\node(H0-dots1) at (\layersep, -2.65 cm) {\vdots};

\node[hidden neuron](H0-i) at (\layersep, -3.5 cm) {$i^{\text{th}}$};

\node(H0-dots2) at (\layersep, -4.15 cm) {\vdots};

\node[hidden neuron](H0-3) at (\layersep, -5 cm) {$\width^{\text{th}}$};

\node[input neuron] (I) at (0, -3 cm) {$x$};

\node[output neuron](O) at (2*\layersep, -3 cm) {$\cN^{\theta}(x)$};

\foreach \dest in {1,2,3}
\path[line width = 0.8 ] (I) edge (H0-\dest);

\path[line width = 0.8 ]
(I) edge node[midway, sloped, draw, line width=0.4pt, fill=red!20, rounded corners, inner sep=1pt, font=\tiny]
{$\ww_i^{\theta} x + \bb_i^{\theta}$} (H0-i);
\foreach \source in {1,2,3}
\path[line width = 0.8 ] (H0-\source) edge (O);

\path[line width = 0.8 ]
(H0-i) edge node[midway, sloped, draw, line width=0.4pt, fill=teal!20, rounded corners, inner sep=1pt, font=\tiny]
{$\vv_i^{\theta} (\ww_i^{\theta} x + \bb_i^{\theta})^d$} (O);

\node[annot,above of=H0-1, node distance=1cm, align=center] (hl) {Hidden layer\\($2^{\text{nd}}$ layer)};
\node[annot,left of=hl, align=center] {Input layer\\ ($1^{\text{st}}$ layer)};
\node[annot,right of=hl, align=center] {Output layer\\($3^{\text{rd}}$ layer)};

\node[annot2,below of=H0-3, node distance=1cm, align=center] (sl) {$\width$ neurons};
\node[annot2,left of=sl, align=center] {$1$ neuron};
\node[annot2,right of=sl, align=center] {$1$ neuron};
\end{tikzpicture}
\end{adjustbox}
\caption{Graphical illustration of the three-layer (shallow) \PNN\ architecture considered in this paper, with $\width \in \N$ hidden neurons and a monomial activation of degree $d \in \N$. This \PNN\ architecture has $1$-dimensional input and output layers and an $\width$-dimensional hidden layer.}
\label{figure_shallow_PNNs_illustration}
\end{figure}
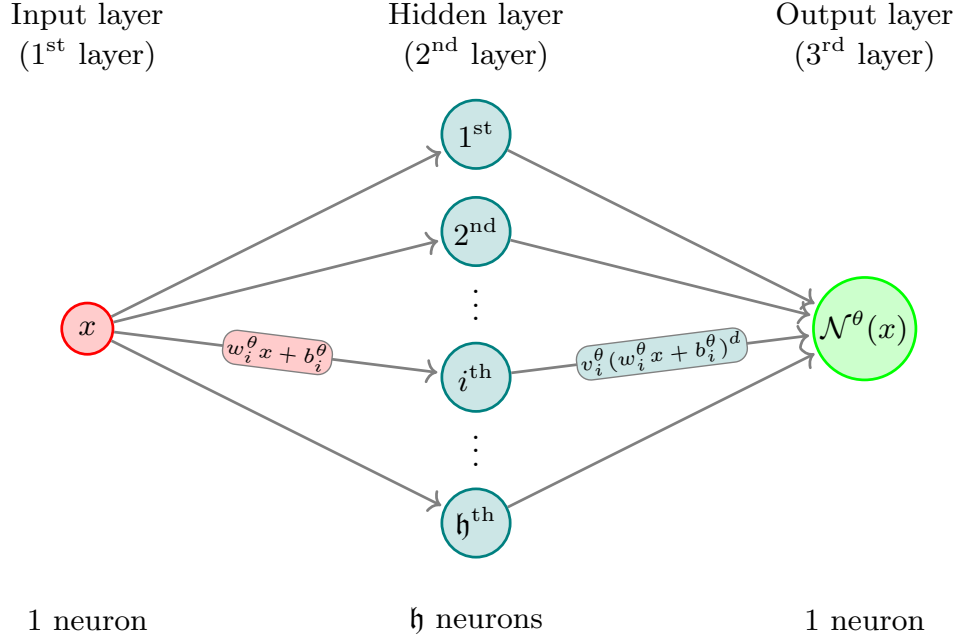

\section{Explicit derivatives of the loss function}
\label{sec:loss_derivatives}

In this section we introduce a mathematical setup for shallow \PNNs\ including the monomial activation function of degree $d \in \N$, the target function $f \in C(\R, \R)$, the realization functions $\cN^{\theta} \colon \R \to \R$, $\theta \in \R^{\dimension}$, and the loss function $\loss \colon \R^{\dimension} \to \R$. This framework, introduced  in \cref{setting:SNNs} in \cref{subsec:PNNs_setting} below,  will be used throughout the paper. We then employ the fact that the loss function is a polynomial in the network parameters to derive explicit formulas for its gradient and Hessian matrix in \cref{lemma:derivatives_of_Risk} and \cref{lemma:Hessian_of_Risk} in \cref{subsec:derivatives_of_loss} below. These are essential in determining and characterizing critical points in the subsequent sections. In particular, the strict convexity of the loss function with respect to the outer bias implies that it admits no local maxima (see \cref{eqn:lemma:Hessian_of_Risk:w_c_c_c} in \cref{lemma:Hessian_of_Risk} below). The proofs of \cref{lemma:derivatives_of_Risk} and \cref{lemma:Hessian_of_Risk} can be derived, e.g., using a special case of the Leibniz integral rule for differentiating parameter dependent Lebesgue integrals with constant boundaries, whose proof in a slightly more general setting can be found, for example, in \cite[Lemma~2.2 and Corollary~2.3]{IbragimovExistenceInfLocMin2022}.


\subsection{Polynomial neural networks (PNNs) with one hidden layer}
\label{subsec:PNNs_setting}

\begin{setting}[\resname{Shallow \PNNs}]\label{setting:SNNs}
Let $d, \dimension, \width \in \N$, $\scra \in \R$, $\scrb \in (\scra, \infty)$, $f \in C(\R, \R)$, $\Width \subseteq \N$ satisfy $\dimension = 3 \width + 1$ and $\Width \allowbreak = \{1, \allowbreak 2, \allowbreak \dots, \allowbreak \width\}$, for every $\theta = \allowbreak (\theta_1, \allowbreak \dots, \allowbreak \theta_{\dimension}) \allowbreak \in \R^{\dimension}$, $i \in \allowbreak \Width$ let $\ww_i^{\theta}, \allowbreak \bb_i^{\theta}, \allowbreak \vv_i^{\theta}, \allowbreak \cc^{\theta} \allowbreak \in \R$, $\cI^{\theta} \subseteq \Width$ satisfy $\ww_i^{\theta} \allowbreak = \theta_i$, $\bb_i^{\theta} \allowbreak = \theta_{\width + i}$, $\vv_i^{\theta} \allowbreak = \theta_{2 \width + i}$, $\cc^{\theta} \allowbreak = \theta_{\dimension}$, and $\cI^{\theta} \allowbreak = \{j \allowbreak \in \Width \colon \allowbreak \theta_j \neq 0\}$ and let $\scrC^{\theta}, \Slope{\theta} \in \R$, $\pvt{i}{\theta} \in \allowbreak \R \cup \{\infty\}$ satisfy
\begin{equation}
\textstyle \scrC^{\theta} = \cc^{\theta} + \sum_{j \in \Width \backslash \cI^{\theta}} \vv_j^{\theta} [\bb_j^{\theta}]^d, \quad \Slope{\theta} = \sum_{j = 1}^{\width} \vv_j^{\theta} [\ww_j^{\theta}]^d, \qandq \pvt{i}{\theta} =
        \begin{cases}
        \nicefrac{- \bb_i^{\theta}}{\ww_i^{\theta}} & \colon i \in \cI^{\theta} \\
        \infty & \colon i \notin \cI^{\theta},
        \end{cases}
\end{equation}
and let $\cN = (\cN^{\theta})_{\theta \in \R^{\dimension}} \colon \allowbreak \R^{\dimension} \allowbreak \to C(\R, \R)$ and $\loss \colon \R^{\dimension} \allowbreak \to \R$ satisfy for all $\theta \allowbreak \in \R^{\dimension}$ that
\begin{equation}\label{eqn:setting:SNNs:realization_risk}
\textstyle \bigl[\Forall x \in \R \colon \cN^{\theta}(x) = \allowbreak \cc^{\theta} + \allowbreak \sum_{i = 1}^{\width} \vv_i^{\theta} (\ww_i^{\theta} x + \bb_i^{\theta})^d\bigr] \qandq \bigl[\loss(\theta) = \int_{\scra}^{\scrb} \abs{\cN^{\theta}(x) - f(x)}^2 \, \d x\bigr].
\end{equation}
\end{setting}

\subsection{Explicit representations of the gradient and Hessian of the loss function}
\label{subsec:derivatives_of_loss}

\cfclear
\begin{lemma}[\resname{Gradient of the loss function}]\label{lemma:derivatives_of_Risk}
Assume \cref{setting:SNNs}. Then $\loss$ is a polynomial which satisfies for all $\theta \allowbreak = (\theta_1, \allowbreak \dots, \allowbreak \theta_{\dimension}) \allowbreak \in \R^{\dimension}$, $i \in \allowbreak \Width$ that
\begin{gather}\label{eqn:1st_derivative_w}
\textstyle (\frac{\partial}{\partial \theta_i} \loss)(\theta) = 2 d \vv_i^{\theta} \int_{\scra}^{\scrb} x (\ww_i^{\theta} x + \bb_i^{\theta})^{d - 1} [\cN^{\theta}(x) - f(x)] \, \d x, \\
\label{eqn:1st_derivative_b}
\textstyle (\frac{\partial}{\partial \theta_{\width + i}} \loss)(\theta) = 2 d \vv_i^{\theta} \int_{\scra}^{\scrb} (\ww_i^{\theta} x + \bb_i^{\theta})^{d - 1} [\cN^{\theta}(x) - f(x)] \, \d x, \\
\label{eqn:1st_derivative_v}
\textstyle (\frac{\partial}{\partial \theta_{2 \width + i}} \loss)(\theta) = 2 \int_{\scra}^{\scrb} (\ww_i^{\theta} x + \bb_i^{\theta})^d [\cN^{\theta}(x) - f(x)] \, \d x, \\
\label{eqn:1st_derivative_c}
\textstyle \andqq (\frac{\partial}{\partial \theta_{\dimension}} \loss)(\theta) = 2 \int_{\scra}^{\scrb} [\cN^{\theta}(x) - f(x)] \, \d x.
\end{gather}
\end{lemma}
\begin{cproof}{lemma:derivatives_of_Risk}
\Nobs that \cref{eqn:setting:SNNs:realization_risk} ensures that
\begin{equation}
\textstyle \R^{\dimension} \times \R \ni (\theta, x) \mapsto \cN^{\theta}(x) \in \R \qqandqq \R^{\dimension} \ni \theta \mapsto \loss(\theta) \in \R
\end{equation}
are both polynomials. Combining this with \cref{eqn:setting:SNNs:realization_risk} and Leibniz integral rule establishes \cref{eqn:1st_derivative_w,eqn:1st_derivative_b,eqn:1st_derivative_v,eqn:1st_derivative_c}.
\end{cproof}

\cfclear
\begin{lemma}[\resname{Hessian entries of the loss function}]\label{lemma:Hessian_of_Risk}
Assume \cref{setting:SNNs}. Then $\loss$ is a polynomial with no local maximum points and satisfies for all $\theta \allowbreak = (\theta_1, \allowbreak \dots, \allowbreak \theta_{\dimension}) \allowbreak \in \R^{\dimension}$, $i, j \in \allowbreak \Width$ that
\begin{equation}\label{eqn:lemma:Hessian_of_Risk:w_c_c_c}
\textstyle (\frac{\partial^2}{\partial \theta_i \partial \theta_{\dimension}} \loss)(\theta) = 2 d \vv_i^{\theta} \int_{\scra}^{\scrb} x (\ww_i^{\theta} x + \bb_i^{\theta})^{d - 1} \, \d x, \qquad (\frac{\partial^2}{\partial \theta_{\dimension}^2} \loss)(\theta) = 2 (\scrb - \scra),
\end{equation}
\begin{equation}\label{eqn:lemma:Hessian_of_Risk:b_c_v_c}
\textstyle (\frac{\partial^2}{\partial \theta_{\width + i} \partial \theta_{\dimension}} \loss)(\theta) = 2 d \vv_i^{\theta} \int_{\scra}^{\scrb} (\ww_i^{\theta} x + \bb_i^{\theta})^{d - 1} \, \d x, \quad (\frac{\partial^2}{\partial \theta_{2 \width + i} \partial \theta_{\dimension}} \loss)(\theta) = 2 \int_{\scra}^{\scrb} (\ww_i^{\theta} x + \bb_i^{\theta})^d \, \d x,
\end{equation}
\begin{multline}\label{eqn:lemma:Hessian_of_Risk:w_w}
\textstyle (\frac{\partial^2}{\partial \theta_i \partial \theta_j} \loss)(\theta) = 2 d^2 \vv_i^{\theta} \vv_j^{\theta} \int_{\scra}^{\scrb} x^2 (\ww_i^{\theta} x + \bb_i^{\theta})^{d - 1} (\ww_j^{\theta} x + \bb_j^{\theta})^{d - 1} \, \d x \\
\textstyle + 2 d (d - 1) \vv_i^{\theta} \indicator{\{i\}}(j) \int_{\scra}^{\scrb} x^2 (\ww_i^{\theta} x + \bb_i^{\theta})^{\max\{0, d - 2\}} [\cN^{\theta}(x) - f(x)] \, \d x,
\end{multline}
\begin{multline}\label{eqn:lemma:Hessian_of_Risk:w_b}
\textstyle (\frac{\partial^2}{\partial \theta_i \partial \theta_{\width + j}} \loss)(\theta) = 2 d^2 \vv_i^{\theta} \vv_j^{\theta} \int_{\scra}^{\scrb} x (\ww_i^{\theta} x + \bb_i^{\theta})^{d - 1} (\ww_j^{\theta} x + \bb_j^{\theta})^{d - 1} \, \d x \\
\textstyle + 2 d (d - 1) \vv_i^{\theta} \indicator{\{i\}}(j) \int_{\scra}^{\scrb} x (\ww_i^{\theta} x + \bb_i^{\theta})^{\max\{0, d - 2\}} [\cN^{\theta}(x) - f(x)] \, \d x,
\end{multline}
\begin{multline}\label{eqn:lemma:Hessian_of_Risk:b_b}
\textstyle (\frac{\partial^2}{\partial \theta_{\width + i} \partial \theta_{\width + j}} \loss)(\theta) = 2 d^2 \vv_i^{\theta} \vv_j^{\theta} \int_{\scra}^{\scrb} (\ww_i^{\theta} x + \bb_i^{\theta})^{d - 1} (\ww_j^{\theta} x + \bb_j^{\theta})^{d - 1} \, \d x \\
\textstyle + 2 d (d - 1) \vv_i^{\theta} \indicator{\{i\}}(j) \int_{\scra}^{\scrb} (\ww_i^{\theta} x + \bb_i^{\theta})^{\max\{0, d - 2\}} [\cN^{\theta}(x) - f(x)] \, \d x,
\end{multline}
\begin{multline}\label{eqn:lemma:Hessian_of_Risk:w_v}
\textstyle (\frac{\partial^2}{\partial \theta_i \partial \theta_{2 \width + j}} \loss)(\theta) = 2 d \vv_i^{\theta} \int_{\scra}^{\scrb} x (\ww_i^{\theta} x + \bb_i^{\theta})^{d - 1} (\ww_j^{\theta} x + \bb_j^{\theta})^d \, \d x \\
\textstyle + 2 d \indicator{\{i\}}(j) \int_{\scra}^{\scrb} x (\ww_i^{\theta} x + \bb_i^{\theta})^{d - 1} [\cN^{\theta}(x) - f(x)] \, \d x,
\end{multline}
\begin{multline}\label{eqn:lemma:Hessian_of_Risk:b_v}
\textstyle (\frac{\partial^2}{\partial \theta_{\width + i} \partial \theta_{2 \width + j}} \loss)(\theta) = 2 d \vv_i^{\theta} \int_{\scra}^{\scrb} (\ww_i^{\theta} x + \bb_i^{\theta})^{d - 1} (\ww_j^{\theta} x + \bb_j^{\theta})^d \, \d x \\
\textstyle + 2 d \indicator{\{i\}}(j) \int_{\scra}^{\scrb} (\ww_i^{\theta} x + \bb_i^{\theta})^{d - 1} [\cN^{\theta}(x) - f(x)] \, \d x,
\end{multline}
\begin{equation}\label{eqn:lemma:Hessian_of_Risk:v_v}
\textstyle \andqq (\frac{\partial^2}{\partial \theta_{2 \width + i} \partial \theta_{2 \width + j}} \loss)(\theta) = 2 \int_{\scra}^{\scrb} (\ww_i^{\theta} x + \bb_i^{\theta})^d (\ww_j^{\theta} x + \bb_j^{\theta})^d \, \d x.
\end{equation}
\end{lemma}
\begin{cproof}{lemma:Hessian_of_Risk}
\Nobs that \cref{eqn:setting:SNNs:realization_risk}, \cref{lemma:derivatives_of_Risk}, and Leibniz integral rule ensure that $\loss$ is a polynomial and establish \cref{eqn:lemma:Hessian_of_Risk:w_c_c_c,eqn:lemma:Hessian_of_Risk:b_c_v_c,eqn:lemma:Hessian_of_Risk:w_w,eqn:lemma:Hessian_of_Risk:w_b,eqn:lemma:Hessian_of_Risk:b_b,eqn:lemma:Hessian_of_Risk:w_v,eqn:lemma:Hessian_of_Risk:b_v,eqn:lemma:Hessian_of_Risk:v_v}. In addition, \nobs that \cref{eqn:lemma:Hessian_of_Risk:w_c_c_c} and the fact that $\scrb > \scra$ prove that $\loss$ has no local maxima. 
\end{cproof}

\section{Existence and non-existence of global minima for mo\-no\-mial activations}
\label{sec:global_minima}

In this section, for an arbitrary activation degree $d \in \N$ and non-constant affine linear target functions, we establish in \cref{prop:existence_of_global_minima} in \cref{subsec:existence_of_global_minima} below a sharp existence/non-existence criterion for global minimizers of the loss function in terms of the relation between the activation degree $d$ and the number $\width \in \N$ of hidden neurons of the \PNNs\ under consideration. In particular, we show that the infimum of the loss function is always zero. Moreover, we prove that this infimum is achievable if and only if the network contains at least $d$ \Active\ and \visible\ hidden neurons with pairwise distinct \pivot s, which subsequently requires $\width \ge d$. In contrast, it can be seen, in \cref{cor:non_existence_of_global_minima} in \cref{subsec:non_existence_of_global_minima} below which is the special case of \cref{prop:existence_of_global_minima}, that if the width $\width$ of the hidden layer is strictly smaller than the activation degree $d$, then the infimum cannot be attained and any sequence of parameters along which the loss function tends to zero necessarily diverges.

The proof of \cref{prop:existence_of_global_minima} relies on non-representability of polynomial target functions in the narrow regime ($\width < d$) in \cref{lemma:target_fn_non_reachable} and a representability of polynomial target functions in the wide regime ($\width \ge d$) in \cref{lemma:target_fn_reachable} in \cref{subsec:polynomial_target_functions_representability} below as well as \cref{lemma:inf_of_loss} in \cref{subsec:inf_of_loss} below where we construct a suitable sequence of parameters to show that the infimum of the loss function is  zero. The proofs of \cref{lemma:target_fn_non_reachable} and \cref{lemma:target_fn_reachable}, in turn, employ their respective sibling results concerning non-representability and representability of polynomials via shifted monomials, presented in \cref{lemma:non_representability_poly} and \cref{lemma:representability_poly}, respectively, in \cref{subsec:polynomials_representability_via_monomials} below. Finally, \cref{lemma:non_representability_poly} and \cref{lemma:representability_poly} make use of a well-known formula for the determinant of square Vandermonde matrices which we include, only for the convenience of the reader, in \cref{lemma:Vandermonde} in \cref{subsec:polynomials_representability_via_monomials} below.

\subsection{On the representability of polynomials via shifted monomials}
\label{subsec:polynomials_representability_via_monomials}

\cfclear
\begin{lemma}[\resname{Vandermonde determinant}]\label{lemma:Vandermonde}
Let $n \in \N$, $x_1, x_2, \dots, x_n \in \R$. Then
\begin{equation}
        \det\!
        \begin{psmallmatrix}
        1 & 1 & \cdots & 1 \\[1ex]
        x_1 & x_2 & \cdots & x_n \\[1ex]
        [x_1]^2 & [x_2]^2 & \cdots & [x_n]^2 \\
        \vdots & \vdots & \ddots & \vdots \\[1ex]
        [x_1]^{n - 1} & [x_2]^{n - 1} & \cdots & [x_n]^{n - 1}
        \end{psmallmatrix}\! = \prod_{1 \le i < j \le n} (x_j - x_i).
\end{equation}
\end{lemma}

\cfclear
\begin{lemma}[\resname{Non-representability}]\label{lemma:non_representability_poly}
Let $d \in \N$, $\width \in \Z \cap [0, d)$ and let $P \in \cup_{i = 1}^{d - \width} \cP(i)$. Then it holds for all $\fc, \allowbreak \fa_i, \fb_i \in \R$, $i \in \{1, 2, \dots, \width\}$, that
\begin{equation}\label{eqn:lemma:non_representability_poly}
\textstyle \# \{x \in \R \colon \fc + \sum_{i = 1}^{\width} \fa_i (x + \fb_i)^d = P(x)\} < \infty.
\end{equation}
\end{lemma}
\begin{cproof}{lemma:non_representability_poly}
Throughout this proof let $\bfd \in \N$, $\fc, \allowbreak \fa_1, \allowbreak \fa_2, \allowbreak \dots, \allowbreak \fa_{\width}, \allowbreak \fb_1, \allowbreak \fb_2, \allowbreak \dots, \allowbreak \fb_{\width} \allowbreak \in \R$ satisfy $\bfd \le d - \width$ and $P \in \cP(\bfd)$, assume\footnote{\Nobs that since $P$ is a non-constant polynomial \cref{eqn:lemma:non_representability_poly} is trivial in the case $\width = 0$. Moreover, if some of the shifts $\fb_i$, $i \in \{1, 2, \dots, \width\}$, are equal, their corresponding coefficients $\fa_i$ can be summed up. This forms an expression of the same form as in \cref{eqn:lemma:non_representability_poly} with pairwise distinct shifts, where the total number of summands is strictly less than $\width$.} w.l.o.g.\ that $\width > 0$ and that $\fb_1, \fb_2, \dots, \fb_{\width}$ are pairwise distinct, and assume for the sake of contradiction that
\begin{equation}\label{eqn:lemma:non_representability_poly:set_R}
\textstyle \{x \in \R \colon \fc + \sum_{i = 1}^{\width} \fa_i (x + \fb_i)^d = P(x)\} = \R.
\end{equation}
\Nobs that the fact that for all $a, b \in \R$ it holds that
\begin{equation}\label{eqn:lemma:non_representability_poly:binomial}
\textstyle (a + b)^d = \sum_{i = 0}^d \binom{d}{i} a^i b^{d - i}
\end{equation}
ensures for all $x \in \R$ that
\begin{equation}
\textstyle \fc + \sum_{i = 1}^{\width} \fa_i (x + \fb_i)^d = \fc + \sum_{i = 0}^d x^{d - i} \binom{d}{i} (\sum_{j = 1}^{\width} \fa_j [\fb_j]^i).
\end{equation}
The fact that $P \in \cP(\bfd)$ and \cref{eqn:lemma:non_representability_poly:set_R} \hence demonstrate that
\begin{equation}
\textstyle \Forall i \in \{0, 1, \dots, d - \bfd - 1\} \colon \sum_{j = 1}^{\width} \fa_j [\fb_j]^i = 0.
\end{equation}
Combining this with \cref{lemma:Vandermonde} (applied for every $i \in \{1, 2, \dots, \width\}$ with $n \with \width$, $x_i \with \fb_i$ in the notation of \cref{lemma:Vandermonde}), the assumption that $\fb_1, \fb_2, \dots, \fb_{\width}$ are pairwise distinct, and the fact that $\bfd \allowbreak \le d \allowbreak - \width$ shows that
\begin{equation}
\textstyle \width \le d - \bfd, \qquad
        \det\!
        \begin{psmallmatrix}
        1 & 1 & \cdots & 1 \\[1ex]
        \fb_1 & \fb_2 & \cdots & \fb_{\width} \\[1ex]
        [\fb_1]^2 & [\fb_2]^2 & \cdots & [\fb_{\width}]^2 \\
        \vdots & \vdots & \ddots & \vdots \\[1ex]
        [\fb_1]^{\width - 1} & [\fb_2]^{\width - 1} & \cdots & [\fb_{\width}]^{\width - 1}
        \end{psmallmatrix}\! = \prod_{1 \le i < j \le \width} (\fb_j - \fb_i) \neq 0,
\end{equation}
\begin{equation}
\andqq
        \begin{psmallmatrix}
        1 & 1 & \cdots & 1 \\[1ex]
        \fb_1 & \fb_2 & \cdots & \fb_{\width} \\[1ex]
        [\fb_1]^2 & [\fb_2]^2 & \cdots & [\fb_{\width}]^2 \\
        \vdots & \vdots & \ddots & \vdots \\[1ex]
        [\fb_1]^{\width - 1} & [\fb_2]^{\width - 1} & \cdots & [\fb_{\width}]^{\width - 1}
        \end{psmallmatrix}
        \begin{psmallmatrix}
        \fa_1 \\[1ex]
        \fa_2 \\[1ex]
        \fa_3 \\
        \vdots\\[1ex]
        \fa_{\width}
        \end{psmallmatrix}\! = 0.
\end{equation}
\Hence that $\fa_1 = \fa_2 = \cdots = \fa_{\width} = 0$. The fact that $P \in \cP(\bfd)$ \hence contradicts the assumption in \cref{eqn:lemma:non_representability_poly:set_R}. This establishes \cref{eqn:lemma:non_representability_poly}.
\end{cproof}

\cfclear
\begin{lemma}[\resname{Representability}]\label{lemma:representability_poly}
Let $d \in \N$, $P \in \cup_{i = 0}^d \cP(i)$. Then there exist $\fc, \allowbreak \fa_1, \allowbreak \fa_2, \allowbreak \dots, \allowbreak \fa_d, \allowbreak \fb_1, \allowbreak \fb_2, \allowbreak \dots, \allowbreak \fb_d \allowbreak \in \R$ such that
\begin{equation}\label{eqn:lemma:representability_poly}
\textstyle \Forall x \in \R \colon \fc + \sum_{i = 1}^d \fa_i (x + \fb_i)^d = P(x).
\end{equation}
\end{lemma}
\begin{cproof}{lemma:representability_poly}
Throughout this proof let $p_0, p_1, \dots, p_d \in \R$ satisfy for all $x \in \R$ that $P(x) = \sum_{i = 0}^d p_i x^i$ and let $\fc, \allowbreak \fa_1, \allowbreak \fa_2, \allowbreak \dots, \allowbreak \fa_d, \allowbreak \fb_1, \allowbreak \fb_2, \allowbreak \dots, \allowbreak \fb_d \allowbreak \in \R$ satisfy
\begin{gather}\label{eqn:lemma:representability_poly:fc}
\textstyle \prod_{1 \le i < j \le d} (\fb_j - \fb_i) \neq 0, \qquad \fc = p_0 - \sum_{j = 1}^{d} \fa_j [\fb_j]^d, \\
\label{eqn:lemma:representability_poly:fa_fb}
\textstyle \andqq
        \begin{psmallmatrix}
        \fa_1 \\[1ex]
        \fa_2 \\[1ex]
        \fa_3 \\
        \vdots \\[1ex]
        \fa_d
        \end{psmallmatrix}\! =\!
        \begin{psmallmatrix}
        1 & 1 & \cdots & 1 \\[1ex]
        \fb_1 & \fb_2 & \cdots & \fb_d \\[1ex]
        [\fb_1]^2 & [\fb_2]^2 & \cdots & [\fb_d]^2 \\
        \vdots & \vdots & \ddots & \vdots \\[1ex]
        [\fb_1]^{d - 1} & [\fb_2]^{d - 1} & \cdots & [\fb_d]^{d - 1}
        \end{psmallmatrix}^{-1}
        \begin{psmallmatrix}
        p_d \binom{d}{0}^{-1} \\[1ex]
        p_{d - 1} \binom{d}{1}^{-1} \\[1ex]
        p_{d - 2} \binom{d}{2}^{-1} \\
        \vdots \\[1ex]
        p_1 \binom{d}{d - 1}^{-1}
        \end{psmallmatrix}
\end{gather}
(cf.\ \cref{lemma:Vandermonde}). \Nobs that \cref{eqn:lemma:representability_poly:fa_fb} ensures that
\begin{equation}
\textstyle
        \begin{psmallmatrix}
        1 & 1 & \cdots & 1 \\[1ex]
        \fb_1 & \fb_2 & \cdots & \fb_d \\[1ex]
        [\fb_1]^2 & [\fb_2]^2 & \cdots & [\fb_d]^2 \\
        \vdots & \vdots & \ddots & \vdots \\[1ex]
        [\fb_1]^{d - 1} & [\fb_2]^{d - 1} & \cdots & [\fb_d]^{d - 1}
        \end{psmallmatrix}
        \begin{psmallmatrix}
        \fa_1 \\[1ex]
        \fa_2 \\[1ex]
        \fa_3 \\
        \vdots \\[1ex]
        \fa_d
        \end{psmallmatrix}\! =\!
        \begin{psmallmatrix}
        p_d \binom{d}{0}^{-1} \\[1ex]
        p_{d - 1} \binom{d}{1}^{-1} \\[1ex]
        p_{d - 2} \binom{d}{2}^{-1} \\
        \vdots \\[1ex]
        p_1 \binom{d}{d - 1}^{-1}
        \end{psmallmatrix}\!.
\end{equation}
This demonstrates for all $i \in \{0, 1, \dots, d - 1\}$ that
\begin{equation}
\textstyle \sum_{j = 1}^d \fa_j [\fb_j]^i = p_{d - i} \binom{d}{i}^{-1}.
\end{equation}
The fact that for all $a, b \in \R$ it holds that
\begin{equation}
\textstyle (a + b)^d = \sum_{i = 0}^d \binom{d}{i} a^i b^{d - i}
\end{equation}
and \cref{eqn:lemma:representability_poly:fc} \hence prove for all $x \in \R$ that
\begin{multline}
\textstyle \fc + \fa_1 (x + \fb_1)^d + \fa_2 (x + \fb_2)^d + \dots + \fa_d (x + \fb_d)^d = \fc + \sum_{i = 0}^d x^{d - i} \binom{d}{i} (\sum_{j = 1}^d \fa_j [\fb_j]^i) \\
\textstyle = \fc + \sum_{j = 1}^d \fa_j [\fb_j]^d + \sum_{i = 0}^{d - 1} x^{d - i} \binom{d}{i} (\sum_{j = 1}^d \fa_j [\fb_j]^i) = p_0 + \sum_{i = 0}^{d - 1} x^{d - i} p_{d - i} = P(x).
\end{multline}
This establishes \cref{eqn:lemma:representability_poly}.
\end{cproof}

\subsection{On the representability of polynomial target functions via PNNs}
\label{subsec:polynomial_target_functions_representability}

\begin{lemma}[\resname{Non-representability}]\label{lemma:target_fn_non_reachable}
Assume \cref{setting:SNNs} with $\width < d$ and $f \in \cup_{i = 1}^{d - \width} \cP(i)$. Then it holds that $\{\theta \in \R^{\dimension} \colon \cN^{\theta} = f\} = \emptyset$.
\end{lemma}
\begin{cproof}{lemma:target_fn_non_reachable}
\Nobs that \cref{eqn:setting:SNNs:realization_risk} ensures for all $\theta \in \R^{\dimension}$, $x \in \R$ that $\# (\cI^{\theta}) \le \width$ and
\begin{multline}
\textstyle \cN^{\theta}(x) = \cc^{\theta} + \sum_{i = 1}^{\width} \vv_i^{\theta} (\ww_i^{\theta} x + \bb_i^{\theta})^d \\
\textstyle = \cc^{\theta} + \sum_{i \in \Width \backslash \cI^{\theta}} \vv_i^{\theta} [\bb_i^{\theta}]^d + \sum_{i \in \cI^{\theta}} \vv_i^{\theta} (\ww_i^{\theta} x + \bb_i^{\theta})^d = \scrC^{\theta} + \sum_{i \in \cI^{\theta}} \vv_i^{\theta} [\ww_i^{\theta}]^d (x - \pvt{i}{\theta})^d.
\end{multline}
\cref{lemma:non_representability_poly} (applied for every $\theta \in \R^{\dimension}$ with $d \with d$, $\width \with \# (\cI^{\theta})$, $P \with f$ in the notation of \cref{lemma:non_representability_poly}) and the assumption that $\width < d$ and $f \in \cup_{i = 1}^{d - \width} \cP(i)$ \hence establishes
\begin{equation}
\textstyle \{\theta \in \R^{\dimension} \colon \cN^{\theta} = f\} = \emptyset.
\end{equation}
\end{cproof}

\begin{lemma}[\resname{Representability}]\label{lemma:target_fn_reachable}
Assume \cref{setting:SNNs} with $\width \ge d$. Then
\begin{enumerate}[label=\textnormal{(\roman*)}]
\item
\label{item1:lemma:target_fn_reachable} $f \in \cup_{i = 0}^d \cP(i)$ implies $\{\theta \in \R^{\dimension} \colon \cN^{\theta} = f\} \neq \emptyset$ and

\item
\label{item2:lemma:target_fn_reachable} $f \in \cP(1)$ implies for all $\theta \in \R^{\dimension}$ with $\cN^{\theta} = f$ that $\# (\cup_{i \in \cI^{\theta}, \vv_i^{\theta} \neq 0} \{\pvt{i}{\theta}\}) \ge d$.
\end{enumerate}
\end{lemma}
\begin{cproof}{lemma:target_fn_reachable}
\Nobs that \cref{eqn:setting:SNNs:realization_risk} ensures for all $\theta \in \R^{\dimension}$, $x \in \R$ that $\# (\cI^{\theta}) \le \width$ and
\begin{multline}\label{eqn:lemma:target_fn_reachable:cN}
\textstyle \cN^{\theta}(x) = \cc^{\theta} + \sum_{i = 1}^{\width} \vv_i^{\theta} (\ww_i^{\theta} x + \bb_i^{\theta})^d \\
\textstyle = \cc^{\theta} + \sum_{i \in \Width \backslash \cI^{\theta}} \vv_i^{\theta} [\bb_i^{\theta}]^d + \sum_{i \in \cI^{\theta}} \vv_i^{\theta} (\ww_i^{\theta} x + \bb_i^{\theta})^d = \scrC^{\theta} + \sum_{i \in \cI^{\theta}} \vv_i^{\theta} [\ww_i^{\theta}]^d (x - \pvt{i}{\theta})^d.
\end{multline}
\cref{lemma:representability_poly} (applied with $d \with d$, $P \with f$ in the notation of \cref{lemma:representability_poly}) and the assumption that $\width \ge d$ \hence imply that, if $f \in \cup_{i = 0}^d \cP(i)$, then there exists $\theta \in \R^{\dimension}$ with $\# (\cI^{\theta}) \ge d$ and $\cN^{\theta} = f$. This establishes \cref{item1:lemma:target_fn_reachable}. Moreover, \nobs that \cref{lemma:non_representability_poly} (applied for every $\theta \in \R^{\dimension}$ with $d \with d$, $\width \with \# (\cup_{i \in \cI^{\theta}, \vv_i^{\theta} \neq 0} \{\pvt{i}{\theta}\})$, $P \with f$ in the notation of \cref{lemma:non_representability_poly}) and \cref{eqn:lemma:target_fn_reachable:cN} show that, if $f \in \cP(1)$, then for all $\theta \in \R^{\dimension}$ with $\# (\cup_{i \in \cI^{\theta}, \vv_i^{\theta} \neq 0} \{\pvt{i}{\theta}\}) < d$ it holds that
\begin{equation}
\textstyle \# \{x \in \R \colon \cN^{\theta}(x) = f(x)\} < \infty.
\end{equation}
This establishes \cref{item2:lemma:target_fn_reachable}.
\end{cproof}

\subsection{Explicit infimum of the loss function for affine target functions}
\label{subsec:inf_of_loss}

\begin{lemma}\label{lemma:inf_of_loss}
Assume \cref{setting:SNNs} with $f \in \Span\{1, \id_{\R}\}$. Then it holds that $\inf_{\theta \in \R^{\dimension}} \loss(\theta) = 0$.
\end{lemma}
\begin{cproof}{lemma:inf_of_loss}
Throughout this proof let $\fa, \fb \in \R$ satisfy for all $x \in \R$ that $f(x) = \fa x + \fb$ and for every $n \in \N$ let $\phi_n \in \R^{\dimension}$ satisfy
\begin{equation}\label{eqn:lemma:inf_of_loss:phi}
\textstyle \cc^{\phi_n} = - \frac{n}{d} + \fb, \quad \ww_1^{\phi_n} = \frac{\fa}{n}, \quad \bb_1^{\phi_n} = 1, \quad \vv_1^{\phi_n} = \frac{n}{d}, \qandq \sum_{i = 2}^{\width} \abs{\vv_i^{\phi_n}} = 0.
\end{equation}
\Nobs that \cref{eqn:setting:SNNs:realization_risk} and \cref{eqn:lemma:inf_of_loss:phi} ensure for all $n \in \N$, $x \in \R$ that
\begin{multline}
\textstyle \cN^{\phi_n}(x) = - \frac{n}{d} + \fb + \frac{n}{d} (1 + \frac{\fa x}{n})^d = \fb + \frac{n}{d} \bigl[(1 + \frac{\fa x}{n})^d - 1\bigr] \\
\textstyle = \fb + \frac{n}{d} \sum_{i = 1}^d \binom{d}{i} (\frac{\fa x}{n})^i = \fb + \fa x + \frac{n}{d} \sum_{i = 2}^d \binom{d}{i} (\frac{\fa x}{n})^i = f(x) + \frac{n}{d} \sum_{i = 2}^d \binom{d}{i} (\frac{\fa x}{n})^i
\end{multline}
\begin{equation}
\textstyle \andqq \loss(\phi_n) = \int_{\scra}^{\scrb} \abs{\cN^{\phi_n}(y) - f(y)}^2 \, \d y = \int_{\scra}^{\scrb} \Abs{\frac{n}{d} \sum_{i = 2}^d \binom{d}{i} (\frac{\fa y}{n})^i}^2 \, \d y = \cO(\frac{1}{n^2}).
\end{equation}
\Hence that $\limsup_{n \to \infty} \loss(\phi_n) = 0$.
\end{cproof}

\subsection{Non-existence of global minima in the narrow regime}
\label{subsec:non_existence_of_global_minima}

\begin{corollary}\label{cor:non_existence_of_global_minima}
Assume \cref{setting:SNNs} with $\width < d$ and $f \in \cP(1)$. Then
\begin{equation}\label{eqn:cor:non_existence_of_global_minima}
\textstyle \inf_{\theta \in \R^{\dimension}} \loss(\theta) = 0 \qqandqq \{\theta \in \R^{\dimension} \colon \loss(\theta) = 0\} = \emptyset.
\end{equation}
\end{corollary}
\begin{cproof}{cor:non_existence_of_global_minima}
\Nobs that \cref{eqn:setting:SNNs:realization_risk} ensures for all $\theta \in \R^{\dimension}$ that $\loss(\theta) = 0$ if and only if $\cN^{\theta} = f$. This implies that
\begin{equation}
\textstyle \{\theta \in \R^{\dimension} \colon \loss(\theta) = 0\} = \{\theta \in \R^{\dimension} \colon \cN^{\theta} = f\}.
\end{equation}
The assumption that $\width < d$, the assumption that $f \in \cP(1)$, \cref{lemma:target_fn_non_reachable}, and \cref{lemma:inf_of_loss} \hence establish \cref{eqn:cor:non_existence_of_global_minima}.
\end{cproof}

\subsection{Existence of global minima in the wide regime}
\label{subsec:existence_of_global_minima}

\begin{proposition}[\resname{Criterion for global minima}]\label{prop:existence_of_global_minima}
Assume \cref{setting:SNNs} with $f \in \cP(1)$. Then
\begin{enumerate}[label=\textnormal{(\roman*)}]
\item
\label{item1:prop:existence_of_global_minima} it holds that $\inf_{\theta \in \R^{\dimension}} \loss(\theta) = 0$,

\item
\label{item2:prop:existence_of_global_minima} it holds that $\{\theta \in \R^{\dimension} \colon \loss(\theta) = 0\} \neq \emptyset$ if and only if $\width \ge d$,

\item
\label{item3:prop:existence_of_global_minima} it holds for all $\theta \in \R^{\dimension}$ with $\loss(\theta) = 0$ that $\# (\cup_{i \in \cI^{\theta}, \vv_i^{\theta} \neq 0} \{\pvt{i}{\theta}\}) - d \ge 0 = (\nabla \loss)(\theta)$, and

\item
\label{item4:prop:existence_of_global_minima} if $\width < d$ it holds for all $(\phi_n)_{n \in \N} \subseteq \R^{\dimension}$ with $\limsup_{n \to \infty} \loss(\phi_n) = 0$ that
\begin{equation}
\textstyle \liminf_{n \to \infty} \norm{\phi_n} = \infty.
\end{equation}
\end{enumerate}
\end{proposition}
\begin{cproof}{prop:existence_of_global_minima}
First \nobs that \cref{lemma:inf_of_loss} establishes \cref{item1:prop:existence_of_global_minima}. Next \nobs that \cref{eqn:setting:SNNs:realization_risk} ensures for all $\theta \in \R^{\dimension}$ that $\loss(\theta) = 0$ if and only if $\cN^{\theta} = f$. This implies that
\begin{equation}\label{eqn:prop:existence_of_global_minima:Loss_cN}
\textstyle \{\theta \in \R^{\dimension} \colon \loss(\theta) = 0\} = \{\theta \in \R^{\dimension} \colon \cN^{\theta} = f\}.
\end{equation}
The assumption that $f \in \cP(1)$ and \cref{item1:lemma:target_fn_reachable} in \cref{lemma:target_fn_reachable} \hence demonstrate that
\begin{equation}\label{eqn:prop:existence_of_global_minima:2_1}
\textstyle [\{\theta \in \R^{\dimension} \colon \loss(\theta) = 0\} \neq \emptyset] \qquad \Leftarrow \qquad [\width \ge d].
\end{equation}
In addition, \nobs that \cref{eqn:prop:existence_of_global_minima:Loss_cN}, \cref{lemma:target_fn_non_reachable}, and the assumption that $f \in \cP(1)$ assure that
\begin{equation}
\textstyle [\{\theta \in \R^{\dimension} \colon \loss(\theta) = 0\} \neq \emptyset] \qquad \Rightarrow \qquad [\width \ge d].
\end{equation}
This and \cref{eqn:prop:existence_of_global_minima:2_1} establish \cref{item2:prop:existence_of_global_minima}. Furthermore, \nobs that \cref{eqn:prop:existence_of_global_minima:Loss_cN}, \cref{lemma:derivatives_of_Risk}, \cref{item2:lemma:target_fn_reachable} in \cref{lemma:target_fn_reachable}, and the assumption that $f \in \cP(1)$ establish \cref{item3:prop:existence_of_global_minima}. In the next step we prove \cref{item4:prop:existence_of_global_minima} by contradiction. Let $(\phi_n)_{n \in \N} \subseteq \R^{\dimension}$ satisfy
\begin{equation}\label{eqn:prop:existence_of_global_minima:phi}
\textstyle \limsup_{n \to \infty} \loss(\phi_n) = 0 \qqandqq \liminf_{n \to \infty} \norm{\phi_n} < \infty
\end{equation}
and assume $\width < d$. \Nobs that \cref{eqn:prop:existence_of_global_minima:phi} and Bolzano-Weierstrass theorem in $\R^{\dimension}$ ensure that $(\phi_n)_{n \in \N}$ has a convergent subsequence. Let $\vartheta \in \R^{\dimension}$ and let $(n_k)_{k \in \N} \subseteq \N$ be an increasing sequence such that
\begin{equation}\label{eqn:prop:existence_of_global_minima:vartheta}
\textstyle \limsup_{k \to \infty} \norm{\phi_{n_k} - \vartheta} = 0.
\end{equation}
\Nobs that \cref{eqn:prop:existence_of_global_minima:phi}, \cref{eqn:prop:existence_of_global_minima:vartheta}, the fact that $\loss$ is a polynomial (cf.\ \cref{lemma:derivatives_of_Risk}), and the triangle inequality show that
\begin{equation}
\textstyle 0 = \limsup_{k \to \infty} (\abs{\loss(\vartheta) - \loss(\phi_{n_k})} + \loss(\phi_{n_k})) \ge \loss(\vartheta) \ge 0.
\end{equation}
\Hence that $\loss(\vartheta) = 0$. This and \cref{item2:prop:existence_of_global_minima} contradict  the assumption that $\width < d$ and \hence establish \cref{item4:prop:existence_of_global_minima}.
\end{cproof}

\section{Determination of critical points for the cubic activation}
\label{sec:determination_cubic_activation}

In this section, for convenience, we consider in \cref{setting1} in \cref{subsec:PNNs_setting_special} below a special case of \cref{setting:SNNs} in \cref{subsec:PNNs_setting} above, where we restrict to the cubic monomial activation, the identity target function, and the unit interval as the domain of integration for the loss function~$\loss$ (corresponding to $d = 3$, $\scra = 0$, $\scrb = 1$, and $f = \id_{\R}$ in \cref{setting:SNNs}). Employing \cref{lemma:derivatives_of_Risk} in \cref{subsec:derivatives_of_loss} above, we then determine all critical points of the loss function in \cref{lemma:realizations} in \cref{subsec:all_realizations} below in terms of the number of distinct \pivot s of \Active\ hidden neurons and explicit realization functions.

In particular, we show that only critical points with all hidden neurons being \inactive\ or \invisible\ \semiactive\ can yield a constant realization function and at least two \Active\ hidden neurons with distinct \pivot s are required for the exact representation of the target function. We prove, in addition, that critical points with two distinct \pivot s of \Active\ hidden neurons must have, in total, at least three \Active\ and \visible\ hidden neurons with pairwise distinct \pivot s. This guarantees the absence of critical points with exactly two distinct \pivot s of \Active\ hidden neurons. Consequently, we demonstrate that there are three different non-trivial realizations of critical points with all \pivot s of \Active\ hidden neurons coinciding.

\cref{lemma:realizations} is a collection of results presented in \cref{lemma:constant_realization} in \cref{subsec:constant_realization} below, \cref{lemma:distinctness_of_kinks} in \cref{subsec:number_of_pivots} below, and \cref{lemma:non_constant_realization} in \cref{subsec:all_realizations} below, where we separately investigate critical points with no \Active\ hidden neurons, with at least three distinct \pivot s of \Active\ hidden neurons, and with all \pivot s of \Active\ hidden neurons coinciding, respectively. In addition, only for the sake of completeness, we also include \cref{cor:kinks_loc_min_saddle} in \cref{subsec:number_of_pivots} below, which is a direct consequence of the latter results and explicitly shows that there are no critical points with exactly two distinct \pivot s of \Active\ hidden neurons.

\subsection{PNNs with the cubic activation and identity target function}
\label{subsec:PNNs_setting_special}

\begin{setting}\label{setting1}
Let $\dimension, \width \in \N$, $\Width \subseteq \N$ satisfy $\dimension = 3 \width + 1$ and $\Width \allowbreak = \{1, \allowbreak 2, \allowbreak \dots, \allowbreak \width\}$, for every $\theta = \allowbreak (\theta_1, \allowbreak \dots, \allowbreak \theta_{\dimension}) \allowbreak \in \R^{\dimension}$, $i \in \allowbreak \Width$ let $\ww^{\theta}_i, \allowbreak \bb^{\theta}_i, \allowbreak \vv^{\theta}_i, \allowbreak \cc^{\theta} \allowbreak \in \R$, $\cI^{\theta} \subseteq \Width$ satisfy $\ww_i^{\theta} = \theta_i$, $\bb_i^{\theta} = \theta_{\width + i}$, $\vv_i^{\theta} = \theta_{2 \width + i}$, $\cc^{\theta} = \theta_{\dimension}$, and $\cI^{\theta} \allowbreak = \{j \allowbreak \in \Width \colon \allowbreak \theta_j \neq 0\}$ and let $\scrC^{\theta}, \Slope{\theta} \in \R$, $\pvt{i}{\theta} \in \allowbreak \R \cup \{\infty\}$ satisfy
\begin{equation}
\textstyle \scrC^{\theta} = \cc^{\theta} + \sum_{j \in \Width \backslash \cI^{\theta}} \vv_j^{\theta} [\bb_j^{\theta}]^3, \quad \Slope{\theta} = \sum_{j = 1}^{\width} \vv_j^{\theta} [\ww_j^{\theta}]^3, \qandq \pvt{i}{\theta} =
        \begin{cases}
        \nicefrac{- \bb_i^{\theta}}{\ww_i^{\theta}} & \colon i \in \cI^{\theta} \\
        \infty & \colon i \notin \cI^{\theta},
        \end{cases}
\end{equation}
and let $\cN = (\cN^{\theta})_{\theta \in \R^{\dimension}} \colon \allowbreak \R^{\dimension} \allowbreak \to C(\R, \R)$ and $\loss \colon \R^{\dimension} \to \R$ satisfy for all $\theta \in \R^{\dimension}$ that
\begin{equation}\label{eqn:setting1:realization_risk}
\textstyle \bigl[\Forall x \in \R \colon \cN^{\theta}(x) = \allowbreak \cc^{\theta} + \allowbreak \sum_{i = 1}^{\width} \vv_i^{\theta} (\ww_i^{\theta} x + \bb_i^{\theta})^3\bigr] \qandq \bigl[\loss(\theta) = \int_0^1 \abs{\cN^{\theta}(x) - x}^2 \, \d x\bigr].
\end{equation}
\end{setting}

\subsection{Structural degeneracies for odd activations}
\label{subsec:constant_realization}

\cfclear
\begin{lemma}[\resname{Constant realizations of critical points}]\label{lemma:constant_realization}
Assume \cref{setting:SNNs} with $\nicefrac{(d + 1)}{2} \in \N$, $\scra = 0$, $\scrb = 1$, and $f = \id_{\R}$ and let $\theta \in \R^{\dimension}$ satisfy $(\nabla \loss)(\theta) = 0$. Then the following statements are equivalent:
\begin{enumerate}[label=\textnormal{(\roman*)}]
\item
\label{item1:lemma:constant_realization} It holds that $\sum_{i = 1}^{\width} (\abs{\ww^{\theta}_i} + \abs{\vv^{\theta}_i \bb^{\theta}_i}) = 0 = \frac{1}{2} - \cc^{\theta}$.

\item
\label{item2:lemma:constant_realization} It holds that $\sum_{i = 1}^{\width} \abs{\ww^{\theta}_i} = 0$.

\item
\label{item3:lemma:constant_realization} It holds that $\cI^{\theta} = \emptyset$.

\item
\label{item4:lemma:constant_realization} It holds that $\sum_{i = 1}^{\width} \abs{\vv_i^{\theta} \ww_i^{\theta}} = 0$.

\item
\label{item5:lemma:constant_realization} It holds for all $x \in \R$ that $\cN^{\theta}(x) = \cN^{\theta}(0)$.

\item
\label{item6:lemma:constant_realization} It holds for all $x \in \R$ that $\cN^{\theta}(x) = \nicefrac{1}{2}$.
\end{enumerate}
\end{lemma}
\begin{cproof}{lemma:constant_realization}
Throughout this proof let $m \in \N$ satisfy $d = 2m - 1$. \Nobs that \cref{eqn:setting:SNNs:realization_risk} ensures that 
\ref{item1:lemma:constant_realization} $\Rightarrow$ \ref{item2:lemma:constant_realization} $\Leftrightarrow$ \ref{item3:lemma:constant_realization} $\Rightarrow$ \ref{item4:lemma:constant_realization} $\Rightarrow$ \ref{item5:lemma:constant_realization}. Moreover, \nobs that \cref{lemma:derivatives_of_Risk} (applied with $\scra \with 0$, $\scrb \with 1$, $f \with \id_{\R}$ in the notation of \cref{lemma:derivatives_of_Risk}), the fact that for all $c \in \R$ it holds that $\int_0^1 [c - x] \d x = c - \nicefrac{1}{2}$, and the fact that $(\nabla \loss)(\theta) = 0$ demonstrate that \ref{item5:lemma:constant_realization} $\Leftrightarrow$ \ref{item6:lemma:constant_realization}. Next it suffices to prove that \ref{item6:lemma:constant_realization} $\Rightarrow$ \ref{item1:lemma:constant_realization}. In the following we assume that
\begin{equation}\label{eqn:lemma:constant_realization:const}
\textstyle \Forall x \in \R \colon \cN^{\theta}(x) = \nicefrac{1}{2}.
\end{equation}
\Nobs that the fact that for all $n \in \N$, $a, b \in \R$ with $a + b > 0$ it holds that $a^{2n - 1} + b^{2n - 1} > 0$ proves for all $n \in \N$, $\fq \in \R$, $y \in (0, \infty)$ that $(y - \frac{1}{2} + \fq)^{2n - 1} + (y + \frac{1}{2} - \fq)^{2n - 1} > 0$. This and the change-of-variables formula ensure for all $\fq \in \R$ that
\begin{equation}
\begin{split}
& \textstyle \int_0^1 (x - \fq)^{2m - 1} (x - \frac{1}{2}) \, \d x = \int_{- \frac{1}{2}}^{\frac{1}{2}} (x + \frac{1}{2} - \fq)^{2m - 1} x \, \d x \\
& \textstyle = \int_{- \frac{1}{2}}^0 (x + \frac{1}{2} - \fq)^{2m - 1} x \, \d x + \int_0^{\frac{1}{2}} (x + \frac{1}{2} - \fq)^{2m - 1} x \, \d x \\
& \textstyle = \int_0^{\frac{1}{2}} (- x + \frac{1}{2} - \fq)^{2m - 1} (- x) \, \d x + \int_0^{\frac{1}{2}} (x + \frac{1}{2} - \fq)^{2m - 1} x \, \d x \\
& \textstyle = \int_0^{\frac{1}{2}} x \bigl[(x - \frac{1}{2} + \fq)^{2m - 1} + (x + \frac{1}{2} - \fq)^{2m - 1}\bigr] \, \d x > 0.
\end{split}
\end{equation}
\cref{lemma:derivatives_of_Risk} (applied with $\scra \with 0$, $\scrb \with 1$, $f \with \id_{\R}$ in the notation of \cref{lemma:derivatives_of_Risk}) and \cref{eqn:lemma:constant_realization:const} \hence show for all $i \in \Width$ that
\begin{equation}
\textstyle \ww_i^{\theta} \neq 0 \qquad \Rightarrow \qquad (\frac{\partial}{\partial \vv_i^{\theta}} \loss)(\theta) \neq 0.
\end{equation}
This and the fact that $(\nabla \loss)(\theta) = 0$ prove that $\sum_{i = 1}^{\width} \abs{\ww_i^{\theta}} = 0$. Combining this with \cref{eqn:lemma:constant_realization:const}, the fact that $\int_0^1 x (\frac{1}{2} - x) \, \d x = -\frac{1}{12} < 0$, and \cref{lemma:derivatives_of_Risk} (applied with $\scra \with 0$, $\scrb \with 1$, $f \with \id_{\R}$ in the notation of \cref{lemma:derivatives_of_Risk}) assures for all $i \in \Width$ that
\begin{equation}
\textstyle \vv_i^{\theta} \bb_i^{\theta} \neq 0 \qquad \Rightarrow \qquad (\frac{\partial}{\partial \ww_i^{\theta}} \loss)(\theta) \neq 0.
\end{equation}
The fact that $(\nabla \loss)(\theta) = 0$ \hence implies that $\sum_{i = 1}^{\width} \abs{\vv_i^{\theta} \bb_i^{\theta}} = 0$. \Hence that
\begin{equation}
\textstyle \sum_{i = 1}^{\width} (\abs{\ww^{\theta}_i} \allowbreak + \abs{\vv^{\theta}_i \bb^{\theta}_i}) = 0.
\end{equation}
This and \cref{eqn:lemma:constant_realization:const} establish that \ref{item6:lemma:constant_realization} $\Rightarrow$ \ref{item1:lemma:constant_realization}.
\end{cproof}

\begin{remark}
In general, \PNNs\ may represent constant functions with some hidden neurons being \Active\ and/or \visible. However, \cref{lemma:constant_realization} above shows in case of odd monomial activations that, at a critical point $\theta \in \R^{\dimension}$ of $\loss$, this situation cannot occur and constancy of the realization function associated with $\theta$ forces each $i^{\text{th}}$ hidden neuron of $\theta$ to be \inactive\ ($\ww_i^{\theta} = 0 = \bb_i^{\theta}$) or \invisible\ \semiactive\ ($\ww_i^{\theta} = \vv_i^{\theta} = 0 \neq \bb_i^{\theta}$).
\end{remark}

We \nobs that while \cref{lemma:constant_realization} above and its direct consequence \cref{lemma:constant_realization_gen} in \cref{subsec:structural_characterization} below hold for all odd-degree monomial activations, in the remainder of this section we restrict to the situation of \cref{setting1}, corresponding to $d = 3$, $\scra = 0$, $\scrb = 1$, and $f = \id_{\R}$ in \cref{setting:SNNs}, and determine the critical points.

\subsection{On the number of distinct \pivot s at critical points}
\label{subsec:number_of_pivots}

\cfclear
\begin{lemma}\label{lemma:distinctness_of_kinks}
Assume \cref{setting1} and let $\theta \in \allowbreak \R^{\dimension}$. Then the following statements are equivalent:
\begin{enumerate}[label=\textnormal{(\roman*)}]
\item
\label{item1:lemma:distinctness_of_kinks} It holds that $\# (\cup_{i \in \cI^{\theta}} \{\pvt{i}{\theta}\}) - 1 > 0 = (\nabla \loss)(\theta)$.

\item
\label{item2:lemma:distinctness_of_kinks} It holds that $\# (\cup_{i \in \cI^{\theta}, \vv_i^{\theta} \neq 0} \{\pvt{i}{\theta}\}) - 3 \ge 0 = (\nabla \loss)(\theta)$.

\item
\label{item3:lemma:distinctness_of_kinks} It holds that $\cN^{\theta} = \id_{\R}$.

\item
\label{item4:lemma:distinctness_of_kinks} It holds that $\loss(\theta) = 0$.
\end{enumerate}
\end{lemma}
\begin{cproof}{lemma:distinctness_of_kinks}
First \nobs that \ref{item2:lemma:distinctness_of_kinks} $\Rightarrow$ \ref{item1:lemma:distinctness_of_kinks}. Next \nobs that \cref{prop:existence_of_global_minima} (applied with $d \with 3$, $\scra \with 0$, $\scrb \with 1$, $f \with \id_{\R}$ in the notation of \cref{prop:existence_of_global_minima}) and \cref{eqn:setting:SNNs:realization_risk} ensure that \ref{item4:lemma:distinctness_of_kinks} $\Leftrightarrow$ \ref{item3:lemma:distinctness_of_kinks} $\Rightarrow$ \ref{item2:lemma:distinctness_of_kinks}. Next it suffices to prove that \ref{item1:lemma:distinctness_of_kinks} $\Rightarrow$ \ref{item3:lemma:distinctness_of_kinks}. In the following we assume that
\begin{equation}\label{eqn:lemma:distinctness_of_kinks:item1_eq}
\textstyle \# (\bigcup_{i \in \cI^{\theta}} \{\pvt{i}{\theta}\}) - 1 > 0 = (\nabla \loss)(\theta).
\end{equation}
\Nobs that \cref{eqn:lemma:distinctness_of_kinks:item1_eq} and \cref{lemma:constant_realization} (applied with $d \with 3$ in the notation of \cref{lemma:constant_realization}) demonstrate that there exist $i_1, i_2 \in \Width$ which satisfy $\pvt{i_1}{\theta} \neq \pvt{i_2}{\theta}$ and $\vv_{i_1}^{\theta} \allowbreak \ww_{i_1}^{\theta} \allowbreak \ww_{i_2}^{\theta} \allowbreak \neq 0$. This, \cref{eqn:lemma:distinctness_of_kinks:item1_eq}, and \cref{lemma:derivatives_of_Risk} (applied with $d \with 3$, $\scra \with 0$, $\scrb \with 1$, $f \with \id_{\R}$ in the notation of \cref{lemma:derivatives_of_Risk}) prove for all $\ell \in \{1, 2\}$ that
\begin{equation}\label{eqn:lemma:classification_for_k=3:derivative=0}
\textstyle \int_0^1 [\cN^{\theta}(x) - x] \, \d x = \int_0^1 (x - \pvt{i_{\ell}}{\theta})^3 [\cN^{\theta}(x) - x] \, \d x = \int_0^1 (x - \pvt{i_1}{\theta})^2 [\cN^{\theta}(x) - x] \, \d x = 0.
\end{equation}
\Hence that
\begin{equation}
\begin{split}
\textstyle 0 & \textstyle = \frac{1}{\pvt{i_2}{\theta} - \pvt{i_1}{\theta}} \int_0^1 [(x - \pvt{i_1}{\theta})^3 - (x - \pvt{i_2}{\theta})^3] [\cN^{\theta}(x) - x] \, \d x \\
& \textstyle = \int_0^1 [(x - \pvt{i_1}{\theta})^2 + (x - \pvt{i_1}{\theta}) (x - \pvt{i_2}{\theta}) + (x - \pvt{i_2}{\theta})^2] [\cN^{\theta}(x) - x] \, \d x \\
& \textstyle = \int_0^1 [3 (x - \pvt{i_1}{\theta})^2 + 3 (x - \pvt{i_1}{\theta}) (\pvt{i_1}{\theta} - \pvt{i_2}{\theta}) + (\pvt{i_1}{\theta} - \pvt{i_2}{\theta})^2] [\cN^{\theta}(x) - x] \, \d x \\
& \textstyle = 3 \int_0^1 (x - \pvt{i_1}{\theta})^2 [\cN^{\theta}(x) - x] \, \d x + 3 (\pvt{i_1}{\theta} - \pvt{i_2}{\theta}) \int_0^1 x [\cN^{\theta}(x) - x] \, \d x \\
& \textstyle \quad + [- 3 \pvt{i_1}{\theta} (\pvt{i_1}{\theta} - \pvt{i_2}{\theta}) + (\pvt{i_1}{\theta} - \pvt{i_2}{\theta})^2] \int_0^1 [\cN^{\theta}(x) - x] \, \d x \\
& \textstyle = 3 (\pvt{i_1}{\theta} - \pvt{i_2}{\theta}) \int_0^1 x [\cN^{\theta}(x) - x] \, \d x = 0.
\end{split}
\end{equation}
Combining this with \cref{eqn:lemma:classification_for_k=3:derivative=0} assures for all $m \in \{0, 1, 2, 3\}$ that
\begin{equation}
\textstyle \int_0^1 x^m [\cN^{\theta}(x) - x] \, \d x = 0.
\end{equation}
The fact that $\cN^{\theta} \in \cup_{i = 0}^3 \cP(i)$ \hence proves that $\cN^{\theta} = \id_{\R}$. This and \cref{eqn:lemma:distinctness_of_kinks:item1_eq} show that \ref{item1:lemma:distinctness_of_kinks} $\Rightarrow$ \ref{item3:lemma:distinctness_of_kinks}. The fact that \ref{item4:lemma:distinctness_of_kinks} $\Leftrightarrow$ \ref{item3:lemma:distinctness_of_kinks} $\Rightarrow$ \ref{item2:lemma:distinctness_of_kinks} $\Rightarrow$ \ref{item1:lemma:distinctness_of_kinks} \hence establishes that \ref{item1:lemma:distinctness_of_kinks} $\Leftrightarrow$ \ref{item2:lemma:distinctness_of_kinks} $\Leftrightarrow$ \ref{item3:lemma:distinctness_of_kinks} $\Leftrightarrow$ \ref{item4:lemma:distinctness_of_kinks}.
\end{cproof}

\cfclear
\begin{corollary}\label{cor:kinks_loc_min_saddle}
Assume \cref{setting1} and let $\theta \in \R^{\dimension}$ satisfy $(\nabla \loss)(\theta) = 0$. Then
\begin{enumerate}[label=\textnormal{(\roman*)}]
\item
\label{item1:cor:kinks_loc_min_saddle} it holds that $\# (\cup_{i \in \cI^{\theta}} \{\pvt{i}{\theta}\}) \neq 2$ and

\item
\label{item2:cor:kinks_loc_min_saddle} it holds that $\# (\cup_{i \in \cI^{\theta}} \{\pvt{i}{\theta}\}) \le 1$ if and only if $\cN^{\theta} \neq \id_{\R}$.
\end{enumerate}
\end{corollary}
\begin{cproof}{cor:kinks_loc_min_saddle}
\Nobs that \cref{lemma:distinctness_of_kinks} establishes \cref{item1:cor:kinks_loc_min_saddle,item2:cor:kinks_loc_min_saddle}.
\end{cproof}

\subsection{Complete description of realization functions at critical points}
\label{subsec:all_realizations}

\cfclear
\begin{lemma}[\resname{Necessity of non-linear realizations}]\label{lemma:non_constant_realization}
Assume \cref{setting1} and let $\theta \in \R^{\dimension}$ satisfy $(\nabla \loss)(\theta) = 0$ and $\# (\cup_{i \in \cI^{\theta}} \{\pvt{i}{\theta}\}) = 1$. Then it holds that $\sum_{i \in \Width \backslash \cI^{\theta}} \abs{\vv^{\theta}_i \bb^{\theta}_i} = 0$ and exactly one of the following holds:
\begin{enumerate}[label=\textnormal{(\roman*)}]
\item\label{item1:lemma:non_constant_realization} It holds for all $x \in \R$ that $\cN^{\theta}(x) = \frac{35 - 4 \sqrt{7}}{70} + \frac{14}{5} [x - \frac{7 - \sqrt{7}}{14}]^3$.

\item\label{item2:lemma:non_constant_realization} It holds for all $x \in \R$ that $\cN^{\theta}(x) = \frac{35 + 4 \sqrt{7}}{70} + \frac{14}{5} [x - \frac{7 + \sqrt{7}}{14}]^3$.

\item\label{item3:lemma:non_constant_realization} It holds for all $x \in \R$ that $\cN^{\theta}(x) = \frac{1}{2} + \frac{28}{5} [x - \frac{1}{2}]^3$.
\end{enumerate}
\end{lemma}
\begin{cproof}{lemma:non_constant_realization}
\Nobs that the assumption that $\# (\cup_{i \in \cI^{\theta}} \{\pvt{i}{\theta}\}) = 1$ and \cref{lemma:constant_realization} (applied with $d \with 3$ in the notation of \cref{lemma:constant_realization}) ensure that there exist $\fc, \fa, \fq \in \R$ which satisfy for all $x \allowbreak \in \R$ that $\cN^{\theta}(x) \allowbreak = \fc + \fa (x - \fq)^3$ and $\fa \neq 0$. This and \cref{lemma:derivatives_of_Risk} (applied with $d \with 3$, $\scra \with 0$, $\scrb \with 1$, $f \with \id_{\R}$ in the notation of \cref{lemma:derivatives_of_Risk}) demonstrate that
\begin{equation}
\begin{gathered}
\textstyle \int_0^1 [\fc + \fa (x - \fq)^3 - x] \, \d x = 0, \qquad \int_0^1 (x - \fq)^2 [\fc + \fa (x - \fq)^3 - x] \, \d x = 0, \\
\textstyle \andqq \int_0^1 (x - \fq)^3 [\fc + \fa (x - \fq)^3 - x] \, \d x = 0.
\end{gathered}
\end{equation}
\Hence that
\begin{equation}
\begin{gathered}
\textstyle \fc + \frac{\fa [(1 - \fq)^4 - \fq^4]}{4} - \frac{1}{2} = 0, \quad \frac{\fc [(1 - \fq)^3 + \fq^3]}{3} + \frac{\fa [(1 - \fq)^6 - \fq^6]}{6} - \bigl[\frac{(1 - \fq)^4 - \fq^4}{4} + \frac{\fq [(1 - \fq)^3 + \fq^3]}{3}\bigr] = 0, \\
\textstyle \andqq \frac{\fc [(1 - \fq)^4 - \fq^4]}{4} + \frac{\fa [(1 - \fq)^7 + \fq^7]}{7} - \big[ \frac{(1 - \fq)^5 + \fq^5}{5} + \frac{\fq [(1 - \fq)^4 - \fq^4]}{4}\big] = 0.
\end{gathered}
\end{equation}
This shows that
\begin{equation}
\begin{gathered}
\textstyle \fc = \frac{1}{4} (2 - \fa + 4 \fa \fq - 6 \fa \fq^2 +4 \fa \fq^3), \qquad \frac{1}{12} (2 \fq - 1) [1 - \fa + 3 \fa \fq (1 - \fq)] = 0, \\
\textstyle \andqq \frac{1}{40} (- 3 + 10 \fq - 10 \fq^2) + \frac{1}{112} \fa (9 - 56 \fq + 140 \fq^2 - 168 \fq^3 + 84 \fq^4) = 0.
\end{gathered}
\end{equation}
\Hence that
\begin{equation}\label{eqn:lemma:non_constant_realization:a_b_c_values}
\begin{gathered}
\textstyle \bigl([\fq = \frac{1}{2}] \wedge [\fa = \frac{28}{5}] \wedge [\fc = \frac{1}{2}]\bigr) \vee \bigl(\bigl[\fq = \frac{7 - \sqrt{7}}{14}\bigr] \wedge [\fa = \frac{14}{5}] \wedge \bigl[\fc = \frac{35 - 4 \sqrt{7}}{70}\bigr]\bigr) \\
\textstyle \vee \bigl(\bigl[\fq = \frac{7 + \sqrt{7}}{14}\bigr] \wedge [\fa = \frac{14}{5}] \wedge \bigl[\fc = \frac{35 + 4 \sqrt{7}}{70}\bigr]\bigr).
\end{gathered}
\end{equation}
In the next step we combine \cref{lemma:derivatives_of_Risk} (applied with $d \with 3$, $\scra \with 0$, $\scrb \with 1$, $f \with \id_{\R}$ in the notation of \cref{lemma:derivatives_of_Risk}) with \cref{eqn:lemma:non_constant_realization:a_b_c_values} and the assumption that $(\nabla \loss)(\theta) = 0$ to obtain for all $i \in \Width \backslash \cI^{\theta} = \{j \in \Width \colon \ww_j^{\theta} = 0\}$ that
\begin{equation}
\begin{split}
\textstyle 0 & \textstyle = (\frac{\partial}{\partial \theta_i} \loss)(\theta) = 6 \vv_i^{\theta} (\bb_i^{\theta})^2 \int_0^1 x [\cN^{\theta}(x) - x] \, \d x = 6 \vv_i^{\theta} (\bb_i^{\theta})^2 \int_0^1 x [\fc + \fa (x - \fq)^3 - x] \, \d x \\
& \textstyle =
\begin{cases}
6 \vv_i^{\theta} (\bb_i^{\theta})^2 (- \frac{1}{75}) & \colon [\fq = \frac{1}{2}] \wedge [\fa = \frac{28}{5}] \wedge [\fc = \frac{1}{2}], \\
6 \vv_i^{\theta} (\bb_i^{\theta})^2 (- \frac{7}{300}) & \colon \bigl[\fq = \frac{7 - \sqrt{7}}{14}\bigr] \wedge [\fa = \frac{14}{5}] \wedge \bigl[\fc = \frac{35 - 4 \sqrt{7}}{70}\bigr], \\
6 \vv_i^{\theta} (\bb_i^{\theta})^2 (- \frac{7}{300}) & \colon \bigl[\fq = \frac{7 + \sqrt{7}}{14}\bigr] \wedge [\fa = \frac{14}{5}] \wedge \bigl[\fc = \frac{35 + 4 \sqrt{7}}{70}\bigr].
\end{cases}
\end{split}
\end{equation}
This ensures for all $i \in \Width \backslash \cI^{\theta}$ that $\vv^{\theta}_i \bb^{\theta}_i = 0$.
\end{cproof}

\cfclear
\begin{lemma}[\resname{Sufficiency of non-linear realizations}]\label{lemma:non_constant_realization_sufficiency}
Assume \cref{setting1}, let $\theta \in \R^{\dimension}$ satisfy $\# (\cup_{i \in \cI^{\theta}} \{\pvt{i}{\theta}\}) = 1$ and $\sum_{i \in \Width \backslash \cI^{\theta}} \abs{\vv^{\theta}_i \bb^{\theta}_i} = 0$, and assume that $\cN^\theta$ has one of the forms in \cref{item1:lemma:non_constant_realization,item2:lemma:non_constant_realization,item3:lemma:non_constant_realization} in \cref{lemma:non_constant_realization}. Then it holds that  $(\nabla \loss)(\theta) = 0$.
\end{lemma}
\begin{cproof}{lemma:non_constant_realization_sufficiency}
Throughout this proof let $\fq \in \R$ satisfy $\{\fq\} = \cup_{i \in \cI^{\theta}} \{\pvt{i}{\theta}\}$. \Nobs that substituting realizations from \cref{item1:lemma:non_constant_realization,item2:lemma:non_constant_realization,item3:lemma:non_constant_realization} in \cref{lemma:non_constant_realization} into \cref{eqn:1st_derivative_b}--\cref{eqn:1st_derivative_c} gives vanishing partial derivatives w.r.t.\ all biases and outer weights
\begin{equation}
\textstyle \int_0^1 (\cN^{\theta}(x) - x) \d x = \int_0^1 (x - \fq)^2 (\cN^{\theta}(x) - x) \d x = \int_0^1 (x - \fq)^3 (\cN^{\theta}(x) - x) \d x = 0.
\end{equation}
This, \cref{eqn:1st_derivative_w}, and the fact that for all $x \in \R$ it holds that $x(x - \fq)^2 = (x - \fq)^3 + \fq (x - \fq)^2$ establish vanishing partial derivatives w.r.t.\ the inner weights of \Active\ hidden neurons
\begin{equation}
\textstyle \int_0^1 x (x - \fq)^2 (\cN^{\theta}(x) - x) \d x = 0.
\end{equation}
Moreover, \cref{eqn:1st_derivative_w} and the condition $\sum_{i \in \Width \backslash \cI^{\theta}} \abs{\vv^{\theta}_i \bb^{\theta}_i} = 0$ annihilates the remaining derivatives w.r.t.\ the inner weights of the hidden neurons which are not \Active. \Hence that $(\nabla \loss)(\theta) = 0$.
\end{cproof}

\cfclear
\begin{lemma}[\resname{All forms of realizations}]\label{lemma:realizations}
Assume \cref{setting1} and let $\theta \in \R^{\dimension}$. Then it holds that $(\nabla \loss)(\theta) = 0$ if and only if one of the following holds:
\begin{enumerate}[label=\textnormal{(\roman*)}]
\item $\cI^{\theta} = \emptyset$, $\sum_{i = 1}^{\width} \abs{\vv^{\theta}_i \bb^{\theta}_i} = 0$, and $\Forall x \in \R \colon \cN^{\theta}(x) = \nicefrac{1}{2}$.

\item $\# (\cup_{i \in \cI^{\theta}} \{\pvt{i}{\theta}\}) = 1$, $\sum_{i \in \Width \backslash \cI^{\theta}} \abs{\vv^{\theta}_i \bb^{\theta}_i} = 0$, and $\Forall x \in \R \colon \cN^{\theta}(x) = \frac{35 - 4 \sqrt{7}}{70} + \frac{14}{5} [x - \frac{7 - \sqrt{7}}{14}]^3$.

\item $\# (\cup_{i \in \cI^{\theta}} \{\pvt{i}{\theta}\}) = 1$, $\sum_{i \in \Width \backslash \cI^{\theta}} \abs{\vv^{\theta}_i \bb^{\theta}_i} = 0$, and $\Forall x \in \R \colon \cN^{\theta}(x) = \frac{35 + 4 \sqrt{7}}{70} + \frac{14}{5} [x - \frac{7 + \sqrt{7}}{14}]^3$.

\item $\# (\cup_{i \in \cI^{\theta}} \{\pvt{i}{\theta}\}) = 1$, $\sum_{i \in \Width \backslash \cI^{\theta}} \abs{\vv^{\theta}_i \bb^{\theta}_i} = 0$, and $\Forall x \in \R \colon \cN^{\theta}(x) = \frac{1}{2} + \frac{28}{5} [x - \frac{1}{2}]^3$.

\item $\# (\cup_{i \in \cI^{\theta}, \vv_i^{\theta} \neq 0} \{\pvt{i}{\theta}\}) \ge 3$ and $\Forall x \in \R \colon \cN^{\theta}(x) = x$.
\end{enumerate}
\end{lemma}

\Nobs that \cref{lemma:realizations} is a direct consequence of \cref{lemma:derivatives_of_Risk}, \cref{lemma:constant_realization}, \cref{lemma:distinctness_of_kinks}, \cref{lemma:non_constant_realization}, and \cref{lemma:non_constant_realization_sufficiency}.

\section{Classification of critical points for the cubic activation}
\label{sec:classification_cubic_activation}

In this section we recall our findings from the previous sections to establish in \cref{thm:complete_classification_cubic} in \cref{subsec:complete_classification} below a complete classification of all critical points of the loss function for the cubic monomial activation and non-constant affine linear target functions. In particular, we first perform a characterization of critical points, determined in \cref{sec:determination_cubic_activation} above, in the situation of \cref{setting1} -- where we consider the identity target function and the unit interval $[0, 1]$ as the domain of integration for the loss function. We then employ the invariance property of the loss landscape under suitable affine transformations, presented in \cref{lemma:canonical_extension} in \cref{subsec:affine_scaling} below, to easily extend the characterization results to the situation of \cref{setting:SNNs} with an arbitrary non-constant affine linear target function $f \colon \R \to \R$ and an arbitrary domain $\emptyset \neq [\scra, \scrb] \subseteq \R$ of integration for the loss function $\loss \colon \R^{\dimension} \to \R$.

In \cref{thm:complete_classification_cubic}, we recall from \cref{sec:loss_derivatives} that the loss function is a polynomial in the \PNN\ parameter vector which admits no local maximizers, and we then show that the target function is exactly represented at global minimizers, which requires at least three \Active\ and \visible\ hidden neurons with pairwise distinct \pivot s, that critical points with no \Active\ hidden neurons correspond to saddle points, and, consequently, that non-global local minimizers and non-trivial saddle points can appear only in networks for which all \pivot s of \Active\ hidden neurons coincide. Moreover, non-global local minimizers require all hidden neurons to be \Active\ and \visible\ with exactly one hidden neuron having a \slope\ sign matching that of the target function.

\cref{thm:complete_classification_cubic} is simply a collection of all structural characterization results for critical points presented separately in \cref{prop:global_minima_main}, \cref{prop:saddle_points}, and \cref{prop:local_minima} in \cref{subsec:structural_characterization} below. The structural characterization of
\begin{itemize}
\item global minimizers, established in \cref{prop:global_minima_main}, is a direct consequence of \cref{lemma:distinctness_of_kinks} in \cref{subsec:number_of_pivots} above,

\item saddle points, established in \cref{prop:saddle_points}, uses \cref{lemma:saddle_points_constant_realization}, \cref{lemma:saddle_point_1}, and \cref{cor:saddle_point_2} in \cref{subsec:saddle_points} below, and 

\item non-global local minimizers, established in \cref{prop:local_minima}, uses \cref{lemma:realizations} in \cref{subsec:all_realizations} above and \cref{lemma:loc_min} in \cref{subsec:local_minima} below.
\end{itemize}

The proofs of \cref{lemma:saddle_points_constant_realization}, \cref{lemma:saddle_point_1}, and \cref{cor:saddle_point_2} are derived via the construction of suitable ascent and descent directions by partially involving the explicit Hessian formulas in \cref{lemma:Hessian_of_Risk} in \cref{subsec:derivatives_of_loss} above.  Our  proof of \cref{lemma:loc_min} employs, for computational simplicity, the reducibility of an arbitrary perturbation of a network with only \Active\ hidden neurons to effective parameters, established in \cref{lemma:invariant_perturbation} in \cref{subsec:reducing_perturbations} below, in order to guarantee the absence of descent directions.

\subsection{Saddle points via construction of ascent and descent directions}
\label{subsec:saddle_points}

\cfclear
\begin{lemma}\label{lemma:saddle_points_constant_realization}
Assume \cref{setting1} and let $\theta \in \R^{\dimension}$ satisfy for all $x \in \R$ that $(\nabla \loss)(\theta) = 0$ and $\cN^{\theta}(x) = \nicefrac{1}{2}$. Then
\begin{enumerate}[label=\textnormal{(\roman*)}]
\item
\label{item1:lemma:saddle_points_constant_realization} it holds that $\sum_{i = 1}^{\width} (\abs{\ww^{\theta}_i} + \abs{\vv^{\theta}_i \bb^{\theta}_i}) = 0 = \frac{1}{2} - \cc^{\theta}$ and

\item
\label{item2:lemma:saddle_points_constant_realization} it holds that $\theta$ is a saddle point of $\loss$.
\end{enumerate}
\end{lemma}
\begin{cproof}{lemma:saddle_points_constant_realization}
First \nobs that \cref{lemma:constant_realization} (applied with $d \with 3$ in the notation of \cref{lemma:constant_realization}) establishes \cref{item1:lemma:saddle_points_constant_realization}. Next we prove that $\theta$ is a saddle point of $\loss$. Let the perturbation $\phi \in C(\R, \R^{\dimension})$ of $\theta$ be given for all $s \in \R$ by
\begin{equation}
\begin{gathered}\label{eqn:lemma:saddle_points_constant_realization:theta_s}
\textstyle \ww_1^{\phi_s} = 
        \begin{cases}
        s \bb_1^{\theta} & \!\!\!\!\colon \vv_1^{\theta} = \ww_1^{\theta} = 0 \neq \bb_1^{\theta} \\
        \sqrt[3]{\frac{s^5}{\vv_1^{\theta}}} & \!\!\!\!\colon \ww_1^{\theta} = \bb_1^{\theta} = 0 \neq \vv_1^{\theta} \\
        \sqrt[3]{s^4} & \!\!\!\!\colon \ww_1^{\theta} = \bb_1^{\theta} = \vv_1^{\theta} = 0,
        \end{cases}
\qquad \textstyle \bb_1^{\phi_s} =
        \begin{cases}
        \bb_1^{\theta} & \!\!\!\!\colon \vv_1^{\theta} = \ww_1^{\theta} = 0 \neq \bb_1^{\theta} \\
        \sqrt[3]{\frac{s^2}{\vv_1^{\theta}}} & \!\!\!\!\colon \ww_1^{\theta} = \bb_1^{\theta} = 0 \neq \vv_1^{\theta} \\
        \sqrt[3]{s} & \!\!\!\!\colon \ww_1^{\theta} = \bb_1^{\theta} = \vv_1^{\theta} = 0,
        \end{cases} \\
\textstyle \vv_1^{\phi_s} =
        \begin{cases}
        \frac{s^2}{[\bb_1^{\theta}]^3} & \!\!\!\!\colon \vv_1^{\theta} = \ww_1^{\theta} = 0 \neq \bb_1^{\theta} \\
        \vv_1^{\theta} & \!\!\!\!\colon \ww_1^{\theta} = \bb_1^{\theta} = 0 \neq \vv_1^{\theta} \\
        s & \!\!\!\!\colon \ww_1^{\theta} = \bb_1^{\theta} = \vv_1^{\theta} = 0,
        \end{cases}
\qquad \bigwedge_{i = 2}^{\width} [(\ww_i^{\phi_s} = \ww_i^{\theta}) \wedge (\bb_i^{\phi_s} = \bb_i^{\theta}) \wedge (\vv_i^{\phi_s} = \vv_i^{\theta})],
\end{gathered}
\end{equation}
and $\cc^{\phi_s} = \cc^{\theta}$. \Nobs that \cref{eqn:setting1:realization_risk,eqn:lemma:saddle_points_constant_realization:theta_s} demonstrate for all $s, x \in \R$ that
\begin{equation}\label{eqn:lemma:saddle_points_constant_realization:gamma}
\textstyle \phi_0 = \theta, \qquad \limsup_{\fs \to 0} \norm{\phi_{\fs} - \theta} = 0, \qqandqq \cN^{\phi_s}(x) = \frac{1}{2} \allowbreak + s^2 (s x \allowbreak + 1)^3.
\end{equation}
This and \cref{eqn:setting1:realization_risk} imply for all $s \in \R$ sufficiently close to $0$ that
\begin{equation}
\begin{split}
& \textstyle \loss(\phi_s) - \loss(\theta) = \int_0^1 (\frac{1}{2} + s^2 (s x + 1)^3 - x)^2 \, \d x - \int_0^1 (\frac{1}{2} - x)^2 \, \d x \\
& \textstyle = \int_0^1 s^2 (s x + 1)^3 (1 - 2x + s^2 (s x + 1)^3)\, \d x \\
& \textstyle = \int_0^1 s^2 (s x + 1)^3 (1 - 2x)\, \d x + \int_0^1 s^4 (s x + 1)^6 \, \d x \\
& \textstyle = \int_0^1 s^2 (1 - 2x)\, \d x + \int_0^1 s^2 (3 s x) (1 - 2x)\, \d x + \cO(s^4) =  - \frac{s^3}{2} + \cO(s^4).
\end{split}
\end{equation}
\Hence for all sufficiently small $s \in (0, 1]$ that $\loss(\phi_s) < \loss(\theta) = \loss(\phi_0)$. This, \cref{eqn:lemma:Hessian_of_Risk:w_c_c_c}, and \cref{eqn:lemma:saddle_points_constant_realization:gamma} prove that $\theta$ is a saddle point of $\loss$.
\end{cproof}

\cfclear
\begin{lemma}\label{lemma:saddle_point_1}
Assume \cref{setting1} and for every $z \in \{-1, 1\}$ let $\theta^z \allowbreak \in \R^{\dimension}$ satisfy for all $x \in \R$ that
\begin{equation}\label{eqn:lemma:saddle_point_1}
\textstyle (\nabla \loss)(\theta^z) = 0 \qqandqq \cN^{\theta^z}(x) = \frac{35 + 4 z \sqrt{7}}{70} + \frac{14}{5} \bigl[x - \frac{7 + z \sqrt{7}}{14}\bigr]^3.
\end{equation}
Then it holds for all $z \in \{-1, 1\}$ that $\cap_{i \in \cI^{\theta^z}} \{\pvt{i}{\theta^z}\} = \{\frac{7 + z \sqrt{7}}{14}\}$ and that $\theta^z$ is a saddle point of $\loss$.
\end{lemma}
\begin{cproof}{lemma:saddle_point_1}
\Nobs that \cref{eqn:lemma:saddle_point_1} and \cref{lemma:realizations} ensure for all $z \in \{-1, 1\}$, $x \in \R$ that
\begin{equation}\label{eqn:lemma:saddle_point_1:cN_new}
\textstyle \scrC^{\theta^z} = \frac{35 + 4 z \sqrt{7}}{70}, \qquad \Slope{\theta^z} = \frac{14}{5}, \qqandqq \bigcap_{i \in \cI^{\theta^z}} \{\pvt{i}{\theta^z}\} = \{\frac{7 + z \sqrt{7}}{14}\}.
\end{equation}
For every $z \in \{-1, 1\}$ let the perturbation $\phi^z \in C(\R, \R^{\dimension})$ of $\theta^z$ be given for all $s \in \R$, $i \in \Width$ by
\begin{equation}
\begin{gathered}\label{eqn:lemma:saddle_point_1:c_w_b_s_v_s}
\textstyle \cc^{\phi^z_s} = \cc^{\theta^z}, \qquad \ww_i^{\phi^z_s} = \ww_i^{\theta^z}, \qquad \bb_i^{\phi^z_s} =
\begin{cases}
\bb_i^{\theta^z} + \frac{z s}{4 \sqrt{7}} \ww_i^{\theta^z} & \colon i \in \cI^{\theta^z} \\
\bb_i^{\theta^z} & \colon i \notin \cI^{\theta^z},
\end{cases} \\
\textstyle \andqq \vv_i^{\phi^z_s} =
        \begin{cases}
        \vv_i^{\theta^z} + \frac{s [\ww_i^{\theta^z}]^{- 3} \indicator{(0, \infty)} (\abs{\vv_i^{\theta^z}})}{\sum_{\kappa \in \cI^{\theta^z}} \indicator{(0, \infty)} (\abs{\vv_{\kappa}^{\theta^z}})} & \colon i \in \cI^{\theta^z} \\
        \vv_i^{\theta^z} & \colon i \notin \cI^{\theta^z}.
        \end{cases}
\end{gathered}
\end{equation}
\Nobs that \cref{eqn:lemma:saddle_point_1:cN_new,eqn:lemma:saddle_point_1:c_w_b_s_v_s} prove for all $z \in \{-1, 1\}$, $s, x \in \R$ that
\begin{equation}\label{eqn:lemma:saddle_point_1:gamma}
\textstyle \phi^z_0 = \theta^z, \qquad \limsup_{\fs \to 0} \norm{\phi^z_{\fs} - \theta^z} = 0, \qqand
\end{equation}
\begin{equation}
\begin{split}
\textstyle \cN^{\phi^z_s}(x) & \textstyle = \cc^{\phi^z_s} + \sum_{i = 1}^{\width} \vv_i^{\phi^z_s} \bigl(\ww_i^{\phi^z_s} x + \bb_i^{\phi^z_s}\bigr)^3 \\
& \textstyle = \cc^{\theta^z} + \sum_{i \in \Width \backslash \cI^{\theta^z}} \vv_i^{\theta^z} [\bb_i^{\theta^z}]^3 + \sum_{i \in \cI^{\theta^z}} \vv_i^{\phi^z_s} [\ww_i^{\theta^z}]^3 \Bigl(x + \frac{\bb_i^{\theta^z}}{\ww_i^{\theta^z}} + \frac{z s}{4 \sqrt{7}}\Bigr)^3 \\
& \textstyle = \scrC^{\theta^z} + \bigl(x - \frac{7 + z \sqrt{7}}{14} + \frac{z s}{4 \sqrt{7}}\bigr)^3 \sum_{i \in \cI^{\theta^z}} \vv_i^{\phi^z_s} [\ww_i^{\theta^z}]^3 \\
& \textstyle = \scrC^{\theta^z} + \bigl(x - \frac{7 + z \sqrt{7}}{14} + \frac{z s}{4 \sqrt{7}}\bigr)^3 \sum_{i \in \cI^{\theta^z}} \Bigl(\vv_i^{\theta^z} [\ww_i^{\theta^z}]^3 + \frac{s \indicator{(0, \infty)} (\abs{\vv_i^{\theta^z}})}{\sum_{\kappa \in \cI^{\theta^z}} \indicator{(0, \infty)} (\abs{\vv_{\kappa}^{\theta^z}})}\Bigr) \\
& \textstyle = \scrC^{\theta^z} + \bigl(x - \frac{7 + z \sqrt{7}}{14} + \frac{z s}{4 \sqrt{7}}\bigr)^3 (\Slope{\theta^z} + s) = \frac{35 + 4 z \sqrt{7}}{70} + (\frac{14}{5} + s) \bigl(x - \frac{7 + z \sqrt{7}}{14} + \frac{z s}{4 \sqrt{7}}\bigr)^3.
\end{split}
\end{equation}
Combining this with \cref{eqn:setting1:realization_risk} assures for all $z \in \{- 1, 1\}$ and all $s \in \R$ sufficiently close to $0$ that
\begin{multline}
\textstyle \loss(\phi^z_s) - \loss(\theta^z) = \int_0^1 \bigl[\frac{35 + 4 z \sqrt{7}}{70} + (\frac{14}{5} + s) \bigl(x - \frac{7 + z \sqrt{7}}{14} + \frac{z s}{4 \sqrt{7}}\bigr)^3 - x\bigr]^2 \, \d x \\
\textstyle - \int_0^1 \bigl[\frac{35 + 4 z \sqrt{7}}{70} + \frac{14}{5} \bigl(x - \frac{7 + z \sqrt{7}}{14}\bigr)^3 - x\bigr]^2 \, \d x = - \frac{73 s^2}{68600} + \cO(s^3).
\end{multline}
\Hence for all $z \in \{-1, 1\}$ and all sufficiently small $s \in (0, 1]$ that $\loss(\phi^z_s) \allowbreak < \loss(\theta^z) \allowbreak = \loss(\phi^z_0)$. Combining this with \cref{eqn:lemma:Hessian_of_Risk:w_c_c_c} and \cref{eqn:lemma:saddle_point_1:gamma} implies for all $z \in \{-1, \allowbreak 1\}$ that $\theta^z$ is a saddle point of $\loss$.
\end{cproof}

\cfclear
\begin{lemma}\label{lemma:saddle_point_2_a}
Assume \cref{setting1} and let $\theta \allowbreak \in \allowbreak \R^{\dimension}$ satisfy for all $x \in \R$ that
\begin{equation}\label{eqn:lemma:saddle_point_2_a}
\textstyle (\nabla \loss)(\theta) = 0, \qquad \cN^{\theta}(x) = \frac{1}{2} + \frac{28}{5} [x - \frac{1}{2}]^3, \qqandqq \prod_{i = 1}^{\width} \vv_i^{\theta} \ww_i^{\theta} = 0.
\end{equation}
Then it holds that $\cap_{i \in \cI^{\theta}} \{\pvt{i}{\theta}\} = \{\frac{1}{2}\}$ and that $\theta$ is a saddle point of $\loss$.
\end{lemma}
\begin{cproof}{lemma:saddle_point_2_a}
Throughout this proof let $\fii \allowbreak \in \Width$ satisfy $\vv_{\fii}^{\theta} \ww_{\fii}^{\theta} = 0$ (cf.\ \cref{eqn:lemma:saddle_point_2_a}). \Nobs that \cref{eqn:lemma:saddle_point_2_a} and \cref{lemma:realizations} demonstrate that
\begin{equation}\label{eqn:lemma:saddle_point_2_a:theta}
\textstyle \scrC^{\theta} = \frac{1}{2}, \qquad \Slope{\theta} = \frac{28}{5}, \qqandqq \bigcap_{i \in \cI^{\theta}} \{\pvt{i}{\theta}\} = \{\frac{1}{2}\}.
\end{equation}
In the following we first prove, in the case $\ww_{\fii}^{\theta} = 0$, that $\theta$ is a saddle point of $\loss$. Let $t \in \R$ and the perturbation $\phi \in C(\R, \R^{\dimension})$ of $\theta$, in the case $\ww_{\fii}^{\theta} = 0$, be given for all $s \in \R$ by
\begin{equation}\label{eqn:lemma:saddle_point_2_a:gamma_s}
\begin{gathered}
\textstyle \cc^{\phi_s} = \cc^{\theta}, \qquad \bb_{\fii}^{\phi_s} = \bb_{\fii}^{\theta}, \qquad \Forall i \in \Width \backslash \{\fii\} \colon [\ww_i^{\phi_s} = \ww_i^{\theta}] \wedge [\bb_i^{\phi_s} = \bb_i^{\theta}] \wedge [\vv_i^{\phi_s} = \vv_i^{\theta}], \\
\textstyle t =
        \begin{cases}
        \bb_{\fii}^{\theta} \sqrt[3]{\vv_{\fii}^{\theta}} & \colon \vv_{\fii}^{\theta} \neq 0 \\
        \bb_{\fii}^{\theta} & \colon \vv_{\fii}^{\theta} = 0,
        \end{cases} \quad
\textstyle \ww_{\fii}^{\phi_s} = 
        \begin{cases}
        \frac{s}{\sqrt[3]{\vv_{\fii}^{\theta}}} & \colon \vv_{\fii}^{\theta} \neq 0 \\
        s & \colon \vv_{\fii}^{\theta} = 0,
        \end{cases} \qandq
\textstyle \vv_{\fii}^{\phi_s} =
        \begin{cases}
        \vv_{\fii}^{\theta} & \colon \vv_{\fii}^{\theta} \neq 0 \\
        s^4 & \colon \vv_{\fii}^{\theta} = 0.
        \end{cases}
\end{gathered}
\end{equation}
\Nobs that \cref{eqn:lemma:saddle_point_2_a:gamma_s} ensures for all $s, x \in \R$ that
\begin{gather}\label{eqn:lemma:saddle_point_2_a:gamma_case_1}
\textstyle \phi_0 = \theta, \quad \limsup_{\fs \to 0} \norm{\phi_{\fs} - \theta} = 0, \qandq
\vv_{\fii}^{\phi_s} (\ww_{\fii}^{\phi_s} x + \bb_{\fii}^{\phi_s})^3 =
        \begin{cases}
        (s x + t)^3 & \colon \vv_{\fii}^{\theta} \neq 0 \\
        s^4 (s x + t)^3 & \colon \vv_{\fii}^{\theta} = 0.
        \end{cases} 
\end{gather}
Combining this with \cref{eqn:setting1:realization_risk} and \cref{eqn:lemma:saddle_point_2_a:theta} shows for all $s, x \in \R$ that
\begin{equation}
\textstyle \cN^{\phi_s} (x) =
        \begin{cases}
        \frac{1}{2} + \frac{28}{5} (x - \frac{1}{2})^3 + (s x + t)^3 - t^3& \colon \vv_{\fii}^{\theta} \neq 0 \\
        \frac{1}{2} + \frac{28}{5} (x - \frac{1}{2})^3 + s^4 (s x + t)^3 & \colon \vv_{\fii}^{\theta} = 0.
        \end{cases}
\end{equation}
This and \cref{eqn:setting1:realization_risk} prove for all $s \in \R$ sufficiently close to $0$ that
\begin{equation}
\begin{split}
\textstyle \loss(\phi_s) - \loss(\theta) & \textstyle = 
        \begin{cases}
        - \frac{2 t^2}{25} s -\frac{s^3}{50} + \frac{s^6}{7} + t \, \cO(s^2) & \colon \vv_{\fii}^{\theta} \neq 0 \\
        - \frac{2 t^2}{25} s^5 - \frac{2 t}{25} s^6 - \frac{s^7}{50} + \frac{s^{14}}{7} + t \, \cO(s^8) & \colon \vv_{\fii}^{\theta} = 0
        \end{cases} \\
& \textstyle =
\begin{cases}
        - \frac{2 t^2}{25} s + \cO(s^2) & \colon \vv_{\fii}^{\theta} \neq 0, t \neq 0 \\
        - \frac{s^3}{50} + \frac{s^6}{7} & \colon \vv_{\fii}^{\theta} \neq 0, t = 0 \\
        - \frac{2 t^2}{25} s^5 + \cO(s^6) & \colon \vv_{\fii}^{\theta} = 0, t \neq 0 \\
        - \frac{s^7}{50} + \frac{s^{14}}{7} & \colon \vv_{\fii}^{\theta} = 0, t = 0.
        \end{cases}
\end{split}
\end{equation}
This ensures for all sufficiently small $s \in (0, 1]$ that $\loss(\phi_s) < \loss(\theta) = \loss(\phi_0)$. Combining this with \cref{eqn:lemma:Hessian_of_Risk:w_c_c_c} and \cref{eqn:lemma:saddle_point_2_a:gamma_case_1} proves, in the case $\ww_{\fii}^{\theta} = 0$, that $\theta$ is a saddle point of $\loss$.

In the next step we prove, in the case $[\ww_{\fii}^{\theta} \allowbreak \neq 0] \allowbreak \wedge [\vv_{\fii}^{\theta} \allowbreak = 0]$, that $\theta$ is a saddle point of $\loss$. Let the perturbation $\phi \in C(\R, \allowbreak \R^{\dimension})$ of $\theta$, in the case $[\ww_{\fii}^{\theta} \allowbreak \neq 0] \allowbreak \wedge [\vv_{\fii}^{\theta} \allowbreak = 0]$, be given for all $s \in \R$ by
\begin{equation}\label{eqn:lemma:saddle_point_2_a:gamma_s_2}
\begin{gathered}
\textstyle \ww_{\fii}^{\phi_s} = \ww_{\fii}^{\theta}, \qquad \bb_{\fii}^{\phi_s} = \bb_{\fii}^{\theta} + \ww_{\fii}^{\theta} s, \qquad \vv_{\fii}^{\phi_s} = \frac{s^2}{[\ww_{\fii}^{\theta}]^3}, \qquad \cc^{\phi_s} = \cc^{\theta}, \\
\textstyle \andqq \bigwedge_{i \in \Width \backslash \{\fii\}} \Bigl([\ww_i^{\phi_s} = \ww_i^{\theta}] \wedge [\bb_i^{\phi_s} = \bb_i^{\theta}] \wedge \Bigl[\vv_i^{\phi_s} = \vv_i^{\theta} - \frac{s^2 [\ww_i^{\theta} + \indicator{\{0\}}(\ww_i^{\theta})]^{-3} \indicator{(0, \infty)}(\vv_i^{\theta} \ww_i^{\theta})}{\sum_{\kappa = 1}^{\width} \indicator{(0, \infty)}(\vv_{\kappa}^{\theta} \ww_{\kappa}^{\theta})}\Bigr]\Bigr).
\end{gathered}
\end{equation}
\Nobs that \cref{eqn:setting1:realization_risk}, \cref{eqn:lemma:saddle_point_2_a:theta}, and \cref{eqn:lemma:saddle_point_2_a:gamma_s_2} ensure for all $s, x \in \R$ that
\begin{gather}\label{eqn:lemma:saddle_point_2_a:gamma_case_2}
\textstyle \phi_0 = \theta, \quad \limsup_{\fs \to 0} \norm{\phi_{\fs} - \theta} = 0, \qandq \cN^{\phi_s} (x) = \frac{1}{2} + (\frac{28}{5} - s^2) (x - \frac{1}{2})^3 + s^2 (x - \frac{1}{2} + s)^3.
\end{gather}
This and \cref{eqn:setting1:realization_risk} prove for all $s \in \R$ sufficiently close to $0$ that
\begin{equation}
\textstyle \loss(\phi_s) - \loss(\theta) = - \frac{2 s^4}{25} + \frac{9 s^6}{80} + \frac{5 s^8}{4} + s^{10} = - \frac{2 s^4}{25} + \cO(s^6).
\end{equation}
\Hence for all sufficiently small $s \in (0, 1]$ that $\loss(\phi_s) < \loss(\theta) = \loss(\phi_0)$. Combining this with \cref{eqn:lemma:Hessian_of_Risk:w_c_c_c} and \cref{eqn:lemma:saddle_point_2_a:gamma_case_2} implies, in the case $[\ww_{\fii}^{\theta} \neq 0] \wedge [\vv_{\fii}^{\theta} = 0]$, that $\theta$ is a saddle point of $\loss$.
\end{cproof}

\cfclear
\begin{lemma}\label{lemma:saddle_point_2_b}
Assume \cref{setting1} and let $\theta \allowbreak \in \allowbreak \R^{\dimension}$ satisfy for all $x \in \R$ that
\begin{equation}\label{eqn:lemma:saddle_point_2_b}
\textstyle (\nabla \loss)(\theta) = 0, \qquad \cN^{\theta}(x) = \frac{1}{2} + \frac{28}{5} [x - \frac{1}{2}]^3, \qqandqq \# \{i \in \Width \colon \vv^{\theta}_i \ww^{\theta}_i > 0\} > 1.
\end{equation}
Then it holds that $\cap_{i \in \cI^{\theta}} \{\pvt{i}{\theta}\} = \{\frac{1}{2}\}$ and that $\theta$ is a saddle point of $\loss$.
\end{lemma}
\begin{cproof}{lemma:saddle_point_2_b}
Throughout this proof for every $i \in \Width$ let $t_i \in \R$ satisfy $t_i = \vv_i^{\theta} [\ww_i^{\theta}]^3$ and let $\theta_{\width + i} \in \R$ be the $(\width + i)^{\text{th}}$ component of $\theta$. \Nobs that \cref{eqn:lemma:saddle_point_2_b} ensures that there exist $\fii, \fj \in \cI^{\theta}$ which satisfy
\begin{equation}\label{eqn:lemma:saddle_point_2_b:fii_fj}
\textstyle \fii \neq \fj \qqandqq \min\{\vv_{\fii}^{\theta} \ww_{\fii}^{\theta}, \vv_{\fj}^{\theta} \ww_{\fj}^{\theta}, t_{\fii}, t_{\fj}\} > 0.
\end{equation}
Moreover, \nobs that \cref{eqn:lemma:saddle_point_2_b} and \cref{lemma:realizations} prove that $\cap_{i \in \cI^{\theta}} \{\pvt{i}{\theta}\} = \{\frac{1}{2}\}$. \cref{lemma:Hessian_of_Risk} (applied with $d \with 3$, $\scra \with 0$, $\scrb \with 1$, $f \with \id_{\R}$ in the notation of \cref{lemma:Hessian_of_Risk}), the fact that
\begin{equation}
\textstyle \int_0^1 (x - \frac{1}{2})^4 \, \d x = \bigl[\frac{1}{5} (x - \frac{1}{2})^5\bigr]_{x = 0}^{x = 1} = \frac{1}{80},
\end{equation}
and the fact that
\begin{multline}
\textstyle \int_0^1 (x - \frac{1}{2}) [\cN^{\theta}(x) - x] \, \d x = \int_0^1 (x - \frac{1}{2}) [\frac{1}{2} + \frac{28}{5} [x - \frac{1}{2}]^3 - x] \, \d x \\
\textstyle = \int_0^1 \bigl[\frac{28}{5} (x - \frac{1}{2})^4 - (x - \frac{1}{2})^2\bigr]\, \d x = \bigl[\frac{28}{25} (x - \frac{1}{2})^5 - \frac{1}{3} (x - \frac{1}{2})^3\bigr]_{x = 0}^{x = 1} = - \frac{1}{75}
\end{multline}
\hence demonstrate for all $\ell \in \{\fii, \fj\}$ that
\begin{equation}\label{eqn:lemma:saddle_point_2_b:Hessian_bb}
\begin{split}
\textstyle (\frac{\partial^2}{\partial \theta_{\width + \ell}^2} \loss)(\theta) & \textstyle = 18 [\vv_{\ell}^{\theta}]^2 \int_0^1 (\ww_{\ell}^{\theta} x + \bb_{\ell}^{\theta})^4 \, \d x + 12 \vv_{\ell}^{\theta} \int_0^1 (\ww_{\ell}^{\theta} x + \bb_{\ell}^{\theta}) [\cN^{\theta}(x) - x] \, \d x \\
& \textstyle = 18 [\vv_{\ell}^{\theta}]^2 [\ww_{\ell}^{\theta}]^4 \int_0^1 (x - \frac{1}{2})^4 \, \d x + 12 \vv_{\ell}^{\theta} \ww_{\ell}^{\theta} \int_0^1 (x - \frac{1}{2}) [\cN^{\theta}(x) - x] \, \d x \\
& \textstyle = \frac{18}{80} [\vv_{\ell}^{\theta}]^2 [\ww_{\ell}^{\theta}]^4 - \frac{12}{75} \vv_{\ell}^{\theta} \ww_{\ell}^{\theta} = \frac{\vv_{\ell}^{\theta} \ww_{\ell}^{\theta}}{200} (45 t_{\ell} - 32)
\end{split}
\end{equation}
\begin{multline}
\textstyle \andqq (\frac{\partial^2}{\partial \theta_{\width + \fii} \partial \theta_{\width + \fj}} \loss)(\theta) = 18 \vv_{\fii}^{\theta} \vv_{\fj}^{\theta} \int_0^1 (\ww_{\fii}^{\theta} x + \bb_{\fii}^{\theta})^2 (\ww_{\fj}^{\theta} x + \bb_{\fj}^{\theta})^2 \, \d x \\
\textstyle = 18 \vv_{\fii}^{\theta} \vv_{\fj}^{\theta} [\ww_{\fii}^{\theta} \ww_{\fj}^{\theta}]^2 \int_0^1 (x - \frac{1}{2})^4 \, \d x = \frac{18}{80} \vv_{\fii}^{\theta} \vv_{\fj}^{\theta} [\ww_{\fii}^{\theta} \ww_{\fj}^{\theta}]^2 = \frac{9}{40} \vv_{\fii}^{\theta} \vv_{\fj}^{\theta} [\ww_{\fii}^{\theta} \ww_{\fj}^{\theta}]^2.
\end{multline}
\Hence that
\begin{equation}
\begin{split}
& \textstyle \det\!
        \begin{psmallmatrix}
        (\frac{\partial^2}{\partial \theta_{\width + \fii}^2} \loss)(\theta) & (\frac{\partial^2}{\partial \theta_{\width + \fii} \partial \theta_{\width + \fj}} \loss)(\theta) \\[1ex]
        (\frac{\partial^2}{\partial \theta_{\width + \fj} \partial \theta_{\width + \fii}} \loss)(\theta) & (\frac{\partial^2}{\partial \theta_{\width + \fj}^2} \loss)(\theta)
        \end{psmallmatrix}\!
\textstyle = \bigl[(\frac{\partial^2}{\partial \theta_{\width + \fii}^2} \loss)(\theta)\bigr] \bigl[(\frac{\partial^2}{\partial \theta_{\width + \fj}^2} \loss)(\theta)\bigr] - \bigl[(\frac{\partial^2}{\partial \theta_{\width + \fii} \partial \theta_{\width + \fj}} \loss)(\theta)\bigr]^2 \\
& \textstyle = \frac{[\vv_{\fii}^{\theta} \ww_{\fii}^{\theta}] [\vv_{\fj}^{\theta} \ww_{\fj}^{\theta}] (45 t_{\fii} - 32) (45 t_{\fj} - 32)}{40000} - \frac{81 [\vv_{\fii}^{\theta} \vv_{\fj}^{\theta}]^2 [\ww_{\fii}^{\theta} \ww_{\fj}^{\theta}]^4}{1600} = \frac{[\vv_{\fii}^{\theta} \ww_{\fii}^{\theta}] [\vv_{\fj}^{\theta} \ww_{\fj}^{\theta}] [(45 t_{\fii} - 32) (45 t_{\fj} - 32) - 45^2 t_{\fii} t_{\fj}]}{40000} \\
& \textstyle = - \frac{32 [\vv_{\fii}^{\theta} \ww_{\fii}^{\theta}] [\vv_{\fj}^{\theta} \ww_{\fj}^{\theta}] [45 t_{\fii} + 45 t_{\fj} - 32]}{40000} = \frac{[\vv_{\fii}^{\theta} \ww_{\fii}^{\theta}] [\vv_{\fj}^{\theta} \ww_{\fj}^{\theta}]}{1250} (32 - 45 t_{\fii} - 45 t_{\fj}).
\end{split}
\end{equation}
This, \cref{eqn:lemma:saddle_point_2_b:fii_fj}, \cref{eqn:lemma:saddle_point_2_b:Hessian_bb}, and the fact that at least one of the quantities $(45 t_{\fii} - 32)$, $(45 t_{\fj} - 32)$, and $(32 - 45 t_{\fii} - 45 t_{\fj})$ is strictly negative show that the $2 \times 2$ principal submatrix of the Hessian of $\loss$ at $\theta$ indexed by $\width + \fii$ and $\width + \fj$, given by
\begin{equation}\label{eqn:lemma:saddle_point_2_b:2x2submatrix}
        \begin{psmallmatrix}
        (\frac{\partial^2}{\partial \theta_{\width + \fii}^2} \loss)(\theta) & (\frac{\partial^2}{\partial \theta_{\width + \fii} \partial \theta_{\width + \fj}} \loss)(\theta) \\[1ex]
        (\frac{\partial^2}{\partial \theta_{\width + \fj} \partial \theta_{\width + \fii}} \loss)(\theta) & (\frac{\partial^2}{\partial \theta_{\width + \fj}^2} \loss)(\theta)
        \end{psmallmatrix}\!,
\end{equation}
has at least one strictly negative principal minor, i.e., either one of the diagonal elements or the determinant of \cref{eqn:lemma:saddle_point_2_b:2x2submatrix} is strictly negative. Combining this with \cref{eqn:lemma:Hessian_of_Risk:w_c_c_c} proves that the Hessian of $\loss$ at $\theta$ has at least one strictly negative and at least one strictly positive eigenvalue. \Hence that $\theta$ is a saddle point of $\loss$.
\end{cproof}

\cfclear
\begin{corollary}\label{cor:saddle_point_2}
Assume \cref{setting1} and let $\theta \allowbreak \in \allowbreak \R^{\dimension}$ satisfy for all $x \in \R$ that
\begin{equation}\label{eqn:cor:saddle_point_2}
\textstyle (\nabla \loss)(\theta) = 0, \qquad \cN^{\theta}(x) = \frac{1}{2} + \frac{28}{5} [x - \frac{1}{2}]^3, \qqandqq \# \{i \in \Width \colon \vv^{\theta}_i \ww^{\theta}_i \ge 0\} > 1.
\end{equation}
Then it holds that $\cap_{i \in \cI^{\theta}} \{\pvt{i}{\theta}\} = \{\frac{1}{2}\}$ and that $\theta$ is a saddle point of $\loss$.
\end{corollary}
\begin{cproof}{cor:saddle_point_2}
\Nobs that \cref{eqn:cor:saddle_point_2} ensures that
\begin{equation}
\textstyle \prod_{i = 1}^{\width} \vv^{\theta}_i \ww^{\theta}_i = 0 \qquad \text{or} \qquad \# \{i \in \Width \colon \vv^{\theta}_i \ww^{\theta}_i > 0\} > 1.
\end{equation}
\cref{lemma:saddle_point_2_a}, \cref{lemma:saddle_point_2_b}, and \cref{eqn:cor:saddle_point_2} \hence establish that $\cap_{i \in \cI^{\theta}} \{\pvt{i}{\theta}\} = \{\frac{1}{2}\}$ and that $\theta$ is a saddle point of $\loss$.
\end{cproof}

\subsection{Reduction of perturbations for PNNs with only \Active\ hidden neurons}
\label{subsec:reducing_perturbations}

\begin{lemma}\label{lemma:invariant_perturbation}
Assume \cref{setting:SNNs} and let $\theta \in \allowbreak \R^{\dimension}$ satisfy $\prod_{i = 1}^{\width} \ww^{\theta}_i \neq 0$. Then it holds for all $s \allowbreak \in \R^{\dimension}$ sufficiently close to $0$ that there exists $\scrs^s = \allowbreak (\scrs_1^s, \allowbreak \dots, \allowbreak \scrs_{2 \width + 1}^s) \allowbreak \in \R^{2 \width + 1}$ such that
\begin{equation}
\textstyle \Forall x \in \R \colon \cN^{\theta + s}(x) = \cc^{\theta} + \scrs_{2 \width + 1}^s + \sum_{i = 1}^{\width} (\vv^{\theta}_i [\ww^{\theta}_i]^d + \scrs_{\width + i}^s) (x + [\ww^{\theta}_i]^{-1} \bb^{\theta}_i + \scrs_i^s)^d
\end{equation}
and $\limsup_{\R^{\dimension} \ni \fs \to 0} \norm{\scrs^{\fs}} = 0$.
\end{lemma}
\begin{cproof}{lemma:invariant_perturbation}
\Nobs that the assumption that $\prod_{i = 1}^{\width} \ww^{\theta}_i \neq 0$ ensures that there exists $\eps \in (0, \infty)$ which satisfies
\begin{equation}
\textstyle \Forall s \in \{t \in \R^{\dimension} \colon \norm{t} \le \eps\} \colon \prod_{i = 1}^{\width} (\ww_i^{\theta} + \ww_i^{s}) \neq 0.
\end{equation}
For every $s \in \{t \in \R^{\dimension} \colon \norm{t} \le \eps\}$, $i \allowbreak \in \Width$ let $\scrs^s = \allowbreak (\scrs_1^s, \allowbreak \dots, \allowbreak \scrs_{2 \width + 1}^s) \allowbreak \in \R^{2 \width + 1}$ satisfy
\begin{equation}\label{eqn:lemma:invariant_perturbation:scrs}
\textstyle \scrs_i^s = \frac{\bb_i^{\theta} + \bb_i^{s}}{\ww_i^{\theta} + \ww_i^{s}} - \frac{\bb_i^{\theta}}{\ww_i^{\theta}}, \quad \scrs_{\width + i}^s = (\vv_i^{\theta} + \vv_i^{s}) [\ww_i^{\theta} + \ww_i^{s}]^{d} - \vv_i^{\theta} [\ww_i^{\theta}]^{d}, \qandq \scrs_{2 \width + 1}^s = \cc^s.
\end{equation}
\Nobs that \cref{eqn:lemma:invariant_perturbation:scrs} shows for all $s \in \R^{\dimension}$, $x \in \R$ with $\norm{s} \le \eps$ that
\begin{equation}
\begin{split}
\textstyle \cN^{\theta + s}(x) & \textstyle = \cc^{\theta} + \cc^{s} + \sum_{i = 1}^{\width} (\vv_i^{\theta} + \vv_i^{s}) ([\ww_i^{\theta} + \ww_i^{s}] x + \bb_i^{\theta} + \bb_i^{s})^d \\
& \textstyle = \cc^{\theta} + \cc^{s} + \sum_{i = 1}^{\width} (\vv_i^{\theta} + \vv_i^{s}) [\ww_i^{\theta} + \ww_i^{s}]^{d} \bigl(x + \frac{\bb_i^{\theta} + \bb_i^{s}}{\ww_i^{\theta} + \ww_i^{s}}\bigr)^d \\
& \textstyle = \cc^{\theta} + \scrs_{2 \width + 1}^s + \sum_{i = 1}^{\width} (\vv^{\theta}_i [\ww^{\theta}_i]^d + \scrs_{\width + i}^s) (x + [\ww^{\theta}_i]^{-1} \bb^{\theta}_i + \scrs_i^s)^d,
\end{split}
\end{equation}
\begin{equation}
\textstyle \scrs_{\width + i}^s = (\vv_i^{\theta} + \vv_i^{s}) [\ww_i^{\theta} + \ww_i^{s}]^{d} - \vv_i^{\theta} [\ww_i^{\theta}]^{d} = \vv_i^{\theta} ([\ww_i^{\theta} + \ww_i^{s}]^{d} - [\ww_i^{\theta}]^{d}) + \vv_i^{s} [\ww_i^{\theta} + \ww_i^{s}]^{d},
\end{equation}
\begin{equation}
\textstyle \andqq \scrs_i^s = \frac{\bb_i^{\theta} + \bb_i^{s}}{\ww_i^{\theta} + \ww_i^{s}} - \frac{\bb_i^{\theta}}{\ww_i^{\theta}} = \frac{(\bb_i^{\theta} + \bb_i^{s}) \ww_i^{\theta} - \bb_i^{\theta} (\ww_i^{\theta} + \ww_i^{s})}{(\ww_i^{\theta} + \ww_i^{s}) \ww_i^{\theta}} = \frac{\bb_i^{s} \ww_i^{\theta} - \bb_i^{\theta} \ww_i^{s}}{(\ww_i^{\theta} + \ww_i^{s}) \ww_i^{\theta}}.
\end{equation}
The assumption that $\prod_{i = 1}^{\width} \ww^{\theta}_i \neq 0$ and the fact that for all $i \in \Width$ it holds that
\begin{equation}
\textstyle \limsup_{\R^{\dimension} \ni s \to 0} (\abs{\ww_i^s} + \abs{\bb_i^s} + \abs{\vv_i^s} + \abs{\cc^s}) = 0
\end{equation}
\hence prove that $\limsup_{\R^{\dimension} \ni s \to 0} \norm{\scrs^s} = 0$.
\end{cproof}

\subsection{Non-global local minima via absence of descent directions}
\label{subsec:local_minima}

\cfclear
\begin{lemma}\label{pre_lemma:Cauchy_Schwarz_type}
Let $n \in \N$, $a_1, a_2, \dots, a_n \in (0, \infty)$, $\fa, \fb, b_1, b_2, \dots, b_n \in \R$ satisfy
\begin{equation}\label{eqn:pre_lemma:Cauchy_Schwarz_type}
\textstyle \fa = \sum_{i = 1}^n a_i, \qquad \fa \fb = \sum_{i = 1}^n a_i b_i, \qqandqq \fa \fb^2 = \sum_{i = 1}^n a_i [b_i]^2.
\end{equation}
Then $\sum_{i = 1}^n \abs{b_i - \fb} = 0$.
\end{lemma}
\begin{cproof}{pre_lemma:Cauchy_Schwarz_type}
\Nobs that \cref{eqn:pre_lemma:Cauchy_Schwarz_type} ensures that
\begin{multline}
\textstyle 0 = \fa^2 \fb^2 - \fa^2 \fb^2 = \fa (\fa \fb^2) - (\fa \fb)^2 \\
\textstyle = \bigl[\sum_{i = 1}^n a_i\bigr] \bigl[\sum_{i = 1}^n a_i [b_i]^2\bigr] - \bigl[\sum_{i = 1}^n a_i b_i\bigr]^2 = \frac{1}{2} \sum_{i = 1}^n \sum_{j = 1}^n a_i a_j (b_i - b_j)^2.
\end{multline}
This and the fact that $\min_{i \in \{1, 2, \dots, n\}} a_i > 0$ show that $\sum_{i = 1}^n \abs{b_i - b_1} = 0$. Combining this with \cref{eqn:pre_lemma:Cauchy_Schwarz_type} demonstrates that
\begin{equation}
\textstyle 0 = \fa \fb - \fa \fb = \bigl[\sum_{i = 1}^n a_i\bigr] \fb - \sum_{i = 1}^n a_i b_i = \bigl[\sum_{i = 1}^n a_i\bigr] (\fb - b_1).
\end{equation}
The fact that $\sum_{i = 1}^n \abs{b_i - b_1} = 0$ and the fact that $\min_{i \in \{1, 2, \dots, n\}} a_i > 0$ therefore establish that $\sum_{i = 1}^n \abs{b_i - \fb} = 0$.
\end{cproof}

\cfclear
\begin{lemma}\label{pre_lemma:equivalency_monomial_0}
Let $d, n \in \N$, $a_1, a_2, \dots, a_n \in (0, \infty)$, $\fa, \fb, \fc, b_1, b_2, \dots, b_n \in \R$ satisfy $d \ge 3$ and
\begin{equation}\label{eqn:pre_lemma:equivalency_monomial_0}
\textstyle \#\{x \in \R \colon \fc + \fa (x - \fb)^d = \sum_{i = 1}^n a_i (x - b_i)^d\} \ge d + 1.
\end{equation}
Then it holds that $\sum_{i = 1}^n \abs{b_i - \fb} = 0 = \fc$ and $\sum_{i = 1}^n a_i = \fa$.
\end{lemma}
\begin{cproof}{pre_lemma:equivalency_monomial_0}
\Nobs that \cref{eqn:pre_lemma:equivalency_monomial_0} ensures for all $x \in \R$ that
\begin{equation}\label{eqn:pre_lemma:equivalency_monomial_0:equiv}
\textstyle \fc + \fa (x - \fb)^d = \sum_{i = 1}^n a_i (x - b_i)^d.
\end{equation}
The fact that $d \ge 3$ hence demonstrates that
\begin{equation}\label{eqn:pre_lemma:equivalency_monomial_0:equiv_params}
\textstyle \fa = \sum_{i = 1}^n a_i, \qquad \fa \fb = \sum_{i = 1}^n a_i b_i, \qqandqq \fa \fb^2 = \sum_{i = 1}^n a_i [b_i]^2.
\end{equation}
Combining this with \cref{pre_lemma:Cauchy_Schwarz_type} (applied for every $i \in \{1, 2, \dots, n\}$ with $n \with n$, $\fa \with \fa$, $\fb \with \fb$, $a_i \with a_i$, $b_i \with b_i$ in the notation of \cref{pre_lemma:Cauchy_Schwarz_type}) shows that $\sum_{i = 1}^n \abs{b_i - \fb} = 0$. This, \cref{eqn:pre_lemma:equivalency_monomial_0:equiv}, and \cref{eqn:pre_lemma:equivalency_monomial_0:equiv_params} establish that $\sum_{i = 1}^n \abs{b_i - \fb} = 0 = \fc$ and $\sum_{i = 1}^n a_i = \fa$.
\end{cproof}

\cfclear
\begin{corollary}\label{pre_lemma:equivalency_monomial}
Let $d, n \in \N$, $\fa \in (0, \infty)$, $\fb, \fc, a_1, a_2, \dots, a_n, b_1, b_2, \dots, b_n \in \R$ satisfy
\begin{equation}\label{eqn:pre_lemma:equivalency_monomial}
\textstyle \#\{i \in \{1, 2, \dots, n\} \colon a_i \ge 0\} = 1, \qquad \Forall x \in \R \colon \fc + \fa (x - \fb)^d = \sum_{i = 1}^n a_i (x - b_i)^d,
\end{equation}
and $d \ge 3$. Then it holds that $\sum_{i = 1}^n \abs{b_i - \fb} = 0 = \fc$ and $\sum_{i = 1}^n a_i = \fa$.
\end{corollary}
\begin{cproof}{pre_lemma:equivalency_monomial}
Throughout this proof let $\pi \colon \{1, 2, \dots, n\} \to \{1, 2, \dots, n\}$ be a permutation which satisfies
\begin{equation}\label{eqn:pre_lemma:equivalency_monomial:pi_0}
\textstyle a_{\pi(n)} \ge 0 \qqandqq \Forall i \in \{1, 2, \dots, n - 1\} \colon a_{\pi(i)} < 0.
\end{equation}
\Nobs that \cref{eqn:pre_lemma:equivalency_monomial} and \cref{eqn:pre_lemma:equivalency_monomial:pi_0} ensure that for all $x \in \R$ it holds that
\begin{equation}\label{eqn:pre_lemma:equivalency_monomial:k_equation}
\textstyle - \fc + a_{\pi(n)} (x - b_{\pi(n)})^d = \fa (x - \fb)^d + \sum_{i = 1}^{n - 1} (- a_{\pi(i)}) (x - b_{\pi(i)})^d.
\end{equation}
\cref{pre_lemma:equivalency_monomial_0} (applied for every $i \in \{1, 2, \dots, n - 1\}$ with $d \with d$, $n \with n$, $a_i \with - a_{\pi(i)}$, $a_n \with \fa$, $b_i \with b_{\pi(i)}$, $b_n \with \fb$, $\fa \with a_{\pi(n)}$, $\fb \with b_{\pi(n)}$, $\fc \with -\fc$ in the notation of \cref{pre_lemma:equivalency_monomial_0}) and \cref{eqn:pre_lemma:equivalency_monomial:pi_0} hence establish that $\sum_{i = 1}^n \abs{b_i - \fb} = 0 = \fc$ and $\sum_{i = 1}^n a_i = \fa$.
\end{cproof}

\cfclear
\begin{lemma}\label{lemma:realization_28_5}
Assume \cref{setting1} and let $\theta \in \allowbreak \R^{\dimension}$ satisfy for all $x \in \R$ that
\begin{equation}\label{eqn:lemma:realization_28_5:cN_theta}
\textstyle \cN^{\theta}(x) = \frac{1}{2} + \frac{28}{5} [x - \frac{1}{2}]^3 \qqandqq \# \{i \in \Width \colon \vv^{\theta}_i \ww^{\theta}_i \ge 0\} = 1.
\end{equation}
Then it holds that
\begin{equation}\label{eqn:lemma:realization_28_5:result}
\textstyle \Slope{\theta} = \frac{28}{5}, \quad \cI^{\theta} = \Width, \quad \scrC^{\theta} = \cc^{\theta} = \frac{1}{2}, \qandq \sum_{i = 1}^{\width} \abs{\ww^{\theta}_i + 2 \bb^{\theta}_i} = 0 = (\nabla \loss)(\theta).
\end{equation}
\end{lemma}
\begin{cproof}{lemma:realization_28_5}
\Nobs that \cref{eqn:setting1:realization_risk,eqn:lemma:realization_28_5:cN_theta} ensure for all $x \in \R$ that
\begin{multline}\label{eqn:lemma:realization_28_5:cN}
\textstyle \frac{1}{2} + \frac{28}{5} [x - \frac{1}{2}]^3 = \cN^{\theta}(x) = \cc^{\theta} + \allowbreak \sum_{i = 1}^{\width} \vv_i^{\theta} (\ww_i^{\theta} x + \bb_i^{\theta})^3 \\
\textstyle = \cc^{\theta} + \allowbreak \sum_{i \in \Width \backslash \cI^{\theta}} \vv_i^{\theta} [\bb_i^{\theta}]^3 + \sum_{i \in \cI^{\theta}} \vv_i^{\theta} [\ww_i^{\theta}]^3 (x - \pvt{i}{\theta})^3 = \scrC^{\theta} + \sum_{i \in \cI^{\theta}} \vv_i^{\theta} [\ww_i^{\theta}]^3 (x - \pvt{i}{\theta})^3,
\end{multline}
\begin{equation}\label{eqn:lemma:realization_28_5:Slope}
\textstyle \Slope{\theta} = \sum_{i = 1}^{\width} \vv_i^{\theta} [\ww_i^{\theta}]^3 = \frac{28}{5}, \qqandqq \Exists j \in \Width \colon \vv_j^{\theta} \ww_j^{\theta} \ge 0 > \max_{i \in \Width \backslash \{j\}} (\vv_i^{\theta} \ww_i^{\theta}).
\end{equation}
\Hence that
\begin{equation}\label{eqn:lemma:realization_28_5:cI}
\textstyle \prod_{i = 1}^{\width} \vv_i^{\theta} \ww^{\theta}_i \neq 0, \qquad \cI^{\theta} = \Width, \qqandqq \scrC^{\theta} = \cc^{\theta}.
\end{equation}
This, \cref{eqn:lemma:realization_28_5:cN}, and \cref{pre_lemma:equivalency_monomial} (applied for every $i \in \Width$ with $d \with \allowbreak 3$, $n \with \allowbreak \width$, $\fa \with \allowbreak \nicefrac{28}{5}$, $\fb \with \allowbreak \nicefrac{1}{2}$, $\fc \with \allowbreak [\nicefrac{1}{2}] - \scrC^{\theta}$, $a_i \with \allowbreak \vv^{\theta}_i [\ww^{\theta}_i]^3$, $b_i \with \allowbreak \pvt{i}{\theta}$ in the notation of \cref{pre_lemma:equivalency_monomial}) show that
\begin{equation}
\textstyle \scrC^{\theta} = \cc^{\theta} = \frac{1}{2} \qqandqq \bigcap_{i \in \Width} \{\pvt{i}{\theta}\} = \{\frac{1}{2}\}.
\end{equation}
Combining this with \cref{lemma:derivatives_of_Risk} (applied with $d \with \allowbreak 3$, $\scra \with 0$, $\scrb \with 1$, $f \with \id_{\R}$ in the notation of \cref{lemma:derivatives_of_Risk}), \cref{eqn:lemma:realization_28_5:cN_theta}, \cref{eqn:lemma:realization_28_5:Slope}, \cref{eqn:lemma:realization_28_5:cI}, and the fact that for all $k \in \{0, 2, 3\}$ it holds that
\begin{equation}
\textstyle \int_0^1 x (x - \frac{1}{2})^2 \bigl[\frac{1}{2} + \frac{28}{5} [x - \frac{1}{2}]^3 - x\bigr] \, \d x = 0 = \int_0^1 (x - \frac{1}{2})^k \bigl[\frac{1}{2} + \frac{28}{5} [x - \frac{1}{2}]^3 - x\bigr] \, \d x
\end{equation}
establishes \cref{eqn:lemma:realization_28_5:result}.
\end{cproof}

\cfclear
\begin{lemma}[\resname{Loss of perturbed local minima}]\label{pre_lemma:loss_perturbed_loc_min}
We have for all $\fa, \fb, \fc, \fe \in \R$ that
\begin{equation}
\textstyle \int_{0}^{1} \bigl[\frac{1}{2} + \frac{28}{5}(x - \frac{1}{2})^3 - x\bigr]^2 \d x = \frac{1}{75}
\end{equation}
and
\begin{equation}
\begin{split}
& \textstyle \int_{0}^{1} \bigl[\frac{1}{2} + (\frac{28}{5} + \fa)(x - \frac{1}{2})^3 + \fb (x - \frac{1}{2})^2 + \fc (x - \frac{1}{2}) + \fe - x\bigr]^2 \d x - \frac{1}{75} \\
& \textstyle = \frac{\fa^2}{448} + \frac{\fb^2}{80} - \frac{2 \fc}{75} + \frac{\fa \fc}{40} + \frac{\fc^2}{12} + \frac{\fb \fe}{6} + \fe^2 = 7 \bigl(\frac{\fa}{56} + \frac{\fc}{10}\bigr)^2 + \bigl(\frac{\fb}{12} + \fe\bigr)^2 + \frac{\fc^2}{75} - \frac{2 \fc}{75} + \frac{\fb^2}{180}.
\end{split}
\end{equation}
\end{lemma}

We \nobs that the proof of \cref{pre_lemma:loss_perturbed_loc_min} is elementary and \hence omitted. It can be verified, e.g., using {\sc Wolfram Mathematica} \cite{Mathematica}.

\cfclear
\begin{lemma}\label{lemma:loc_min}
Assume \cref{setting1} and let $\theta \in \allowbreak \R^{\dimension}$ satisfy for all $x \in \R$ that
\begin{equation}\label{eqn:lemma:loc_min:cN_theta}
\textstyle \cN^{\theta}(x) = \frac{1}{2} + \frac{28}{5} [x - \frac{1}{2}]^3 \qqandqq \# \{i \in \Width \colon \vv^{\theta}_i \ww^{\theta}_i \ge 0\} = 1.
\end{equation}
Then
\begin{enumerate}[label=\textnormal{(\roman*)}]
\item\label{item:params:lemma:loc_min} it holds that $\cI^{\theta} = \Width$, $\scrC^{\theta} = \cc^{\theta} = \frac{1}{2}$, $\Slope{\theta} = \frac{28}{5}$, and $\cap_{i \in \Width} \{\pvt{i}{\theta}\} = \{\frac{1}{2}\}$,

\item\label{item:grad:lemma:loc_min} it holds that $(\nabla \loss)(\theta) = 0$, and

\item\label{item:loc_min:lemma:loc_min} it holds that $\theta$ is a non-global local minimum point of $\loss$.
\end{enumerate}
\end{lemma}
\begin{cproof}{lemma:loc_min}
Throughout this proof let $(\ft_1, \ft_2, \dots, \ft_{\width}) \in (- \infty, 0)^{\width - 1} \times [0, \infty)$, let $\pi \colon \allowbreak \Width \allowbreak \to \Width$ be a permutation, assume for all $i \in \Width$ that
\begin{equation}\label{eqn:lemma:loc_min:ft}
\textstyle \ft_i = \vv^{\theta}_{\pi(i)} \allowbreak [\ww^{\theta}_{\pi(i)}]^3,
\end{equation}
and for every $\fs = (\fs_1, \dots, \fs_{2 \width + 1}) \in \R^{2 \width + 1}$ let $f_{\fs} \colon \R \to \R$ satisfy for all $x \in \R$ that
\begin{equation}\label{eqn:lemma:loc_min:f}
\textstyle f_{\fs}(x) = \frac{1}{2} + \fs_{2 \width + 1} + \sum_{i = 1}^{\width} (\ft_i + \fs_{\width + i})[x - \frac{1}{2} + \fs_i]^3
\end{equation}
and let $a_{\fs}, b_{\fs}, c_{\fs}, e_{\fs} \in \R$ satisfy
\begin{equation}\label{eqn:lemma:loc_min:a_b_c_d}
\begin{gathered}
\textstyle \sum_{i = 1}^{\width} \fs_{\width + i} = a_{\fs}, \quad \sum_{i = 1}^{\width} \fs_i (\ft_i + \fs_{\width + i}) = \frac{b_{\fs}}{3}, \quad \sum_{i = 1}^{\width} [\fs_i]^2 (\ft_i + \fs_{\width + i}) = \frac{c_{\fs}}{3}, \\
\textstyle \andqq \fs_{2 \width + 1} + \sum_{i = 1}^{\width} [\fs_i]^3 (\ft_i + \fs_{\width + i}) = e_{\fs}.
\end{gathered}
\end{equation}
\Nobs that \cref{eqn:lemma:loc_min:cN_theta} and \cref{lemma:realization_28_5} establish \cref{item:params:lemma:loc_min,item:grad:lemma:loc_min}. Next \nobs that \cref{lemma:invariant_perturbation} (applied with $d \with 3$, $\scra \with 0$, $\scrb \with 1$, $f \with \id_{\R}$ in the notation of \cref{lemma:invariant_perturbation}), \cref{eqn:lemma:loc_min:cN_theta}, \cref{eqn:lemma:loc_min:ft}, and \cref{eqn:lemma:loc_min:f} show that there exists $\scrs^s = \allowbreak (\scrs^s_1, \allowbreak \dots, \allowbreak \scrs^s_{2 \width + 1}) \allowbreak \in \R^{2 \width + 1}$, $s \in \R^{\dimension}$, which satisfies for all $x \in \allowbreak \R$ and all $s \allowbreak \in \R^{\dimension}$ sufficiently close to $0$ that
\begin{equation}\label{eqn:lemma:loc_min:cN_theta_s}
\textstyle \cN^{\theta + s}(x) = \frac{1}{2} + \scrs^s_{2 \width + 1} + \sum_{i = 1}^{\width} (\ft_i + \scrs^s_{\width + i})[x - \frac{1}{2} + \scrs^s_i]^3 = f_{\scrs^s}(x)
\end{equation}
and $\limsup_{\norm{\fs} \to 0} \norm{\scrs^{\fs}} = 0$. In the next step \nobs that \cref{item:params:lemma:loc_min}, \cref{eqn:lemma:loc_min:cN_theta}, \cref{eqn:lemma:loc_min:f}, and \cref{eqn:lemma:loc_min:a_b_c_d} ensure for all $x \in \R$, $\fs = \allowbreak (\fs_1, \allowbreak \dots, \allowbreak \fs_{2 \width + 1}) \allowbreak \in \allowbreak \R^{2 \width + 1}$ that
\begin{equation}\label{t_h_positive}
\textstyle \ft_{\width} = \Slope{\theta} - \sum_{i = 1}^{\width - 1} \ft_i = \frac{28}{5} - \sum_{i = 1}^{\width - 1} \ft_i \ge \frac{28}{5} > 0
\end{equation}
and
\begin{equation}\label{cN_theta_s_explicit}
\begin{split}
& \textstyle f_{\fs}(x) = \frac{1}{2} + \fs_{2 \width + 1} + \sum_{i = 1}^{\width} (\ft_i + \fs_{\width + i})[x - \frac{1}{2} + \fs_i]^3 \\
& \textstyle = \frac{1}{2} + \frac{28}{5} [x - \frac{1}{2}]^3 + \bigl[\sum_{i = 1}^{\width} \fs_{\width + i}\bigr] [x - \frac{1}{2}]^3 + 3 \bigl[\sum_{i = 1}^{\width} \fs_i (\ft_i + \fs_{\width + i})\bigr] [x - \frac{1}{2}]^2 \\
& \textstyle \quad + 3 \bigl[\sum_{i = 1}^{\width} [\fs_i]^2 (\ft_i + \fs_{\width + i})\bigr] [x - \frac{1}{2}] + \fs_{2 \width + 1} + \sum_{i = 1}^{\width} [\fs_i]^3 (\ft_i + \fs_{\width + i}) \\
& \textstyle = \frac{1}{2} + \frac{28}{5} [x - \frac{1}{2}]^3 + a_{\fs} [x - \frac{1}{2}]^3 + b_{\fs} [x - \frac{1}{2}]^2 + c_{\fs} [x - \frac{1}{2}] + e_{\fs}.
\end{split}
\end{equation}
This and \cref{pre_lemma:loss_perturbed_loc_min} prove for all $\fs \in \R^{2 \width + 1}$ that
\begin{equation}\label{cR_theta_s}
\begin{split}
& \textstyle \int_0^1 [f_{\fs}(x) - x]^2 \, \d x - \int_0^1 [\cN^{\theta}(x) - x]^2 \, \d x \\
& \textstyle = \int_{0}^{1} (\frac{1}{2} + (\frac{28}{5} + a_{\fs})[x - \frac{1}{2}]^3 + b_{\fs} [x - \frac{1}{2}]^2 + c_{\fs} [x - \frac{1}{2}] + e_{\fs} - x)^2 \d x - \frac{1}{75} \\
& \textstyle = 7 \bigl(\frac{a_{\fs}}{56} + \frac{c_{\fs}}{10}\bigr)^2 + \bigl(\frac{b_{\fs}}{12} + e_{\fs}\bigr)^2 + \frac{[c_{\fs}]^2}{75} - \frac{2 c_{\fs}}{75} + \frac{[b_{\fs}]^2}{180}.
\end{split}
\end{equation}
Next \nobs that \cref{eqn:lemma:loc_min:a_b_c_d}, \cref{t_h_positive}, and the fact that $[\Forall i \in \{1, 2, \dots, \width - 1\} \colon \ft_i < 0]$ prove that for all $\fs = (\fs_1, \allowbreak \dots, \allowbreak \fs_{2 \width + 1}) \allowbreak \in \R^{2 \width + 1}$ sufficiently close to $0$ we have that
\begin{equation}\label{t_s_small}
\textstyle \min\{\ft_{\width} + \fs_{2 \width}, \frac{20}{5} + a_{\fs}\} > 0 \qqandqq \Forall i \in \{1, 2, \dots, \width - 1\} \colon [\ft_i + \fs_{\width + i} < 0].
\end{equation}
In the next step \nobs that \cref{eqn:lemma:loc_min:a_b_c_d} shows for all $\fs = (\fs_1, \dots, \fs_{2 \width + 1}) \in \R^{2 \width + 1}$ that
\begin{equation}
\textstyle \bigl[\frac{b_{\fs}}{3} - \sum_{i = 1}^{\width - 1} \fs_i (\ft_i + \fs_{\width + i})\bigr]^2 = [\fs_{\width}]^2 (\ft_{\width} + \fs_{2 \width})^2 = (\ft_{\width} + \fs_{2 \width}) \bigl[\frac{c_{\fs}}{3} - \sum_{i = 1}^{\width - 1} [\fs_i]^2 (\ft_i + \fs_{\width + i})\bigr].
\end{equation}
Combining this and Cauchy-Schwarz inequality with \cref{eqn:lemma:loc_min:a_b_c_d,t_h_positive,t_s_small} ensures for all $\fs = (\fs_1, \allowbreak \dots, \allowbreak \fs_{2 \width + 1}) \allowbreak \in \R^{2 \width + 1}$ sufficiently close to $0$ that
\begin{equation}
\begin{split}
& \textstyle [\frac{b_{\fs}}{3}]^2 - \frac{2 b_{\fs}}{3} \sum_{i = 1}^{\width - 1} \fs_i (\ft_i + \fs_{\width + i}) + \bigl[\sum_{i = 1}^{\width - 1} \fs_i (\ft_i + \fs_{\width + i})\bigr]^2 = \bigl[\frac{b_{\fs}}{3} - \sum_{i = 1}^{\width - 1} \fs_i (\ft_i + \fs_{\width + i})\bigr]^2 \\
& \textstyle = (\ft_{\width} + \fs_{2 \width}) \frac{c_{\fs}}{3} - (\ft_{\width} + \fs_{2 \width}) \sum_{i = 1}^{\width - 1} [\fs_i]^2 (\ft_i + \fs_{\width + i}) \\
& \textstyle = (\ft_{\width} + \fs_{2 \width}) \frac{c_{\fs}}{3} - \bigl(\frac{28}{5} + a_{\fs} - \sum_{i = 1}^{\width - 1} (\ft_i + \fs_{\width + i})\bigr) \sum_{i = 1}^{\width - 1} [\fs_i]^2 (\ft_i + \fs_{\width + i}) \\
& \textstyle = (\ft_{\width} + \fs_{2 \width}) \frac{c_{\fs}}{3} - (\frac{28}{5} + a_{\fs}) \bigl[\sum_{i = 1}^{\width - 1} [\fs_i]^2 (\ft_i + \fs_{\width + i})\bigr] \\
& \textstyle \quad + \bigl[\sum_{i = 1}^{\width - 1} (\underbrace{\ft_i + \fs_{\width + i}}_{< 0})\bigr] \bigl[\sum_{i = 1}^{\width - 1} [\fs_i]^2 (\underbrace{\ft_i + \fs_{\width + i}}_{< 0})\bigr] \\
& \textstyle \ge (\ft_{\width} + \fs_{2 \width}) \frac{c_{\fs}}{3} - (\frac{28}{5} + a_{\fs}) \bigl[\sum_{i = 1}^{\width - 1} [\fs_i]^2 (\ft_i + \fs_{\width + i})\bigr] + \bigl[\sum_{i = 1}^{\width - 1} \fs_i (\ft_i + \fs_{\width + i})\bigr]^2.
\end{split}
\end{equation}
\Hence for all $\fs = (\fs_1, \dots, \fs_{2 \width + 1}) \in \R^{2 \width + 1}$ sufficiently close to $0$ that
\begin{equation}
\textstyle [\frac{b_{\fs}}{3}]^2 - \frac{2 b_{\fs}}{3} \sum_{i = 1}^{\width - 1} \fs_i (\ft_i + \fs_{\width + i}) \ge (\ft_{\width} + \fs_{2 \width}) \frac{c_{\fs}}{3} - (\frac{28}{5} + a_{\fs}) \bigl[\sum_{i = 1}^{\width - 1} [\fs_i]^2 (\ft_i + \fs_{\width + i})\bigr].
\end{equation}
This and \cref{t_s_small} show for all $\fs = (\fs_1, \dots, \fs_{2 \width + 1}) \in \R^{2 \width + 1}$ sufficiently close to $0$ that
\begin{equation}
\begin{split}
& \textstyle [\frac{b_{\fs}}{3}]^2 \ge (\ft_{\width} + \fs_{2 \width}) \frac{c_{\fs}}{3} - (\frac{28}{5} + a_{\fs}) \bigl[\sum_{i = 1}^{\width - 1} (\ft_i + \fs_{\width + i}) \bigl([\fs_i]^2 - \frac{2 \fs_i b_{\fs}}{3 (a_{\fs} + [\nicefrac{28}{5}])}\bigr)\bigr] \\
& \textstyle = (\ft_{\width} + \fs_{2 \width}) \frac{c_{\fs}}{3} - (\underbrace{\textstyle \frac{28}{5} + a_{\fs}}_{> 0}) \bigl[\sum_{i = 1}^{\width - 1} (\underbrace{\ft_i + \fs_{\width + i}}_{< 0}) \underbrace{\bigl(\textstyle \fs_i - \frac{b_{\fs}}{3 (a_{\fs} + [\nicefrac{28}{5}])}\bigr)^2}_{\ge 0}\bigr] + [\frac{b_{\fs}}{3}]^2 \Bigl[\frac{\sum_{i = 1}^{\width - 1} (\ft_i + \fs_{\width + i})}{a_{\fs} + [\nicefrac{28}{5}]}\Bigr] \\
& \textstyle \ge (\ft_{\width} + \fs_{2 \width}) \frac{c_{\fs}}{3} + [\frac{b_{\fs}}{3}]^2 \Bigl[\frac{\sum_{i = 1}^{\width - 1} (\ft_i + \fs_{\width + i})}{a_{\fs} + [\nicefrac{28}{5}]}\Bigr].
\end{split}
\end{equation}
Combining this with \cref{eqn:lemma:loc_min:a_b_c_d,t_h_positive,t_s_small} demonstrates for all $\fs = (\fs_1, \dots, \fs_{2 \width + 1}) \in \R^{2 \width + 1}$ sufficiently close to $0$ that
\begin{multline}
\textstyle (\ft_{\width} + \fs_{2 \width}) \frac{c_{\fs}}{3} \le [\frac{b_{\fs}}{3}]^2 \Bigl[1 - \frac{\sum_{i = 1}^{\width - 1} (\ft_i + \fs_{\width + i})}{a_{\fs} + [\nicefrac{28}{5}]}\Bigr] = [\frac{b_{\fs}}{3}]^2 \Bigl[\frac{a_{\fs} + [\nicefrac{28}{5}] - \sum_{i = 1}^{\width - 1} (\ft_i + \fs_{\width + i})}{a_{\fs} + [\nicefrac{28}{5}]}\Bigr] \\
\textstyle = [\frac{b_{\fs}}{3}]^2 \bigl[\frac{\ft_{\width} + \fs_{2 \width}}{a_{\fs} + [\nicefrac{28}{5}]}\bigr] \le [\frac{b_{\fs}}{3}]^2 \bigl[\frac{\ft_{\width} + \fs_{2 \width}}{\nicefrac{8}{5}}\bigr].
\end{multline}
This and \cref{t_s_small} show for all $\fs \in \R^{2 \width + 1}$ sufficiently close to $0$ that
\begin{equation}
\textstyle \frac{[b_{\fs}]^2}{180} - \frac{2 c_{\fs}}{75} = \frac{1}{20} \cdot [\frac{b_{\fs}}{3}]^2 - \frac{2 c_{\fs}}{75} \ge \frac{1}{20} \cdot \frac{8}{5} \cdot \frac{c_{\fs}}{3} - \frac{2 c_{\fs}}{75} = 0.
\end{equation}
Combining this with \cref{lemma:inf_of_loss} (applied with $d \with 3$, $\scra \with 0$, $\scrb \with 1$, $f \with \id_{\R}$ in the notation of \cref{lemma:inf_of_loss}), \cref{pre_lemma:loss_perturbed_loc_min}, \cref{eqn:lemma:loc_min:cN_theta}, \cref{eqn:lemma:loc_min:cN_theta_s}, and \cref{cR_theta_s} demonstrates that there exists $\eps \in (0, \infty)$ which satisfies for all $s \in [-\eps, \eps]^{\dimension}$ that
\begin{equation}
\begin{split}
& \textstyle \loss(\theta + s) - \loss(\theta) = \int_0^1 [\cN^{\theta + s}(x) - x]^2 \, \d x - \int_0^1 [\cN^{\theta}(x) - x]^2 \, \d x \\
& \textstyle = \int_0^1 [f_{\scrs^s}(x) - x]^2 \, \d x - \int_0^1 [\cN^{\theta}(x) - x]^2 \, \d x \\
& \textstyle = 7 \bigl(\frac{a_{\scrs^s}}{56} + \frac{c_{\scrs^s}}{10}\bigr)^2 + \bigl(\frac{b_{\scrs^s}}{12} + e_{\scrs^s}\bigr)^2 + \frac{[c_{\scrs^s}]^2}{75} - \frac{2 c_{\scrs^s}}{75} + \frac{[b_{\scrs^s}]^2}{180} \ge - \frac{2 c_{\scrs^s}}{75} + \frac{[b_{\scrs^s}]^2}{180} \ge 0
\end{split}
\end{equation}
and
\begin{equation}
\textstyle \loss(\theta) = \frac{1}{75} > 0 = \inf_{\vartheta \in \R^{\dimension}} \loss(\vartheta).
\end{equation}
This and \cref{item:grad:lemma:loc_min} ensure that $\theta$ is a non-global local minimum point of $\loss$. This establishes \cref{item:loc_min:lemma:loc_min}.
\end{cproof}

\begin{remark}
A key feature of \cref{lemma:loc_min} is that for an arbitrary parameter $\theta \in \R^{\dimension}$, in order to be a non-global local minimum point of $\loss$, it is sufficient that exactly one of its hidden neurons has a non-negative \slope, and that its realization function is as in \cref{eqn:lemma:loc_min:cN_theta} above. These two assumptions alone force $\theta$ to have only \Active\ and \visible\ hidden neurons with all \pivot s coinciding with $\nicefrac{1}{2}$, and to be a critical point of $\loss$, in addition to being a non-global local minimum point of $\loss$.
\end{remark}

\subsection{From the canonical to the general setting via affine rescaling}
\label{subsec:affine_scaling}

\cfclear
\begin{lemma}[\resname{Affine invariance of the loss landscape}]\label{lemma:canonical_extension}
Assume \cref{setting:SNNs} with $f \in \cP(1)$, for every $\theta = (\theta_1, \allowbreak \dots, \allowbreak \theta_{\dimension}) \allowbreak \in \R^{\dimension}$ let $\theta^{\can} = (\theta_1^{\can}, \allowbreak \dots, \allowbreak \theta_{\dimension}^{\can}) \allowbreak \in \R^{\dimension}$ satisfy for all $i \in \Width$ that
\begin{equation}\label{eqn:lemma:canonical_extension:theta_can}
\begin{gathered}
\textstyle \theta_i^{\can} = (\scrb - \scra) \theta_i, \quad \theta_{\width + i}^{\can} = \theta_{\width + i} + \scra \theta_i, \quad \theta_{2 \width + i}^{\can} = \frac{\theta_{2 \width + i}}{(\scrb - \scra) f'(\scra)}, \qandq \theta_{\dimension}^{\can} = \frac{\theta_{\dimension} - f(\scra)}{(\scrb - \scra) f'(\scra)},
\end{gathered}
\end{equation}
and let $\losss \colon \R^{\dimension} \to \R$ satisfy for all $\theta \in \R^{\dimension}$ that $\losss(\theta) = \int_0^1 \abs{\cN^{\theta}(x) - x}^2 \, \d x$. Then
\begin{enumerate}[label=\textnormal{(\roman*)}]
\item
\label{item1:lemma:canonical_extension} it holds for all $\theta \in \R^{\dimension}$, $i \in \Width$ that the $i^{\text{th}}$ hidden neuron of $\theta^{\can}$ is
        \begin{enumerate}[label=\textnormal{(\alph*)}]
        \item \Active\ if and only if the $i^{\text{th}}$ hidden neuron of $\theta$ is \Active,
        
        \item \semiactive\ if and only if the $i^{\text{th}}$ hidden neuron of $\theta$ is \semiactive,
        
        \item \inactive\ if and only if the $i^{\text{th}}$ hidden neuron of $\theta$ is \inactive,
        
        \item \visible\ if and only if the $i^{\text{th}}$ hidden neuron of $\theta$ is \visible, and
        
        \item \invisible\ if and only if the $i^{\text{th}}$ hidden neuron of $\theta$ is \invisible,
        \end{enumerate}

\item
\label{item2:lemma:canonical_extension} it holds for all $\theta \in \R^{\dimension}$, $y \in \R$, $i \in \Width$ that
\begin{equation}
\begin{gathered}
\textstyle \cN^{\theta^{\can}}(\frac{y - \scra}{\scrb - \scra}) = \frac{\cN^{\theta}(y) - f(\scra)}{(\scrb - \scra) f'(\scra)}, \qquad \losss(\theta^{\can}) = \frac{\loss(\theta)}{(\scrb - \scra)^3 \abs{f'(\scra)}^2}, \qquad \cI^{\theta^{\can}} = \cI^{\theta}, \\
\textstyle \scrC^{\theta^{\can}} = \frac{\scrC^{\theta} - f(\scra)}{(\scrb - \scra) f'(\scra)}, \qquad \Slope{\theta^{\can}} = \frac{(\scrb - \scra)^{d - 1} \Slope{\theta}}{f'(\scra)}, \qqandqq \pvt{i}{\theta^{\can}} = \frac{\pvt{i}{\theta} - \scra}{\scrb - \scra},
\end{gathered}
\end{equation}

\item
\label{item3:lemma:canonical_extension} it holds for all $\theta = (\theta_1, \dots, \theta_{\dimension}) \in \R^{\dimension}$, $i \in \Width$ that
\begin{equation}
\begin{gathered}
\textstyle \frac{(\frac{\partial}{\partial \theta_i} \loss)(\theta)}{(\scrb - \scra)^4 \abs{f'(\scra)}^2} = (\frac{\partial}{\partial \theta_i^{\can}} \losss)(\theta^{\can}) + \frac{\scra}{\scrb - \scra} (\frac{\partial}{\partial \theta_{\width + i}^{\can}} \losss)(\theta^{\can}), \quad \frac{(\frac{\partial}{\partial \theta_{\dimension}} \loss)(\theta)}{(\scrb - \scra)^2 f'(\scra)} = (\frac{\partial}{\partial \theta_{\dimension}^{\can}} \losss)(\theta^{\can}) \\
\textstyle \frac{(\frac{\partial}{\partial \theta_{\width + i}} \loss)(\theta)}{(\scrb - \scra)^3 \abs{f'(\scra)}^2} = (\frac{\partial}{\partial \theta_{\width + i}^{\can}} \losss)(\theta^{\can}), \qqandqq \frac{(\frac{\partial}{\partial \theta_{2 \width + i}} \loss)(\theta)}{(\scrb - \scra)^2 f'(\scra)} = (\frac{\partial}{\partial \theta_{2 \width + i}^{\can}} \losss)(\theta^{\can}),
\end{gathered}
\end{equation}
and
\item
\label{item4:lemma:canonical_extension} it holds for all $\theta \in \R^{\dimension}$ that
        \begin{enumerate}[label=\textnormal{(\alph*)}]
        \item $(\nabla \loss)(\theta) = 0$ if and only if $(\nabla \losss)(\theta^{\can}) = 0$,
        
        \item $\theta$ is a global minimum point of $\loss$ if and only if $\theta^{\can}$ is a global minimum point of $\losss$,
        
        \item $\theta$ is a non-global local minimum point of $\loss$ if and only if $\theta^{\can}$ is a non-global local minimum point of $\losss$, and
        
        \item $\theta$ is a saddle point of $\loss$ if and only if $\theta^{\can}$ is a saddle point of $\losss$.
        \end{enumerate}
\end{enumerate}
\end{lemma}
\begin{cproof}{lemma:canonical_extension}
\Nobs that \cref{eqn:lemma:canonical_extension:theta_can} establishes that $\R^{\dimension} \ni \xi \mapsto \xi^{\can} \in \R^{\dimension}$ is a homeomorphism. Moreover, \nobs that \cref{eqn:lemma:canonical_extension:theta_can} ensures for all $\theta \in \R^{\dimension}$, $i \in \Width$ that
$\ww_i^{\theta^{\can}} \neq 0$ if and only if $\ww_i^{\theta} \neq 0$, that $\ww_i^{\theta^{\can}} = 0 \neq \bb_i^{\theta^{\can}}$ if and only if $\ww_i^{\theta} = 0 \neq \bb_i^{\theta}$, that $\ww_i^{\theta^{\can}} = 0 = \bb_i^{\theta^{\can}}$ if and only if $\ww_i^{\theta} = 0 = \bb_i^{\theta}$, and that $\vv_i^{\theta^{\can}} = 0$ if and only if $\vv_i^{\theta} = 0$. This establishes \cref{item1:lemma:canonical_extension}. In the next step we \nobs that \cref{eqn:setting:SNNs:realization_risk}, \cref{eqn:lemma:canonical_extension:theta_can}, \cref{item1:lemma:canonical_extension}, the fact that for all $g \in \cP(1)$, $x, y \in \R$ it holds that $g'(y) \neq 0$ and $g(y) + (x - y) g'(y) = g(x)$, the fact that $\scra < \scrb$, and the change-of-variables formula show for all $\theta = (\theta_1, \allowbreak \dots, \allowbreak \theta_{\dimension}) \allowbreak \in \R^{\dimension}$, $y \in \R$ that
\begin{equation}
\begin{split}
\textstyle \cN^{\theta^{\can}}(\frac{y - \scra}{\scrb - \scra}) & \textstyle = \theta_{\dimension}^{\can} + \sum_{i = 1}^{\width} \theta_{2 \width + i}^{\can} (\theta_{\width + i}^{\can} + \theta_i^{\can} \frac{y - \scra}{\scrb - \scra})^d \\
& \textstyle = \frac{\theta_{\dimension} - f(\scra)}{(\scrb - \scra) f'(\scra)} + \sum_{i = 1}^{\width} \frac{\theta_{2 \width + i}}{(\scrb - \scra) f'(\scra)} (\theta_{\width + i} + \scra \theta_i + (\scrb - \scra) \theta_i \frac{y - \scra}{\scrb - \scra})^d \\
& \textstyle = \frac{1}{(\scrb - \scra) f'(\scra)} \bigl[\theta_{\dimension} - f(\scra) + \sum_{i = 1}^{\width} \theta_{2 \width + i} (\theta_{\width + i} + \theta_i y)^d\bigr] = \frac{\cN^{\theta}(y) - f(\scra)}{(\scrb - \scra) f'(\scra)},
\end{split}
\end{equation}
\begin{equation}
\begin{split}
\textstyle \losss(\theta^{\can}) & \textstyle = \int_0^1 \abs{\cN^{\theta^{\can}}(x) - x}^2 \, \d x = \int_0^1 \Abs{\frac{\cN^{\theta}(x (\scrb - \scra) + \scra) - f(\scra)}{(\scrb - \scra) f'(\scra)} - x}^2 \, \d x \\
& \textstyle = \frac{1}{\scrb - \scra} \int_{\scra}^{\scrb} \Abs{\frac{\cN^{\theta}(x) - f(\scra)}{(\scrb - \scra) f'(\scra)} - \frac{x - \scra}{\scrb - \scra}}^2 \, \d x = \frac{1}{\scrb - \scra} \int_{\scra}^{\scrb} \Abs{\frac{\cN^{\theta}(x) - f(\scra) - (x - \scra) f'(\scra)}{(\scrb - \scra) f'(\scra)}}^2 \, \d x \\
& \textstyle = \frac{1}{(\scrb - \scra)^3 \abs{f'(\scra)}^2} \int_{\scra}^{\scrb} \abs{\cN^{\theta}(x) - f(x)}^2 \, \d x = \frac{\loss(\theta)}{(\scrb - \scra)^3 \abs{f'(\scra)}^2},
\end{split}
\end{equation}
\begin{equation}
\textstyle \cI^{\theta^{\can}} = \cI^{\theta}, \qquad \Forall i \in \Width \backslash \cI^{\theta} \colon \bb_i^{\theta^{\can}} = \bb_i^{\theta},
\end{equation}
\begin{equation}
\textstyle \scrC^{\theta^{\can}} = \cc^{\theta^{\can}} + \sum_{i \in \Width \backslash \cI^{\theta^{\can}}} \vv_i^{\theta^{\can}} [\bb_i^{\theta^{\can}}]^d = \frac{\cc^{\theta} - f(\scra)}{(\scrb - \scra) f'(\scra)} + \sum_{i \in \Width \backslash \cI^{\theta}} \frac{\vv_i^{\theta} [\bb_i^{\theta}]^d}{(\scrb - \scra) f'(\scra)} = \frac{\scrC^{\theta} - f(\scra)}{(\scrb - \scra) f'(\scra)},
\end{equation}
\begin{equation}
\textstyle \Slope{\theta^{\can}} = \sum_{i = 1}^{\width} \vv_i^{\theta^{\can}} [\ww_i^{\theta^{\can}}]^d = \sum_{i = 1}^{\width} \frac{\vv_i^{\theta} [(\scrb - \scra) \ww_i^{\theta}]^d}{(\scrb - \scra) f'(\scra)} = \frac{(\scrb - \scra)^{d - 1}}{f'(\scra)} \sum_{i = 1}^{\width} \vv_i^{\theta} [\ww_i^{\theta}]^d = \frac{(\scrb - \scra)^{d - 1} \Slope{\theta}}{f'(\scra)},
\end{equation}
\begin{equation}
\textstyle \andqq  \Forall (i, j) \in \{(\fii, \fj) \in \Width^2 \colon \ww_{\fii}^{\theta} = 0 \neq \ww_{\fj}^{\theta}\} \colon [\pvt{i}{\theta^{\can}} = \infty = \pvt{i}{\theta}] \wedge \bigl[\pvt{j}{\theta^{\can}} = \frac{\pvt{j}{\theta} - \scra}{\scrb - \scra}\bigr].
\end{equation}
This establishes \cref{item2:lemma:canonical_extension}. Furthermore, \nobs that \cref{eqn:setting:SNNs:realization_risk}, \cref{eqn:lemma:canonical_extension:theta_can}, \cref{item2:lemma:canonical_extension}, the fact that $\scra < \scrb$, and the fact that $f'(\scra) \neq 0$ establish \cref{item3:lemma:canonical_extension,item4:lemma:canonical_extension}.
\end{cproof}

\subsection{Structural characterization of critical points}
\label{subsec:structural_characterization}

\cfclear
\begin{lemma}[\resname{Structural degeneracy for odd activations}]\label{lemma:constant_realization_gen}
Assume \cref{setting:SNNs} with $\nicefrac{(d + 1)}{2} \in \N$ and $f \in \cP(1)$ and let $\theta \in \allowbreak (\nabla \loss)^{-1}(\{0\})$. Then the following statements are equivalent:
\begin{enumerate}[label=\textnormal{(\roman*)}]
\item It holds that $\sum_{i = 1}^{\width} (\abs{\ww^{\theta}_i} + \abs{\vv^{\theta}_i \bb^{\theta}_i}) = 0 = f(\frac{\scrb + \scra}{2}) - \cc^{\theta}$.

\item It holds for all $x \in \R$ that $\cN^{\theta}(x) = \cN^{\theta}(0)$.

\item It holds for all $x \in \R$ that $\cN^{\theta}(x) = f(\frac{\scrb + \scra}{2})$.
\end{enumerate}
\end{lemma}

\Nobs that \cref{lemma:constant_realization_gen} is a direct consequence of \cref{lemma:constant_realization}, \cref{lemma:canonical_extension}, and the fact that for all $g \in \cP(1)$ it holds that $\Forall x, y \in \R \colon \frac{1}{2} (y - x) g'(x) + g(x) = g(\frac{y + x}{2})$.

\cfclear
\begin{proposition}[\resname{Global minima}]\label{prop:global_minima_main}
Assume \cref{setting:SNNs} with $d = 3$ and $f \in \cP(1)$ and let $\theta \allowbreak \in \R^{\dimension}$. Then the following statements are equivalent:
\begin{enumerate}[label=\textnormal{(\roman*)}]
\item It holds that $\theta$ is a global minimum point of $\loss$.

\item It holds that $\loss(\theta) = 0$.

\item It holds that $\# (\cup_{i \in \cI^{\theta}} \{\pvt{i}{\theta}\}) - 1 > 0 = (\nabla \loss)(\theta)$.

\item It holds that $\# (\cup_{i \in \cI^{\theta}, \vv_i^{\theta} \neq 0} \{\pvt{i}{\theta}\}) - 3 \ge 0 = (\nabla \loss)(\theta)$.
\end{enumerate}
\end{proposition}

\Nobs that \cref{prop:global_minima_main} is a direct consequence of \cref{lemma:distinctness_of_kinks} and \cref{lemma:canonical_extension}.

\begin{proposition}[\resname{Saddle points}]\label{prop:saddle_points}
Assume \cref{setting:SNNs} with $d = 3$ and $f \in \cP(1)$ and let $\theta \allowbreak \in \R^{\dimension}$. Then $\theta$ is a saddle point of $\loss$ if and only if $\sum_{i \in \Width \backslash \cI^{\theta}} \abs{\vv^{\theta}_i \bb^{\theta}_i} = 0$ and one of the following holds:
\begin{enumerate}[label=\textnormal{(\roman*)}]
\item
\label{item1:prop:saddle_points} It holds that $\sum_{i = 1}^{\width} (\abs{\ww^{\theta}_i} + \abs{\vv^{\theta}_i \bb^{\theta}_i}) = 0 = \cc^{\theta} - f(\frac{\scrb + \scra}{2})$.

\item
\label{item2:prop:saddle_points} It holds that $\scrC^{\theta} = f(\frac{\scrb + \scra}{2}) - \frac{2 (\scrb - \scra) f'(\scra)}{5 \sqrt{7}}$, $\Slope{\theta} = \frac{14 f'(\scra)}{5 (\scrb - \scra)^2}$, and $\cap_{i \in \cI^{\theta}} \{\pvt{i}{\theta}\} = \{\frac{\scrb + \scra}{2} - \frac{\scrb - \scra}{2 \sqrt{7}}\}$.

\item
\label{item3:prop:saddle_points} It holds that $\scrC^{\theta} = f(\frac{\scrb + \scra}{2}) + \frac{2 (\scrb - \scra) f'(\scra)}{5 \sqrt{7}}$, $\Slope{\theta} = \frac{14 f'(\scra)}{5 (\scrb - \scra)^2}$, and $\cap_{i \in \cI^{\theta}} \{\pvt{i}{\theta}\} = \{\frac{\scrb + \scra}{2} + \frac{\scrb - \scra}{2 \sqrt{7}}\}$.

\item
\label{item4:prop:saddle_points} It holds that
\begin{equation}
\begin{gathered}
\textstyle \scrC^{\theta} = f(\frac{\scrb + \scra}{2}), \qquad \Slope{\theta} = \frac{28 f'(\scra)}{5 (\scrb - \scra)^2}, \qquad \bigcap_{i \in \cI^{\theta}} \{\pvt{i}{\theta}\} = \{\frac{\scrb + \scra}{2}\}, \\
\textstyle \andqq \# \{i \in \Width \colon f'(\scra) \vv^{\theta}_i \ww^{\theta}_i \ge 0\} > 1.
\end{gathered}
\end{equation}
\end{enumerate}
\end{proposition}
\begin{cproof}{prop:saddle_points}
First \nobs that if $\sum_{i \in \Width \backslash \cI^{\theta}} \abs{\vv^{\theta}_i \bb^{\theta}_i} = 0$ and one of \cref{item1:prop:saddle_points,item2:prop:saddle_points,item3:prop:saddle_points,item4:prop:saddle_points} holds then \cref{lemma:realizations}, \cref{lemma:saddle_points_constant_realization}, \cref{lemma:saddle_point_1}, \cref{cor:saddle_point_2}, and \cref{lemma:canonical_extension} ensure that $\theta$ is a saddle point of $\loss$. In the next step, to prove the opposite direction of the implication, we assume that $\theta$ is a saddle point of $\loss$. \Nobs that the assumption that $\theta$ is a saddle point of polynomial $\loss$ proves that $(\nabla \loss)(\theta) = 0$. \cref{lemma:realizations}, \cref{cor:saddle_point_2}, \cref{lemma:loc_min}, \cref{lemma:canonical_extension}, \cref{prop:global_minima_main}, and the fact that $\# \{i \in \Width \colon f'(\scra) \vv^{\theta}_i \ww^{\theta}_i \ge 0\} \ge 1$ \hence show that $\sum_{i \in \Width \backslash \cI^{\theta}} \abs{\vv^{\theta}_i \bb^{\theta}_i} = 0$ and one of \cref{item1:prop:saddle_points,item2:prop:saddle_points,item3:prop:saddle_points,item4:prop:saddle_points} holds.
\end{cproof}

\begin{proposition}[\resname{Non-global local minima}]\label{prop:local_minima}
Assume \cref{setting:SNNs} with $d = 3$ and $f \in \cP(1)$ and let $\theta \in \allowbreak \R^{\dimension}$. Then the following statements are equivalent:
\begin{enumerate}[label=\textnormal{(\roman*)}]
\item
\label{item1:prop:local_minima} It holds that $\theta$ is a non-global local minimum point of $\loss$.

\item
\label{item2:prop:local_minima} It holds for all $x \in \R$ that
\begin{equation}
\textstyle \cN^{\theta}(x) = f(\frac{\scrb + \scra}{2}) + \frac{28 f'(\scra)}{5 (\scrb - \scra)^2} \bigl(x - \frac{\scrb + \scra}{2}\bigr)^3 \qqandqq \# \{i \in \Width \colon f'(\scra) \vv^{\theta}_i \ww^{\theta}_i \ge 0\} = 1.
\end{equation}

\item
\label{item3:prop:local_minima} It holds that
\begin{equation}
\begin{gathered}
\textstyle \cI^{\theta} = \Width, \qquad \scrC^{\theta} = \cc^{\theta} = f(\frac{\scrb + \scra}{2}), \qquad \Slope{\theta} = \frac{28 f'(\scra)}{5 (\scrb - \scra)^2}, \\
\textstyle \bigcap_{i \in \cI^{\theta}} \{\pvt{i}{\theta}\} = \{\frac{\scrb + \scra}{2}\}, \qqandqq \# \{i \in \Width \colon f'(\scra) \vv^{\theta}_i \ww^{\theta}_i \ge 0\} = 1.
\end{gathered}
\end{equation}
\end{enumerate}
\end{proposition}
\begin{cproof}{prop:local_minima}
First \nobs that \cref{eqn:setting:SNNs:realization_risk}, \cref{lemma:loc_min}, and \cref{lemma:canonical_extension} ensure that \ref{item3:prop:local_minima} $\Leftrightarrow$ \ref{item2:prop:local_minima} $\Rightarrow$ \ref{item1:prop:local_minima}. Next it suffices to prove that \ref{item1:prop:local_minima} $\Rightarrow$ \ref{item2:prop:local_minima}. In the following we assume that $\theta$ is a non-global local minimum point of $\loss$. \Nobs that the assumption that $\theta$ is a non-global local minimum point of the polynomial $\loss$ proves that $(\nabla \loss)(\theta) = 0$. Combining this with \cref{lemma:realizations}, \cref{lemma:canonical_extension}, \cref{prop:global_minima_main}, \cref{prop:saddle_points}, the fact that $f'(\scra) \neq 0$, and the fact that
\begin{equation}
\textstyle \Slope{\theta} = \frac{28 f'(\scra)}{5 (\scrb - \scra)^2} \qquad \Rightarrow \qquad f'(\scra) \sum_{i = 1}^{\width} \vv_i^{\theta} (\ww_i^{\theta})^3 = f'(\scra) \Slope{\theta} = \frac{28 [f'(\scra)]^2}{5 (\scrb - \scra)^2} > 0
\end{equation}
establishes $\# \{i \in \Width \colon f'(\scra) \vv^{\theta}_i \ww^{\theta}_i \ge 0\} \ge 1$ and \ref{item1:prop:local_minima} $\Rightarrow$ \ref{item2:prop:local_minima}.
\end{cproof}

\subsection{Complete classification of critical points}
\label{subsec:complete_classification}

\begin{theorem}[\textcolor{darkgreen}{Complete classification}]\label{thm:complete_classification_cubic}
Assume \cref{setting:SNNs} with $d = 3$ and $f \in \cP(1)$, let $\theta \in \allowbreak \R^{\dimension}$, and for every $j \in \{-1, 0, 1\}$ let $\scrp_j, \scrf_j \allowbreak \in \R$, $\eta \in \N_0$ satisfy
\begin{equation}\label{eqn:thm:complete_classification_cubic:scrp_scrf}
\textstyle \scrp_j = \frac{\scrb + \scra}{2} + \frac{j (\scrb - \scra)}{2 \sqrt{7}}, \quad \scrf_j = \sup_{x \in [\scra, \scrb]} \! \Abs{f(\frac{\scrb + \scra}{2}) + \frac{2 j (\scrb - \scra) f'(\scra)}{5 \sqrt{7}} + \frac{28 f'(\scra) (x - \scrp_j)^3}{5 (1 + \abs{j}) (\scrb - \scra)^2} - \cN^{\theta}(x)},
\end{equation}
and $\eta = \#\{i \in \Width \colon \allowbreak f'(\scra) \vv^{\theta}_i \ww^{\theta}_i \ge 0\}$. Then the following statements hold:
\begin{enumerate}[label=\textnormal{(\roman*)}]
\item
\label{item:loc_max:thm:complete_classification_cubic} $\loss$ is a polynomial with no local maximum points.

\item
\label{item:global_min:thm:complete_classification_cubic} The following statements are equivalent:
\begin{enumerate}[label=\textnormal{(\alph*)}]
\item $\theta$ is a global minimum point of $\loss$ with $\loss(\theta) = 0$.

\item It holds that $\# (\cup_{i \in \cI^{\theta}} \{\pvt{i}{\theta}\}) - 1 > 0 = (\nabla \loss)(\theta)$.

\item It holds that $\# (\cup_{i \in \cI^{\theta}, \vv_i^{\theta} \neq 0} \{\pvt{i}{\theta}\}) - 3 \ge 0 = (\nabla \loss)(\theta)$.
\end{enumerate}

\item
\label{item:loc_min:thm:complete_classification_cubic} $\theta$ is a non-global local minimum point of $\loss$ if and only if $\eta - 1 = 0 = \scrf_0$.

\item
\label{item:saddle:thm:complete_classification_cubic} $\theta$ is a saddle point of $\loss$ if and only if $\sum_{i = 1}^{\width} (\abs{\ww^{\theta}_i} + \abs{\vv^{\theta}_i \bb^{\theta}_i}) = 0 = f(\frac{\scrb + \scra}{2}) - \cc^{\theta}$ or
\begin{equation}
\prod_{j = -1}^1 \textstyle \bigl(\indicator{\{0\}}(j) \indicator{\{0\}}(\eta - 1) + \scrf_j + \sum_{i \in \cI^{\theta}} \abs{\pvt{i}{\theta} - \scrp_j}\bigr) = 0 = \! \displaystyle \sum_{i \in \Width \backslash \cI^{\theta}} \abs{\vv^{\theta}_i \bb^{\theta}_i}.
\end{equation}
\end{enumerate}
\end{theorem}

\Nobs that \cref{thm:complete_classification_cubic} is a direct consequence of \cref{lemma:Hessian_of_Risk}, \cref{prop:global_minima_main}, \cref{prop:saddle_points}, and \cref{prop:local_minima} above.

\begin{remark}
We \nobs that while \cref{thm:complete_classification_cubic} establishes the complete classification explicitly for the cubic activation ($d = 3$), the algebraic and geometric structures uncovered in the proof seem to provide the necessary bridge to completely characterize the optimization landscape for \PNNs\ with arbitrary odd and even degree monomial activations, a generalization we leave for future work.
\end{remark}

\subsubsection*{Acknowledgements}
We gratefully acknowledge Prof.\ Dr.\ Arnulf Jentzen for his fruitful suggestions and valuable discussions during the preparation of this manuscript. We also gratefully acknowledge the Cluster of Excellence EXC 2044/2-390685587, Mathematics Münster: Dynamics-Geometry-Structure funded by the Deutsche Forschungsgemeinschaft (DFG, German Research Foundation) and DFG – Project-ID 499552394 – SFB 1597. Some useful \LaTeX\ commands from \cite{Bennoargumentcommand} are used.

\subsubsection*{Use of large language models (LLMs)}

We acknowledge {\sc Google Gemini} for assistance with the stylistic and linguistic refinement of the introductory texts, based on the authors' original structural outline, and {\sc ChatGPT (GPT-5.5 and 5.6)} for verifying mathematical statements and results upon completion of the manuscript. All AI-generated suggestions were strictly monitored, verified, and edited by the authors. The authors remain fully responsible for the final content.

\newpage
\bibliographystyle{acm}
\bibliography{bibfile}

\end{document}

%% file: commands.tex
\usepackage{amsmath,amssymb,amsthm,geometry,nicefrac,enumerate,mathtools,comment,bbm,mathrsfs,mathrsfs,
color,verbatim,amsfonts,graphicx,bm,dsfont,listings}

\DeclareFontEncoding{LS1}{}{}
\DeclareFontSubstitution{LS1}{stix}{m}{n}
\DeclareMathAlphabet{\mathscr}{LS1}{stixscr}{m}{n}

\usepackage{textcomp}
\usepackage{float}
\usepackage{tikz}
\usetikzlibrary{matrix,chains,positioning,decorations.pathreplacing,arrows}
\usetikzlibrary{shapes,fadings,arrows.meta}
\tikzset{font={\fontsize{9pt}{12}\selectfont}}
\usepackage{adjustbox}

\usepackage{datetime}
\usepackage[inline]{enumitem}
\usepackage{colortbl}
\usepackage{footnote}
\usepackage[most]{tcolorbox}
\usepackage[table]{xcolor}
\usepackage[hypertexnames=false]{hyperref}
\hypersetup{
  breaklinks=true,
  colorlinks,
  linkcolor={blue!90!black},  
  citecolor={green},
  urlcolor={green!60!black}   
}
\usepackage[hypcap=false]{caption}
\usepackage[capitalise,sort,nameinlink]{cleveref}

\crefname{enumi}{item}{items}
\crefname{equation}{}{}
\crefname{subsection}{Subsection}{Subsections}
\Crefname{subsection}{Subsection}{Subsections}
\crefname{figure}{Figure}{Figures}
\Crefname{figure}{Figure}{Figures}

\usepackage{xparse}
\usepackage{etoolbox}
\usepackage{cite}
\usepackage[english]{babel}
\usepackage[utf8]{inputenc}
\usepackage{subfig}
\usepackage{thmtools}

\theoremstyle{plain}
\newtheorem{theorem}{Theorem}[section]
\newtheorem{lemma}[theorem]{Lemma}

\newtheorem{corollary}[theorem]{Corollary}
\newtheorem{proposition}[theorem]{Proposition}
\newtheorem{setting}[theorem]{Setting}
\newtheorem{definition}[theorem]{Definition}

\newtheorem{remark}[theorem]{Remark}

\numberwithin{equation}{section}
\numberwithin{figure}{section}
\usepackage{mleftright}
\usepackage{mathtools}

\newcommand{\R}{\mathbb{R}}
\newcommand{\N}{\mathbb{N}}

\newcommand{\Z}{\mathbb{Z}}

\newcommand{\norm}[1]{\cfadd{def:p-norm}\lVert #1 \rVert}

\newcommand{\normmm}[1]{{\left\vert\kern-0.25ex\left\vert\kern-0.25ex\left\vert #1
    \right\vert\kern-0.25ex\right\vert\kern-0.25ex\right\vert}}
\newcommand{\abs}[1]{\lvert #1\rvert}
\newcommand{\Abs}[1]{\big| #1 \big|}

\newcommand{\eps}{\varepsilon}

\newcommand{\cA}{\mathcal{A}}

\newcommand{\cI}{\mathcal{I}}

\newcommand{\cN}{\mathcal{N}}
\newcommand{\cO}{\mathcal{O}}
\newcommand{\cP}{\mathcal{P}}

\newcommand{\fa}{\mathfrak{a}}
\newcommand{\fb}{\mathfrak{b}}
\newcommand{\fc}{\mathfrak{c}}

\newcommand{\fe}{\mathfrak{e}}

\newcommand{\fii}{\mathfrak{i}}
\newcommand{\fj}{\mathfrak{j}}

\newcommand{\fq}{\mathfrak{q}}

\newcommand{\fs}{\mathfrak{s}}
\newcommand{\ft}{\mathfrak{t}}

\newcommand{\bfd}{\mathbf{d}}

\newcommand{\scrC}{\mathscr{C}}

\newcommand{\scra}{\mathscr{a}}
\newcommand{\scrb}{\mathscr{b}}

\newcommand{\scrd}{\mathscr{d}}

\newcommand{\scrf}{\mathscr{f}}

\newcommand{\scrp}{\mathscr{p}}

\newcommand{\scrs}{\mathscr{s}}

\definecolor{codegreen}{rgb}{0,0.6,0}
\definecolor{codegray}{rgb}{0.5,0.5,0.5}
\definecolor{codepurple}{rgb}{0.58,0,0.82}
\definecolor{backcolour}{rgb}{0.95,0.95,0.92}

\definecolor{darkgreen}{rgb}{0.0, 0.6, 0.0}
\definecolor{cyan(process)}{rgb}{0.0, 0.72, 0.92}
\definecolor{cinnamon}{rgb}{0.82, 0.41, 0.12}
\definecolor{darkred}{rgb}{0.55, 0.0, 0.0} 
\definecolor{persianblue}{rgb}{0.11, 0.22, 0.73}
\definecolor{resnamecolor}{rgb}{0.0, 0.5, 0.0}
\newcommand{\resname}[1]{\textcolor{resnamecolor}{#1}}
\newcommand{\with}{\curvearrowleft}

\DeclareMathOperator{\id}{id}
\DeclareMathOperator{\Span}{span}

\newcommand{\qandq}{\quad \text{and} \quad}
\newcommand{\andqq}{\text{and} \qquad}
\newcommand{\qqand}{\qquad \text{and}}
\newcommand{\qqandqq}{\qquad \text{and} \qquad}

\newcommand{\loss}{\mathcal{L}}
\newcommand{\losss}{\mathfrak{L}}

\newcommand{\width}{\mathfrak{h}}
\newcommand{\Width}{\mathfrak{H}}
\newcommand{\can}{\mathrm{can}}
\newcommand{\Slope}[1]{\mathscr{S}^{#1}}
\newcommand{\ww}{w}
\newcommand{\bb}{b}
\newcommand{\vv}{v}
\newcommand{\cc}{c}

\newcommand{\Active}{active}
\newcommand{\semiactive}{semi-active}
\newcommand{\inactive}{inactive}

\newcommand{\visible}{visible}
\newcommand{\invisible}{invisible}
\newcommand{\slope}{slope}
\newcommand{\tslope}{total slope}
\newcommand{\pivot}{pivot}
\newcommand{\pvt}[2]{\mathfrak{p}_{#1}^{#2}}

\renewcommand{\d}{\mathrm{d}} 
\renewcommand{\emptyset}{\varnothing}

\newcommand{\dimension}{\mathfrak{d}}

\newcommand{\indicator}[1]{\mathbbm{1}_{\smash{#1}}}

\newcommand{\Exists}{\exists\,}
\newcommand{\Forall}{\forall\,}

\NewDocumentEnvironment{cproof}{m}
{\begin{proof}[Proof of \cref{#1}]}%
{\noindent The proof of \cref{#1} is thus complete.
\end{proof}}

\NewDocumentEnvironment{cproof2}{m}
{\begin{proof}[Proof of \cref{#1}]}%
{\noindent This completes the proof of \cref{#1}.
\end{proof}}

\ExplSyntaxOn
\NewDocumentCommand{\enum}{ O{;} m o }
 {
  \my_enum:nnn { #1 } { #2 } { #3 }
 }

\seq_new:N \l__my_enum_seq
\tl_new:N \l__my_enum_item_tl
\prop_const_from_keyval:Nn \l__verbs
{
show = shows ,
imply = implies ,
demonstrate = demonstrates ,
prove = proves ,
establish = establishes ,
ensure = ensures ,
assure = assures
}
\int_new:N \l__number_of_args

\cs_new_protected:Nn \my_enum:nnn
 {
  \seq_set_split:Nnn \l__my_enum_seq { #1 } { #2 }
  \seq_remove_all:Nn \l__my_enum_seq {}
  \int_set_eq:NN \l__number_of_args { \seq_count:N \l__my_enum_seq }
  \seq_use:Nnnn \l__my_enum_seq { ~and~ } { ,~ } { ,~and~ }
  \IfNoValueTF{#3}{}{
    \space
    \int_compare:nNnTF{ \l__number_of_args } < {2}{ \prop_item:Nn \l__verbs {#3} }{ #3 }
  }
 }

\seq_new:N \g_cflist_loaded
\seq_new:N \g_cflist_pending

\NewDocumentCommand{\cfadd}{ m }
{
  \seq_if_in:NnF \g_cflist_loaded { #1 } {
    \seq_if_in:NnF \g_cflist_pending { #1 } {
      \seq_gput_right:Nn \g_cflist_pending { #1 }
    }
  }
}

\NewDocumentCommand{\cfconsiderloaded}{ m }{
  \seq_gput_right:Nn \g_cflist_loaded {#1}
}

\NewDocumentCommand{\cfremove}{ m }
{
  \seq_gremove_all:Nn \g_cflist_pending { #1 }
}

\NewDocumentCommand{\cfload}{ o }
{
  \seq_if_empty:NTF \g_cflist_pending {\unskip} {
    (cf.\ \cref{\seq_use:Nn \g_cflist_pending {,}})\IfValueTF{#1}{#1~}{\unskip}
    \seq_gconcat:NNN \g_cflist_loaded \g_cflist_loaded \g_cflist_pending
    \seq_gclear:N \g_cflist_pending
  }
}

\NewDocumentCommand{\cfclear} {} {
  \seq_gclear:N \g_cflist_loaded
  \seq_gclear:N \g_cflist_pending
}

\NewDocumentCommand{\cfout}{ o }
{
  \seq_if_empty:NTF \g_cflist_pending {\unskip} {
    (cf.\ \cref{\seq_use:Nn \g_cflist_pending {,}})\IfValueTF{#1}{#1~}{\unskip}
    \seq_gclear:N \g_cflist_pending
  }
}

\NewDocumentCommand{\ifnocf} { m } {
  \seq_if_empty:NT \g_cflist_pending { #1 }
}

\ExplSyntaxOff

\ExplSyntaxOn

\bool_new:N \g_noteobserve

\NewDocumentCommand{\nobs}{}{
  \bool_if:nTF { \g_noteobserve } {
    \bool_gset_false:N \g_noteobserve
    note~
  } {
    \bool_gset_true:N \g_noteobserve
    observe~
  }
}

\NewDocumentCommand{\Nobs}{}{
  \bool_if:nTF { \g_noteobserve } {
    \bool_gset_false:N \g_noteobserve
    Note~
  } {
    \bool_gset_true:N \g_noteobserve
    Observe~
  }
}

\ExplSyntaxOff

\ExplSyntaxOn

\bool_new:N \g_hencetherefore

\NewDocumentCommand{\hence}{}{
  \bool_if:nTF { \g_hencetherefore } {
    \bool_gset_false:N \g_hencetherefore
    hence~
  } {
    \bool_gset_true:N \g_hencetherefore
    therefore~
  }
}

\NewDocumentCommand{\Hence}{}{
  \bool_if:nTF { \g_hencetherefore } {
    \bool_gset_false:N \g_hencetherefore
    Hence,~we~obtain~
  } {
    \bool_gset_true:N \g_hencetherefore
    Therefore,~we~obtain~
  }
}

\ExplSyntaxOff

\NewDocumentEnvironment {athm} {m m} {%
\begin{#1}\label{#2}\global\def\loc{#2}%
}{%
\end{#1}%
}

\NewDocumentEnvironment {adef} {m} {%
\begin{definition}\label{#1}\global\def\loc{#1}%
}{%
\end{definition}%
}

\NewDocumentEnvironment{aproof} {} {%
\begin{proof}[Proof~of~\cref{\loc}]%
}{%
\end{proof}%
}

\ExplSyntaxOff

\lstdefinestyle{mystyle}{
    backgroundcolor=\color{backcolour},   
    commentstyle=\color{codegreen},
    keywordstyle=\color{magenta},
    numberstyle=\tiny\color{codegray},
    stringstyle=\color{codepurple},
    basicstyle=\ttfamily\footnotesize,
    breakatwhitespace=false,         
    breaklines=true,                 
    captionpos=b,                    
    keepspaces=true,                 
    numbers=left,                    
    numbersep=5pt,                  
    showspaces=false,                
    showstringspaces=false,
    showtabs=false,                  
    tabsize=2
}
\ExplSyntaxOn

\int_new:N \g_furthermore

\NewDocumentCommand{\Moreover}{ o o }{
  \IfValueT{#1}{
    \str_case:nn {#1} {
      {Furthermore} {\int_set:Nn {\g_furthermore} {0}}
      {Moreover} {\int_set:Nn {\g_furthermore} {1}}
      {In~addition} {\int_set:Nn {\g_furthermore} {2}}
      {note} {\bool_gset_true:N \g_noteobserve}
      {observe} {\bool_gset_false:N \g_noteobserve}
    }
    \IfValueT{#2}{
      \str_case:nn {#2} {
        {Furthermore} {\int_set:Nn {\g_furthermore} {0}}
        {Moreover} {\int_set:Nn {\g_furthermore} {1}}
        {In~addition} {\int_set:Nn {\g_furthermore} {2}}
        {note} {\bool_gset_true:N \g_noteobserve}
        {observe} {\bool_gset_false:N \g_noteobserve}
      }
    }
  }
  \int_case:nn { \int_mod:nn {\g_furthermore} {3} } {
    { 0 } { Furthermore,~\nobs that}
    { 1 } { Moreover,~\nobs that}
    { 2 } { In~addition,~\nobs that}
  }
  \int_incr:N \g_furthermore
  \IfValueF{#1}{~}
}

\ExplSyntaxOff

\makeatletter
\ExplSyntaxOn
\seq_new:N \g_abbrs
\prop_new:N \g_abbr_counts
\tl_new:N \l_abbr_count_tl

\ExplSyntaxOn

\NewDocumentCommand{\abbr}{m m O{#1} m m O{#4} m}{
	\expandafter\newcommand\csname#3\endcsname[1][]{
		\seq_if_in:NnTF \g_abbrs {#1} {
			\prop_get:NnN \g_abbr_counts {#1} \l_abbr_count_tl
			\prop_gput:Nnx \g_abbr_counts {#1} {\int_eval:n {\l_abbr_count_tl + 1}}
			\hyperref[#1]{#7}
		} {
			\seq_gput_left:Nn \g_abbrs {#1}
			\prop_gput:Nnn \g_abbr_counts {#1} {1}
			\expandafter\gdef\csname#1@def\endcsname{#2}
			\phantomsection\label{#1}
			\str_if_eq:nnTF{##1}{}{\emph{#2}}{##1}~(\hyperref[#1]{#7})
		}
	}
	\expandafter\newcommand\csname#6\endcsname[1][]{
		\seq_if_in:NnTF \g_abbrs {#1} {
			\prop_get:NnN \g_abbr_counts {#1} \l_abbr_count_tl
			\prop_gput:Nnx \g_abbr_counts {#1} {\int_eval:n {\l_abbr_count_tl + 1}}
			\hyperref[#1]{#4}
		} {
			\expandafter\gdef\csname#1@def\endcsname{#5}
			\seq_gput_left:Nn \g_abbrs {#1}
			\prop_gput:Nnn \g_abbr_counts {#1} {1}
			\phantomsection\label{#1}
			\str_if_eq:nnTF{##1}{}{\emph{#5}}{##1}~(\hyperref[#1]{#4})
		}
	}
}

\ExplSyntaxOff

\makeatother
\abbr{OP}{optimization problem}{OPs}{optimization problems}{OP}
\abbr{ANN}{artificial neural network}{ANNs}{artificial neural networks}{ANN}
\abbr{DNN}{deep artificial neural network}{DNNs}{deep artificial neural networks}{DNN}
\abbr{PNN}{polynomial neural network}{PNNs}{polynomial neural networks}{PNN}
\abbr{DL}{deep learning}{DLs}{deep learning}{DL}
\abbr{GD}{gradient descent}{GDs}{gradient descent}{GD}
\abbr{GF}{gradient flow}{GFs}{gradient flow}{GF}
\abbr{SGD}{stochastic gradient descent}{SGDs}{stochastic gradient descent}{SGD}
\abbr{AI}{artificial intelligence}{AIs}{artificial intelligence}{AI}
\abbr{LLM}{large language model}{LLMs}{large language models}{LLM}
\abbr{GeLU}{Gaussian error linear unit}{GeLUs}{Gaussian error linear unit}{GeLU}
\abbr{ReLUs}{rectified linear unit}{ReLU}{rectified linear unit}
\makeatother